\documentclass[11pt,a4paper]{article}

\usepackage[T1]{fontenc}
\usepackage[latin1]{inputenc}
\usepackage{graphicx}
\usepackage{epsfig}
\usepackage{fancyhdr,fancybox}
\usepackage{amsfonts}
\usepackage{amsmath}
\usepackage{amssymb}
\usepackage{amsthm}
\usepackage{verbatim}
\usepackage{subfigure}
\usepackage[english]{babel}
\usepackage[latin1]{inputenc}
\usepackage{latexsym}
\usepackage{multicol}
\usepackage{booktabs}
\usepackage{color}
\usepackage{bm}
\usepackage{natbib}
\bibpunct{[}{]}{,}{n}{,}{,}

\newcommand{\cec}{{(\widehat{\sigma}_{ij}^{\,(1)} })_k^{\phantom{\big|}}}
\newcommand{\cau}{{(\widehat{\sigma}_{ij}^{\,(2)} })_k^{\phantom{\big|}}}
\newcommand{\cc}{{(\widehat{\sigma}_{ij}^{\,(3)}})_k^{\phantom{\big|}}}
\newcommand{\covM}{\bm \Sigma}

\newcommand{\B}{\mathcal{B}(\R^d)}

\newcommand{\Var}{{\rm Var}}

\newcommand{\cube}{K}
\renewcommand{\b}[1]{\mbox{\boldmath$#1$}}

\newcommand{\bbM}{\mathbb{M}}

\newcommand{\bo}{\bm o}
\newcommand{\R}{{\mathbb R}}
\newcommand{\E}{{\mathbb E}}

\newcommand{\bbP}{\mathbb{P}}

\newcommand{\pom}{P_M^{\bm o}}

\renewcommand{\o}{\mathbf{o}}
\newcommand{\e}{\varepsilon}

\newcommand{\be}{\begin{equation}}
\newcommand{\ee}{\end{equation}}
\newcommand{\bea}{\begin{eqnarray}}
\newcommand{\eea}{\end{eqnarray}}
\newcommand{\beastar}{\begin{eqnarray*}}
\newcommand{\eeastar}{\end{eqnarray*}}
\newcommand{\bee}{\begin{equation*}}
\newcommand{\eee}{\end{equation*}}

\newcommand{\diam}{\|\Xi_0\|}

\newcommand{\halmos}{\unskip\nobreak\hfill\hskip2em\null\nobreak\hfill$\Box$%
\vskip3mm}

\newcommand{\calA}{{\mathcal{A}}}
\newcommand{\calB}{{\cal B}}

\newcommand{\calH}{{\cal H}}

\newcommand{\ind}{{1\hspace{-1mm}{\rm I}}}

\newcommand{\kommentar}[1]{}
\newcommand{\lonk}{\mathop{\longrightarrow}\limits_{k \to \infty}}
\newcommand{\Lonk}{\stackrel{\mathrm{D}}{\lonk}}
\newcommand{\PLonk}{\stackrel{\mathbb{P}}{\lonk}}
\newcommand{\PasLonk}{\stackrel{\mathrm{\mathbb{P}-a.s.}}{\lonk}}
\newcommand{\Lond}{\stackrel{\mathrm{D}}{\longrightarrow}}
\newcommand{\eqd}{\stackrel{\mathrm{D}}{=}}

\newcommand{\Su}{\mathop{{\sum}^{\neq}}}
\newcommand{\lona}{\mathop{\longrightarrow}\limits_{a \to \infty}}

\newcommand{\lonr}{\mathop{\longrightarrow}\limits_{r \to \infty}}

\newcommand{\timesop}{\mathop{\times}\limits}
\newcommand{\timesim}{\timesop_{i=1}^m}
\newcommand{\timesimz}{\timesop_{i=2}^m}

\def\acknowledgement{\par\addvspace{17pt}\small\rmfamily
\trivlist\if!\ackname!\item[]\else
\item[\hskip\labelsep
{\it\ackname}]\fi}

\newcommand\ackname{Acknowledgements\runinend}
\def\runinend{.}

\newlength{\MyNumberwidth}
\newlength{\MyHeaderWidth}

\newtheorem{Th}{Theorem}[section]

\newtheorem{Sa}{Satz}[section]
\newtheorem{Lemma}[Sa]{Lemma}
\newtheorem{Cor}[Th]{Corollary}

\newtheorem{Def}{Definition}[section]

\numberwithin{equation}{section} 

\begin{document}
\pagestyle{fancy}
\title{\vspace*{-1cm}
Non-parametric asymptotic statistics for the Palm mark distribution of $\beta$-mixing marked point processes
}

\author{Lothar Heinrich$^{1\ast}$ \and Sebastian L\"uck$^2$
 \and Volker Schmidt$^2$}

\renewcommand{\thefootnote}{\arabic{footnote}}

\date{\today}
\maketitle
{\footnotesize
$^1$Institute of
Mathematics, University of Augsburg, D-86135 Augsburg, Germany\\
$^2$ Institute of Stochastics, Ulm
University, D-89069 Ulm, Germany\\}
$^\ast$ {\footnotesize Corresponding author\\[1mm]
}



\begin{abstract}
\noindent We consider spatially homogeneous marked point patterns in an unboundedly expanding convex 
sampling window. Our main objective is to identify the distribution of the typical mark
by constructing an asymptotic $\chi^2$-goodness-of-fit test. The corresponding test statistic is 
based on a natural empirical version of the Palm mark distribution and a smoothed 
covariance estimator which turns out to be mean-square consistent. Our approach does not
require  independent marks and allows dependences between the mark field and the point 
pattern. Instead we impose a suitable $\beta$-mixing condition on the underlying stationary 
marked point process which can be checked for a number of Poisson-based models and, in particular, 
in the case of geostatistical marking. Our method needs a central limit theorem for $\beta$-mixing 
random fields which is proved by extending Bernstein's blocking technique to non-cubic index sets
and seems to be of interest in its own right.     
By large-scale model-based simulations the performance of our test is
studied in dependence of the model parameters which determine the range of spatial correlations.\\

\vspace{0.02cm}

\vspace{0.02cm}
\begin{tabbing}
\noindent {\it MSC 2000}$\,$: {\sc Primary 62 G 10, 60 G 55; Secondary 60 F 05, 62 G 20}
\end{tabbing}
\end{abstract}

\smallskip\noindent
\section{Introduction}\label{sec.int.int}

Marked point processes (MPPs) are versatile models for the statistical analysis of data recorded at 
irregularly scattered locations. The simplest marking scenario is independent marking, where marks are given by a sequence of independent 
and identically distributed random elements, which is also independent of the underlying point pattern of locations.
A more complex class of models considers a so-called geostatistical marking, where the marks are determined by the values of a random
field at the given locations. Although the random field usually exhibits intrinsic spatial correlations, it  
is assumed to be independent of the location point process (PP).
However, in many real datasets interactions between locations and marks occur. 
Moreover, many marked point patterns arising in models from stochastic geometry
such as edge centers  in (anisotropic) Voronoi-tessellations marked by orientation or PPs marked by nearest-neighbour distances do not fit the setting of geostatistical marking. Statistical tests for independence between marks and points are e.g. discussed in \cite{Guan06, Guan07a, Schlather04, Schoenberg04}.
A frequent approach to investigate dependences in marked point patterns is based on 
 mark variogram and mark covariance functions. Recently, asymptotic normality of empirical 
versions of these functions with applications to mark correlation analysis has been studied 
 in \cite{Guan04, Guan07b, HeMo12}.
The main goal of this paper is to investigate estimators of the Palm mark distribution 
$\pom$ in point patterns exhibiting correlations between different marks as well as between
marks and locations. The probability measure $\pom$ can be interpreted as the distribution of the typical mark which denotes the mark of a randomly chosen point of the pattern.
For any mark set $C$ we consider the scaled deviations $Z_k(C)=\sqrt{|W_k|}
\bigl(\,(\widehat{P}^{\bm o}_M)^{\phantom{o}}_k(C) - \pom(C)\,\bigr)$ as measure
of the distance between  
 $\pom$ and an empirical Palm mark distribution $(\widehat{P}^{\bm o}_M)^{\phantom{o}}_k\,$.
Under appropriate strong mixing conditions we are able to prove asymptotic normality of 
 the scaled deviation vector
$\mathbf{Z}_k = (Z_k(C_1),\ldots,Z_k(C_\ell))^T$  when the
observation window $W_k$ with volume $|W_k|$ grows unboundedly in all directions as $k \to \infty$. 
The proof relies on Bernstein's blocking method, see e.g. \cite{Bul07, Naha91},
which so far has been applied only to sequences of cubic or cubelike  windows $W_k$, see e.g.  \cite{Guan07b, Hei94}. By means of some convex-geometric arguments it turns out that the blocking method is indeed applicable to any increasing sequence of convex  observation windows $W_k$ with unboundedly growing inball radii.
In addition  we  discuss consistent estimators for the covariance matrix of the Gaussian limit of $\mathbf{Z}_k$. This enables us to construct asymptotic $\chi^2$-goodness-of-fit tests for the Palm mark distribution $\pom$.
By means of computer simulations we study the convergence of first and
second type errors of the tests for growing observation windows in relation
to the range of dependence of the MPP. In this way we demonstrate the
practicability 
of the tests in analysis of real data. 
A promising field of application of our testing methodology  could be the directional analysis 
of random surfaces. 
Based on our results one can e.g. consider Cox processes on the boundary of Boolean models,
marking them with the local outer normal direction and testing for a hypothetical
 directional distribution. This allows to identify the rose of directions of the surface
 process associated with the Boolean model and represents an alternative to a Monte-Carlo test for the rose of direction suggested in \cite{benes01}.
The occurring marked point patterns differ basically from the setting of 
independent and geostatistical marking, for which functional central limit theorems (CLTs) and 
corresponding tests have been derived in \cite{Heinrich08, Pawlas09}.
Our paper is organized as follows. Section~\ref{sub.bas.not} introduces basic notation and 
definitions. In Section~\ref{sec.results} we present our main results, which are proved
in Section~\ref{sec:proofs}. In Section~\ref{sec:examples} we briefly discuss some 
models satisfying the assumptions needed to prove our asymptotic results.
In the final Section~\ref{sec.sta.app} we study the performance of the
proposed tests  by large-scale simulations.

\section{Stationary marked point processes}\label{sub.bas.not}

An MPP $X_M=  \sum_{n \ge 1}\delta_{(X_n, M_n)}$ is a random locally finite counting measure acting on the Borel sets of $\R^d \times \mathbb{M}$ 
with atoms $(X_n,M_n)\,$, where the marks $M_n$ belong to some Polish mark 
space $\mathbb{M}$ endowed with the Borel $\sigma$-algebra $\mathcal{B}(\mathbb{M})$. 
Throughout we assume that $X_M$ is simple, i.e. all locations $X_n$ in
$\R^d$ have multiplicity $1$ regardless which mark they have.
Mathematically spoken, $X_M$ is a measurable mapping $X_M : \Omega \longrightarrow 
\mathsf{N}_{\mathbb{M}}$ from some probability space $(\Omega, \mathcal{A},\mathbb{P})$ 
into the set $\mathsf{N}_{\mathbb{M}}$ of  counting measures $\varphi(\cdot)$ on 
${\mathcal B}(\R^d\times {\mathbb M})$ satisfying $\varphi(B\times\mathbb{M})< \infty$ 
for all bounded $B\in \B\,$, where $\mathsf{N}_{\mathbb{M}}$ is 
endowed with the smallest  $\sigma$-algebra $\mathcal{N}_{\mathbb{M}}$ containing 
all sets of the form $\{\varphi \in {\mathsf N}_{{\mathbb M}}: \varphi(B \times C)=j\}$ 
for $j \ge 0$, bounded $B \in \B\,$, and $C \in {\mathcal B}({\mathbb M})\,$.  
In what follows we only consider  {\it stationary} MPPs,
which means that the distribution $P_{X_M}(\cdot) = \mathbb{P}(X_M \in
(\cdot))$ of  $X_M$ on $\mathcal{N}_{\mathbb{M}}$ is invariant 
under location shifts of the atoms, i.e.,
\bee 
X_M \eqd 
 \sum_{n\geq 1} \delta_{(X_n-x,M_n
  )}\quad\mbox{for all}\quad x \in \R^d\,.
\eee

Provided that $X_M$ is stationary and the {\it intensity} $\lambda =
\E X_M([0,1)^d \times \mathbb{M})$ is finite we have 
$\E X_M(B\times C) = |B|\,\E X_M([0,1)^d \times C)$ for all bounded
$B\in \mathcal{B}(\R^d)$ and $C\in \mathcal{B}(\mathbb{M})\,$, 
where $|\cdot|$ denotes $d$-dimensional Lebesgue measure.

\subsection{Palm mark distribution}

For a stationary MPP $X_M$ the probability measure
$P_M^{\o}$ on $\mathcal{B}(\mathbb{M})$ defined by
\begin{equation}
\label{def.mar.dis}
P_M^{\o}(C) = \frac{1}{\lambda}\,
\mathbb{E} X_M([0,1)^d \times C)\,,\quad
C\in\mathcal{B}(\mathbb{M})\,,
\end{equation}
is called the {\it Palm mark distribution} of $X_M$. It can
be interpreted as the conditional distribution of the mark  
of an atom of $X_M$ located at the origin  $\o\,$. 
A random element $M_0$ in $\mathbb{M}$ with distribution $P_M^{\o}$ is 
 called {\it typical mark} of $X_M$.

\begin{Def}
{\rm An increasing sequence  $\{W_k\}$ of convex and compact sets in $\R^d$ 
such that $\varrho(W_k)=\sup\{ r>0 : B(x,r) \subset W_k$ for some $x \in W_k\}
 \rightarrow  \infty$ as $k \to \infty$ is called a {\it convex averaging
sequence} (briefly CAS). Here $ B(x,r)$ denotes the closed ball (w.r.t. the 
Euclidean norm $\|\cdot\|$) with midpoint at $x\in\R^d$ and  radius 
$r \ge 0\,$.}
\end{Def}

Some results from convex geometry applied to  CAS $\{W_k\}$ yield the following  
inequalities
\be 
\label{ineq:HeinrichPawlas} 
\frac{1}{\varrho(W_k)} \le \frac{\calH_{d-1}(\partial W_k)}{|W_k|} \le \frac{d}{\varrho(W_k)}
\qquad\mbox{and}\qquad      
1 - \frac{|W_k \cap(W_k-x)|}{|W_k|} \le \frac{d\,\|x\|}{\varrho(W_k)}
\qquad\quad
\ee

\medskip\noindent
for $\|x\| \le \varrho(W_k)\,$, where $\calH_{d-1}(\partial W_k)$ is the surface content 
(i.e. $(d-1)$-dimensional Hausdorff measure) of the boundary 
$\partial W_k\,$, see 
\cite{Boehm04} and \cite{Heinrich08} for details.

If $X_M$ is ergodic (for a precise definition see \cite{Daley0307}, p. 194), the individual ergodic theorem applied to MPPs (see Theorem 12.2.IV and Corollary 12.2.V in
\cite{Daley0307}) provides the $\mathrm{\mathbb{P}-a.s.}$ limits
\be 
\label{erg.emp.mark}
\widehat{\lambda}_k = \frac{X_M(W_k \times {\mathbb M})}{|W_k|} 
\PasLonk \lambda \quad\mbox{and}\quad
(\widehat{P}^{\o}_M)^{\phantom{o}}_{k}(C) = \frac{X_M(W_k \times C)}{X_M(W_k \times {\mathbb M})} 
\PasLonk P^{\o}_M(C)
\ee
for any $C \in \mathcal{B}(\mathbb{M})$ and an arbitrary CAS $\{W_k\}\,$.

\subsection{Factorial moment measures and the covariance measure}\label{sec:mom_cum_meas}

For any integer $m \ge 1$, the  $m$th {\it factorial moment measure} $\alpha_{X_M}^{(m)}$ of
the MPP  $X_M$ is defined on $\mathcal{B}((\R^d\times {\mathbb M})^m)$  by 
\be
\label{def.alf.isk}
\alpha_{X_M}^{(m)}\bigl(\,\timesim ( B_i \times C_i ) \,\bigr) =
\E{\Su\limits_{n_1,\ldots,n_m 
\geq 1}}\prod\limits_{i=1}^m \bigl(\,
\ind_{B_i}(X_{n_i})\ind_{C_i}(M_{n_i})\,\bigr)\;, 
\ee
where the sum $\sum_{n_1,\ldots, n_m \geq 1}^{\neq}$ runs over all $m$-tuples of pairwise distinct indices $n_1,\ldots ,n_m\ge 1$ for bounded  $B_i \in \B$ and $C_i\in \mathcal{B}(\mathbb{M})\,,\,i=1,\ldots, m$. We also need 
the $m$th factorial moment measure $\alpha_X^{(m)}$ of the {\it unmarked} PP  $X(\cdot) = X_M((\cdot)\times \bbM) = \sum_{n \ge 1}\delta_{X_n}(\cdot)$ defined on  $\mathcal{B}((\R^d)^m)$ by
\bee
\alpha_X^{(m)}\left( \timesim B_i \right) = 
\alpha_{X_M}^{(m)}\bigl(\,\timesim ( B_i \times \mathbb{M} ) \,\bigr)
\quad\mbox{for\;\;bounded}\quad B_1,\ldots, B_m \in \B\,.
\eee
The stationarity of $X_M$ implies that $\alpha_{X}^{(m)}$ is invariant under diagonal shifts, 
which allows to define the $m$th {\it reduced factorial moment measure} $\alpha_{X, red}^{(m)}$ 
uniquely determined by the following desintegration formula
\be
\label{def.red.alf}
\alpha_X^{(m)}\left( \timesim B_i \right)=\lambda\int_{B_1} 
\alpha_{X,red}^{(m)}\bigl(\,\timesimz (B_i-x)\,\bigr)\,{\rm d}x\;.
\ee 

We need a condition of weak dependence between parts of $X$ defined over distant Borel sets which can be expressed by the (factorial) {\it covariance measure}
$\gamma_X^{(2)}$ on $\mathcal{B}((\R^d)^2)$ defined by 
\bee
\gamma_X^{(2)}\bigl(B_1\times B_2\bigr) = \alpha_X^{(2)}\bigl(B_1 \times B_2\bigr) - 
\lambda^2\, |B_1|\;|B_2|\,. 
\eee

The {\it reduced covariance measure} $\gamma^{(2)}_{X,red}: \B \to [-\infty,\infty]$ is in general a
signed measure defined by (\ref{def.red.alf}) with $\gamma_X^{(2)}$ instead of $\alpha_X^{(2)}$, 
which means that
 \bee
\gamma^{(2)}_{X,red}(B)=\alpha^{(2)}_{X,red}(B) - \lambda\,|B|\,\quad\mbox{for\;\;bounded}\quad
B\in\B\,.
 \eee
For more details on factorial moment measures and measures related with them 
we refer to Chapters 8 and 12 in \cite{Daley0307}.

\subsection{$\b{m}$-point Palm mark distribution}\label{sub.two.poi}

For fixed mark sets $C_1,\ldots,C_m \in \mathcal{B}(\mathbb{M})\,,\,m\ge 1\,$, 
the $m$th factorial moment measure $\alpha_{X_M}^{(m)}$  defined by 
(\ref{def.alf.isk}) can be regarded as a measure on the Borel sets 
$\mathcal{B}((\R^d)^m)$ depending on $C_1,\ldots,C_m$. This new measure is 
absolutely continuous w.r.t. the $m$th factorial moment measure $\alpha_X^{(m)}$.  
Thus, the Radon-Nikodym theorem (cf. \cite{Folland99}, p. 90) implies the existence 
of a density  $P_M^{x_1,\ldots, x_m}(C_1\times\cdots\times C_m)$, which is 
uniquely determined for $\alpha_X^{(m)}$-almost all $(x_1,\ldots,x_m)\in (\R^d)^m$, 
such that for any $B_1,\ldots, B_m \in \B$,
\be
\label{cond:RadonNikodym}
 \alpha_{X_M}^{(m)}\bigl(\,\timesim (B_i\times C_i)\,\bigr)
=\int_{\timesim B_i }P_M^{x_1,\ldots, x_m} \left(\timesim C_i\right)\,
\alpha_X^{(m)}({\rm d}(x_1,\ldots, x_m)).
\ee
Since the mark space $\mathbb{M}$ is Polish, this Radon-Nikodym density can be extended 
to a regular conditional distribution of the mark vector $(M_1, \ldots, M_m)$ given that 
the corresponding  atoms $X_1, \ldots, X_m$ are located at pairwise distinct
points $x_1, \ldots, x_m$, i.e., 

\smallskip
\bee 
P_M^{x_1,\ldots, x_m}( C ) = 
\mathbb{P}((M_1,\ldots, M_m)\in C \mid X_1=x_1, \ldots , X_m=x_m)\quad\mbox{for}\quad C \in\mathcal{B}(\bbM^m)\,.
\eee

\smallskip\noindent
This means that the mapping $(x_1,\ldots,x_m, C) \mapsto P_M^{x_1,\ldots, x_m}(C)$ is a {stochastic kernel} , i.e., 
$ P_M^{x_1,\ldots, x_m}(C)$ is  $\mathcal{B}((\R^{d})^m)$-measurable in $(x_1,\ldots, x_m) \in (\R^d)^m$ for fixed $C \in \mathcal{B}(\bbM^m)$ 
and a probability measure in $C\in\mathcal{B}(\bbM^m)$ for fixed $(x_1,\ldots, x_m) \in (\R^d)^m$. For details we refer to \cite{Kallenberg86}, p. 164.
The regular conditional distribution 
$P_M^{x_1,\ldots, x_m}(C)$ for $C \in \mathcal{B}(\mathbb{M}^m)$  is called the {\it $m$-point Palm
mark distribution} of $X_M\,$. This  stochastic kernel is only of interest for
$m$-tuples $(x_1,\ldots,x_m)$ of  pairwise distinct points $x_i\in \R^d\,,\,i=1,...,m$. 
In case of a stationary simple MPP $X_M$ it can be shown that

\smallskip
\bee 
P_M^{x_1,\ldots, x_m}( C )= P_M^{\o,x_2-x_1\ldots, x_m-x_1} ( C )\quad\mbox{for }\quad 
C \in \mathcal{B}(\mathbb{M}^m)\,,\,m \ge 1
\eee

\smallskip\noindent
and any  $x_1,\ldots, x_m \in\R^d$ with $x_i \ne x_j$ for $i\ne j\,$. In this 
way the Palm mark distribution defined in (\ref{def.mar.dis})
can be considered as one-point Palm mark distribution.

The following result is crucial to prove asymptotic properties of variances estimators of the empirical mark distribution. It generalizes an analogous result stated for unmarked PPs 
in \cite{Hei10} to MPPs by involving the notion $m$-point Palm mark distribution for $m=2,3,4\,$. The proof 
is just a slight extension of the one of Lemma 5 in \cite{Hei10} by using the relation (\ref{cond:RadonNikodym}) for $m=2,3,4\,$. The details are left to the reader.

\begin{Lemma} \label{Lemma:FormelMonster}
Let $X_M=\sum_{n\ge 1}\delta_{(X_n,M_n)}$ be an MPP satisfying $\E
\bigl(X_M(B\times {\mathbb M})\bigr)^4 < \infty$ for all bounded $B\in \mathcal{B}(\R^d)$, and let $f:\,\R^d\times\R^d\times \mathbb{M}^2 \mapsto \R$ be a Borel-measurable function such that the second moment of $\sum\nolimits_{p,q \geq 1}^{\neq} |\,f(X_p, X_q, M_p, M_q)\,|$ exists. Then, 
\be
\label{streuung}
\Var \,\Big(\,\Su\limits_{p,q \ge 1} f(X_p, X_q, M_p, M_q)\,\Big) 
\ee
\be
= \int\limits_{(\R^d)^2} \int\limits_{{{\mathbb{M}}^{\phantom{}}}^2}
f(x_1, x_2, u_1, u_2)\Bigl[f(x_1, x_2, u_1, u_2)+f(x_2, x_1, u_2,
  u_1)\Bigr] P_M^{x_1, x_2}\bigl({\rm
  d}(u_1,u_2)\bigr)\alpha_X^{(2)}\bigl({\rm d}(x_1,x_2)\bigr)\nonumber\ee 
\bea
+ \int\limits_{(\R^d)^3} \int\limits_{{{\mathbb{M}}^{\phantom{}}}^3} f(x_1, x_2, u_1, u_2) \Bigl[f(x_1, x_3, u_1, u_3)+f(x_3, x_1, u_3, u_1)
\qquad\qquad\qquad\qquad\qquad\qquad\qquad\qquad\qquad\qquad\qquad\qquad&& \nonumber \\
+\; f(x_2, x_3, u_2, u_3)+ f(x_3, x_2, u_3, u_2)\Bigr]P_M^{x_1, x_2,
  x_3}\bigl({\rm d}(u_1, u_2, u_3)\bigr)\alpha_X^{(3)}\bigl({\rm d}(x_1,x_2, x_3)\bigr)\;\;\qquad\qquad\qquad\qquad\qquad&& 
\nonumber\\ 
+ \int\limits_{(\R^d)^4}\int\limits_{{{\mathbb{M}}^{\phantom{}}}^4} f(x_1, x_2, u_1, u_2)f(x_3, x_4, u_3, u_4)
\Bigl[
P_M^{x_1, x_2, x_3, x_4}\bigl({\rm d}(u_1, u_2, u_3,
u_4)\bigr)\alpha_X^{(4)}\bigl({\rm d}(x_1,x_2, x_3, x_4)\bigr)\qquad\qquad\qquad\qquad\qquad\;\;&&\nonumber\\  
-\;P_M^{x_1, x_2}\bigl({\rm d}(u_1, u_2)\bigr)P_M^{x_3, x_4}\bigl({\rm
  d}(u_3, u_4)\bigr) \alpha_X^{(2)}\bigl({\rm d}(x_1,x_2)\bigr)
\alpha_X^{(2)}\bigl({\rm d}(x_3,x_4)\bigr) \Bigr].\qquad\qquad\qquad\qquad\qquad&&\nonumber
\eea

\end{Lemma}

\subsection{$\b{\beta}$-mixing coefficient and covariance inequality}\label{sub.bet.mix}

For any $B \in \B$, let $\calA_{X_M}(B)$  denote the
sub-$\sigma$-algebra of $\calA$  generated by the
restriction of the MPP $X_M$ 
to the set $B \times {\mathbb M}$.  For any $B, B^\prime \in \B$ a natural
 measure of dependence between $\calA_{X_M}(B)$ and $\calA_{X_M}(B^\prime)$ can be 
formulated in terms of the $\beta-${\em mixing} (or {\em absolute regularity}, respectively 
{\em weak Bernoulli}) {\em coefficient}   
\be
\label{def.bet.emm}
\beta\bigl(\calA_{X_M}(B),\calA_{X_M}(B^\prime)\bigr)=\;\frac{1}{2}
 \sup_{\{A_i\}, \{A'_j\}}\, \sum_{i,j}\;
\bigl|\;\bbP(A_i \cap A_j^{\prime})\, - \, \bbP(A_i)\,\bbP(A_j^{\prime})\;\bigr|\,,
\ee 
where the supremum is taken over all finite partitions
$\{A_i\}$ and $\{A_j^{\prime}\}$ of $\Omega$ such that
$A_i \in \calA_{X_M}(B)$ and $A_j^{\prime} \in \calA_{X_M}(B^\prime)$ for all $i,j\,$,
see e.g. \cite{Dou94}, \cite{Hei94} or \cite{Yosh76}. It should be noticed that the
supremum in (\ref{def.bet.emm}) does not change if the sets $A_i$ and
$A_j^{\prime}$ belong to semi-algebras generating  $\calA_{X_M}(B)$ and
$\calA_{X_M}(B^\prime)$, respectively.
To express the  degree of  dependence of the MPP $X_M$ for disjoint sets 
$\cube_a = [-a,a]^d$ and 
$\cube^c_{a+b}=\R^d \setminus \cube_{a+b}$, where $b \ge 0$, we consider non-increasing functions 
$\beta_{X_M}^*,\; \beta_{X_M}^{**}:[\frac{1}{2},\infty) \to [0,\infty)$ such that
\bea  
\beta\bigl(\calA_{X_M}(\cube_a),\calA_{X_M}(\cube^c_{a+b})\bigr) \le 
\left\{\begin{array}{ll}\beta_{X_M}^*(b) & \quad\mbox{for}\quad 
\frac{1}{2}\le a \le b\;,\\  
& \phantom{000}\label{bet.sta.one}\\
a^{d-1}\,\beta_{X_M}^{**}(b) & 
\quad\mbox{for}\quad \frac{1}{2} \le b \le a\;.
\end{array}\right.
\eea

A stationary MPP $X_M$ is called $\beta$-{\it mixing} or {\it absolutely 
regular}, respectively {\em weak Bernoulli} if 
both $\beta$-{\it mixing rates} $\beta_{X_M}^*(r)$ and $\beta_{X_M}^{**}(r)$ tend to $0$ 
as $r \to \infty$. By standard measure-theoretic approximation arguments it is easily seen
that any stationary $\beta$-mixing MPP $X_M$ is  mixing in the usual sense and therefore also 
ergodic, see Lemma 12.3.II and Proposition~12.3.III in \cite{Daley0307} Vol. II p.~206. 
In order to prove CLTs we need further conditions on the decay of the $\beta$-{\it mixing rates}
$\beta_{X_M}^*(r)$ and $\beta_{X_M}^{**}(r)$ on the right-hand side (rhs)  of (\ref{bet.sta.one}). 
For this we formulate

{\bf Condition $\b{\beta(\delta)}$:} There exists some $\delta > 0$ such that
$\E\bigl(X_M(\,[0,1]^d \times \mathbb{M}\,)\bigr)^{2+\delta} < \infty\,$,
\be 
\int_1^\infty r^{d-1}\bigl(\beta_{X_M}^*(r)\bigr)^{\delta/(2+\delta)}\,{\rm d}r < \infty     
\qquad\mbox{and}\qquad r^{2d-1}\beta_{X_M}^{**}(r) \lonr 0\;.
\nonumber
\ee





The following type of covariance bound in terms of the $\beta$-mixing coefficient 
 (\ref{def.bet.emm}) was first stated in \cite{Yosh76}.

\begin{Lemma}\label{lem.cov.bet} Let $Y$ and $Y^\prime$ denote the restrictions of the MPP $X_M$ 
to $B\times \mathbb{M}\,$ and $B^\prime\times \mathbb{M}$  for some $B, B^\prime \in \B\,$, respectively. 
Furthermore, let $\widetilde Y$ and $\widetilde Y^\prime$ be independent copies of $Y$ and $Y^\prime$,
respectively. Then, for any ${\mathcal N}_{\mathbb{M}}\otimes {\mathcal N}_{\mathbb{M}}$-measurable function 
$f:\mathsf{N}_\mathbb{M}\times\mathsf{N}_\mathbb{M} \to [0,\infty)$
and for any $\eta>0$
\begin{eqnarray}
\bigl|\E f(Y,Y^\prime)-\E f(\widetilde Y,\widetilde
Y^\prime)\bigr|&\le & 2\,
\beta({\calA_{X_M}(B)},{\calA_{X_M}}(B^\prime))^{\frac{\eta}{1+\eta}}\nonumber\\
&\times& \max\Bigl\{\bigl(\E
f^{1+\eta}(Y,Y^\prime)\bigr)^{\frac{1}{1+\eta}},\,\bigl(\E
f^{1+\eta}(\widetilde Y,\widetilde Y^\prime)
\bigr)^{\frac{1}{1+\eta}}\Bigr\}\,.\label{cov.ine.bet}
\end{eqnarray}
\end{Lemma}
If $f$ is bounded, then {\rm (\ref{cov.ine.bet})} remains valid for $\eta = \infty\,$. 
In the particular case $f(y,y') = f_1(y)\,f_2(y')$ and $\eta=\delta/2$ for 
$\delta > 0$, the Cauchy-Schwarz inequality applied to the expectations on the rhs
of {\rm (\ref{cov.ine.bet})} yields
\be  
\bigl|\,{\rm Cov}\bigl(\,f_1(Y),f_2(Y^\prime)\,\bigr)\,\bigr|\le 2
\,\|\,f_1(Y)\,\|_{2+\delta}\,\|\,f_2(Y')\,\|_{2+\delta}\,
\bigl(\,\beta({\calA_{X_M}(B)},{\calA_{X_M}}(B^\prime)\,)\,\bigr)^{\frac{\delta}{2+\delta}}\,, 
\label{cov.ine.bet2} 
\ee
where $\|\,Z\,\|_q = (\E|Z|^q)^{1/q}$ is the ${L}^q-$norm $(q \ge 1)$ of a 
random variable $Z\,$.

\section{Results}\label{sec.results}


\subsection{Central limit theorem}

We consider a sequence of set-indexed empirical processes
$\{\,Y_k(C)\,,\,C \in \mathcal{B}(\mathbb{M})\,\}$ defined by 
\be
\label{yps.kah.til}
 Y_k(C) = \frac{1}{\sqrt{|W_k|}}\sum_{n \ge 1} \ind_{W_k}(X_n) \bigl( \ind_C(M_n) -
P^{\o}_M(C) \bigr) = \sqrt{|W_k|}\,\widehat{\lambda}_k\,\bigl( (\widehat{P}^{\o}_M)^{\phantom{o}}_k(C)
- P_M^{\o}(C) \bigr)\,,  
\ee
 
where $\{W_k\}$ is a CAS of observation windows in $\R^d$. 
We will first state a multivariate CLT for the joint distribution of $Y_k(C_1),\ldots,Y_k(C_\ell)$. For this, let  \lq$\Lond$\rq$\;$  denote {\em convergence in distribution} and ${\mathcal N}_{\ell}(a,\mathbf{\Sigma})$ be an $\ell$-dimensional Gaussian
vector with expectation vector $a\in\R^\ell$ and covariance matrix
$\mathbf{\Sigma} = (\sigma_{ij})_{i,j=1}^\ell$.

\begin{Th}\label{the.asy.nor}
Let $X_M$ be a stationary MPP with $\lambda > 0$ satisfying Condition $\beta(\delta)$. Then
\be
\label{asy.equ.nor}
\mathbf{Y}_k=\bigl(Y_k(C_1),\ldots,Y_k(C_\ell)\bigr)^\top \Lonk \;{\mathcal N}_\ell(\o_{\ell},\mathbf{\Sigma})  \quad\mbox{for any}\quad C_1,\ldots,C_\ell\in\calB(\mathbb{M})\,,\quad 
\ee
where $\o_{\ell}=(0,\ldots,0)^\top$ and  the asymptotic covariance matrix $\mathbf{\Sigma} =
(\sigma_{ij})_{i,j=1}^\ell$ is  given by the limits
\be 
\sigma_{ij}=\lim_{k \to \infty} \E Y_k(C_i)Y_k(C_j). \label{cov.mat.tau}
\ee
\end{Th}
The above result can also be stated in terms of the empirical set-indexed process 
$\{Z_k(C),\, C \in \mathcal{B}(\mathbb{M})\}$, where
\bee 
 Z_k(C) = (\,\widehat{\lambda}_k\,)^{-1} Y_k(C)=\sqrt{|W_k|}\,
\bigl(\,(\widehat{P}^{\o}_M)^{\phantom{o}}_k(C) - P^{\o}_M(C)\,\bigr)\,.
 \eee 
In other words, as refinement of the ergodic theorem (\ref{erg.emp.mark}), 
we derive  asymptotic normality of a suitably scaled deviation of the 
ratio-unbiased empirical Palm mark probabilities
$(\widehat{P}^{\o}_M)^{\phantom{o}}_k(C)$ from  $P^{\o}_M(C)$ defined by
(\ref{def.mar.dis}) for any $C \in \mathcal{B}(\mathbb{M})\,$.  
Since Condition $\beta(\delta)$ ensures the ergodicity of $X_M\,$, the first limiting relation in (\ref{erg.emp.mark}) combined with Slutsky's lemma
yields the following result as a corollary of Theorem~\ref{the.asy.nor}. 
\begin{Cor} 
The conditions of Theorem{\rm~\ref{the.asy.nor}} imply the CLT
\bee \mathbf{Z}_k = (Z_k(C_1),\ldots,Z_k(C_\ell))^\top \Lonk \;
{\mathcal N}_\ell(\o_{\ell},\lambda^{-2}\,\mathbf{\Sigma})\,. 
\label{asy.asy.nor} 
\eee
\end{Cor}

\subsection{$\b{\beta}$-mixing and integrability conditions}\label{sec:beta_integrability}

In this subsection we give a condition in terms of  the mixing rate $\beta_{X_M}^*(r)$ which implies finite total variation of the reduced covariance measure $\gamma^{(2)}_{X,red}$  and a certain integrability condition (\ref{int.con.alf}) which expresses
weak dependence between any two marks located at far distant sites. Both of these conditions are needed to get  the  asymptotic unbiasedness resp. $L^2$-consistency  of some estimators for the asymptotic covariances (\ref{cov.mat.tau}). 

Note that the total variation measure $|\gamma^{(2)}_{X,red}|$ of $\gamma^{(2)}_{X,red}$ 
is defined as sum of the positive part $\gamma^{(2)+}_{X,red}$ and negative part 
$\gamma^{(2)-}_{X,red}$ of the Jordan decomposition of $\gamma^{(2)}_{X,red}$, 
i.e., 
\be
\gamma^{(2)}_{X,red}= \gamma^{(2)+}_{X,red}-\gamma^{(2)-}_{X,red}
\qquad\mbox{and}\qquad
|\gamma^{(2)}_{X,red}| = \gamma^{(2)+}_{X,red} + \gamma^{(2)-}_{X,red}\,,
\nonumber
\ee
where the positive measures $\gamma^{(2)+}_{X,red}$ and $\gamma^{(2)-}_{X,red}$ 
are mutually singular, see \cite{Folland99}, p. 87.

\begin{Lemma}\label{lem.int.gam} 
Let $X_M$ be a stationary MPP satisfying  
\be
\E\bigl(X_M(\,[0,1]^d \times \mathbb{M}\,)\bigr)^{2+\delta} < \infty
\quad\mbox{and}\quad 
\int_1^\infty r^{d-1}\bigl(\beta_{X_M}^*(r)\bigr)^{\delta/(2+\delta)}\,
{\rm d}r < \infty\quad\mbox{for some}\quad \delta > 0\,.
\nonumber
\ee 
Then $\gamma_{X,red}^{(2)}$ has finite total variation on $\R^d\,$, i.e.,  
\be 
\label{int.cum.mea}
|\gamma_{X, red}^{(2)}|(\R^d) < \infty\,.  
\ee
Furthermore, for any $C_1,C_2 \in \cal{B}(\mathbb M)\,$
\be 
\label{int.con.alf} 
\int_{\R^d} \Bigl|\,P_M^{\o,x}(C_1 \times C_2) - P_M^{\o}(C_1)\,P_M^{\o}(C_2)\,\Bigr|\,
\alpha^{(2)}_{X,red}({\rm d}x) < \infty\,.
\ee 
\end{Lemma}

\medskip\noindent
\subsection{Representation of the asymptotic covariance matrix}
\label{sec:covrep}

In Theorem~\ref{the.asy.nor} we stated conditions for
asymptotic normality of the random vector
$\mathbf{Y}_k$. Clearly, (\ref{def.mar.dis}) and (\ref{yps.kah.til}) 
immediately imply that $\E Y_k(C) = 0$ for any $C\in\calB(\mathbb{M})$.
The following theorem gives a representation formula for the asymptotic 
covariance matrix ${\mathbf \Sigma}\,$.

\begin{Th}
\label{the.rep.tau}
Let $X_M$ be a stationary MPP satisfying $(\ref{int.con.alf})$ and let
$\{W_k\}$ be a CAS. Then, the limits in $(\ref{cov.mat.tau})$ exist
and take the form

\noindent
\begin{eqnarray}
\label{dar.tau.til}
\sigma_{ij} &=& \lambda \bigl(\,P_M^{\o}(C_i \cap C_j) 
- P_M^{\o}(C_i)\,P_M^{\o}(C_j)\,\bigr) + \lambda\,
\int_{\R^d}\bigl(\,P_M^{\o,x}(C_i \times C_j) \\
&&\nonumber\\
&-& P_M^{\o,x}(C_i \times \mathbb{M})\,P_M^{\o}(C_j) - P_M^{\o,x}(C_j
\times \mathbb{M})\,P_M^{\o}(C_i) + P_M^{\o}(C_i)\,P_M^{\o}(C_j)\,\bigr) \,
\alpha^{(2)}_{X,red}({\rm d}x)\,.\nonumber
\end{eqnarray}

In particular, if $X_M$ is an independently MPP, then 

\be 
\label{tau.ind.mar} 
\sigma_{ij}=\lambda \bigl(\,P_M^{\o}( C_i \cap C_j )- 
P_M^{\o}(C_i)\,P_M^{\o}(C_j)\,\bigr). 
\ee
\end{Th}

\medskip\noindent
\subsection{Estimation of the asymptotic covariance matrix}\label{sec:covest}

In Section~\ref{sec.sta.app} we will exploit the normal convergence (\ref{asy.equ.nor}) for
statistical inference of the typical mark distribution.
More precisely, assuming that the asymptotic covariance matrix ${\mathbf \Sigma}$ 
is invertible, we consider asymptotic $\chi^2$-goodness-of-fit tests, which are based on the distributional limit 
\bee
{\bf Y}_k^\top \widehat{\covM}_k^{-1}{\bf Y}_k \Lonk \chi^2_\ell,
\eee
which is an immediate consequence of (\ref{asy.asy.nor}) and Slutsky's lemma, given that 
${\widehat{\mathbf \Sigma}}_k$ is a consistent estimator for ${\mathbf \Sigma}$.
Here we use the notation ${\mathbf Y}_k=\bigl( Y_k(C_1),\ldots, Y_k(C_\ell)\bigr)^\top$ 
(see (\ref{yps.kah.til}))
and the random variable $\chi^2_\ell$ is $\chi^2$-distributed with $\ell$
degrees of freedom. In the following we will discuss several estimators for 
${\mathbf \Sigma}$. Our first observation is that the simple plug-in estimator
$\widehat{\mathbf \Sigma}_k^{(0)}=\bigl(Y_k(C_i)Y_k(C_j)\bigr)_{i,j=1}^{\ell}$ for ${\mathbf \Sigma}$
 is useless, since the determinant of $\widehat{\mathbf \Sigma}_k^{(0)}$
vanishes. Instead of $\widehat{\mathbf \Sigma}_k^{(0)}$ we take  the
edge-corrected estimator  $\widehat{\mathbf \Sigma}_k^{(1)}=\bigl(\cec\bigr)_{i,j=1}^{\ell}$  with
\begin{eqnarray}
\label{est.tau.one}
\cec &=&  \frac{1}{|W_k|} \;\sum_{p \ge 1}
\ind_{W_k}(X_p) \bigl(\,\ind_{C_i \cap
C_j}(M_p)-P_M^{\o}(C_i)\,P_M^{\o}(C_j)\,\bigr)\\
&+& \Su\limits_{p,q \ge 1}\;
\frac{\ind_{W_k}(X_p)\ind_{W_k}(X_q)\bigl(\,\ind_{C_i}(M_p)-P_M^{\o}(C_i)\,\bigr)
\bigl(\ind_{C_j}(M_q)-P_M^{\o}(C_j)\bigr)}{|(W_k-X_p)\cap
(W_k-X_q)|}\,.\qquad
\nonumber
\end{eqnarray}

As an alternative, which can be implemented in a more efficient way, we neglect the edge correction and consider the naive estimator $\widehat{\mathbf \Sigma}_k^{(2)}=\bigl(\cau\bigr)_{i,j=1}^{\ell}$ for ${\mathbf \Sigma}$ with
\beastar
\cau  &=&  \frac{1}{|W_k|}\,\sum_{p \ge 1} \ind_{W_k}(X_p) \bigl(\,\ind_{C_i \cap
C_j}(M_p)-P_M^{\o}(C_i)\,P_M^{\o}(C_j)\,\bigr)\\
&+& \frac{1}{|W_k|}\,\Su_{p,q\ge 1}\;\ind_{W_k}(X_p)\ind_{W_k}(X_q)\bigl(\,\ind_{C_i}(M_p) - 
P_M^{\o}(C_i)\,\bigr)\,\bigl(\,\ind_{C_j}(M_q) - P_M^{\o}(C_j)\,\bigr)\,.\quad
\nonumber
\eeastar

\smallskip
\begin{Th}\label{the.con.var}
Let $X_M$ be a stationary MPP satisfying $(\ref{int.con.alf})$ and let
$\{W_k\}$ be a CAS.  Then $\cec$ is an unbiased estimator, whereas  $\cau$  
is an asymptotically unbiased estimator for $\sigma_{ij}\,,\,i,j=1,...,\ell\,$.  
\end{Th}

{\bf Remark:} In general, neither $\cec$ nor  $\cau$ are $L^2$-consistent estimators 
for $\sigma_{ij}$, even if stronger moment and mixing conditions are supposed.

According to Lemma~\ref{lem.int.gam}, the integrability condition (\ref{int.con.alf}) in
Theorems~\ref{the.rep.tau} and~\ref{the.con.var} can be replaced by the stronger Condition $\beta(\delta)$. 
In order to obtain an $L^2$-consistent estimator,  
 we introduce a smoothed version of the unbiased estimator in (\ref{est.tau.one}), which is based on some kernel function and a sequence of bandwidths depending on the CAS $\{W_k\}$.

{\bf Condition $\b{(wb)}$:}  Let $w: \R \mapsto \R$ be a non-negative,
symmetric,     
Borel-measurable {\em kernel function} satisfying $w(x)\longrightarrow
w(0)=1$ as $x \to 0\,$.  In addition, assume that $w(\cdot)$ is bounded by 
$m_w<\infty$ and vanishes outside $B({\bf o},r_w)\,$ for some $r_w\in (0,\infty)$. 
Further, associated with $w(\cdot)$ and some given CAS $\{W_k\}$, let $\{b_k\}$ be a 
sequence of positive {\it bandwidths}  such that

\be
\label{cond:bk1} 
\frac{\varrho(W_k)}{2\,d\,r_w\,|W_k|^{1/d}}\geq b_k \lonk 0\;\;\;,\;\;\; 
b_k^d\,|W_k| \lonk \infty\;\;\;\mbox{and}\;\;\;
b_k^{\frac{3}{2}\,d}\,|W_k| \lonk 0\,. 
\ee

\medskip
\begin{Th} 
\label{thm:estTauThree} 
Let $\{W_k\}$ be an arbitrary CAS and $w(\cdot)$ be a 
kernel function with an associated sequence of bandwidths $\{b_k\}$ satisfying 
Condition $(wb)$. If the MPP $X_M$ satisfies 
\be 
\E\bigl(X_M(\,[0,1]^d \times \mathbb{M}\,)\bigr)^{4+\delta} < \infty\quad\mbox{and}\quad  
\int_1^\infty r^{d-1}\bigl(\beta_{X_M}^*(r)\bigr)^{\delta/(4+\delta)}\,{\rm d}r < \infty
\label{cond:ConsistencyTau}
\ee
for some $\delta >0\,$, then
\bee
\E\bigl(\,\sigma_{ij} - \cc\bigr)^2 \lonk 0 \;,
\eee
where $\cc$ is a smoothed covariance estimator defined by 
\beastar
&& \cc\,=\, \frac{1}{|W_k|} \sum_{p \ge 1} \ind_{W_k}(X_p)\,\bigl(\,\ind_{C_i \cap C_j}(M_p) - 
P^{\o}_M(C_i)\,P^{\o}_M(C_j)\,\bigr)\qquad  \\ 
&& + \; \Su_{p,q \ge 1}\;
\frac{\ind_{W_k}(X_p)\,\ind_{W_k}(X_q)\,\bigl(\,\ind_{C_i}(M_p) - P_M^{\o}(C_i)\,\bigr)
\bigl(\,\ind_{C_j}(M_q) - P_M^{\o}(C_j)\,\bigr)}{|(W_k-X_p) \cap (W_k-X_q)|}\,
w\Bigl(\frac{\|X_q-X_p\|}{b_k |W_k|^{1/d}}\Bigr)\,.
\nonumber
\eeastar
\end{Th}

\smallskip\noindent
{\bf Remark:}
The full strength of condition (\ref{cond:ConsistencyTau}) on the 
$\beta$-mixing rate $\beta_{X_M}^*(r)$ introduced in (\ref{bet.sta.one})
is only necessary to prove the consistency result of the preceding
Theorem~\ref{thm:estTauThree}. However, the $\beta$-mixing rate
$\beta_{X_M}^*(r)$ in Condition $\beta(\delta)$, which is needed to prove 
(\ref{int.cum.mea}) and (\ref{int.con.alf}) as well as
Theorem~\ref{the.asy.nor}, can be defined by the slightly smaller
non-increasing $\beta$-mixing rate function 
\be 
\beta_{X_M}^*(r) = \beta\bigl(\calA_{X_M}(\cube_a),\calA_{X_M}(\cube_{a+r}^c)\bigr) \quad\mbox{for}\quad 
r \ge a = 1/2\;.\label{cond:betaDeltaA}
\ee

\noindent
Moreover, in order to prove Theorem~\ref{the.asy.nor},
condition $\beta(\delta)$ relying on the $\beta$-mixing coefficient 
considered in (\ref{def.bet.emm}) with  $\beta_{X_M}^*(r)$ and $\beta_{X_M}^{**}(r)$
given in (\ref{cond:betaDeltaA}) and (\ref{bet.sta.one}), respectively, can be 
relaxed by using the slightly smaller $\alpha$-mixing coefficient
\be
\alpha\bigl(\calA_{X_M}(B),\calA_{X_M}(B^\prime)\bigr)=
\sup\{\,\bigl|\bbP(A\cap A^\prime)-\bbP(A)\bbP(A^\prime)\bigr|\,:
 A\in\calA_{X_{M}}(B),A^\prime\in\calA_{X_{M}}(B^\prime)\} 
\nonumber
\ee
instead of (\ref{def.bet.emm}). The corresponding  $\alpha$-mixing rates $\alpha_{X_M}^*(r)$ and
$\alpha_{X_M}^{**}(r)$ are then defined in analogy to (\ref{cond:betaDeltaA}) and
(\ref{bet.sta.one}), respectively. 
A covariance inequality for the $\alpha$-mixing case similar to $(\ref{cov.ine.bet2})$
can be found in {\rm \cite{Dou94}}, see \cite{Hei12a} for an improved 
version.
Despite of the subtle differences between the discussed mixing conditions,
we prefer to present our results under the unified assumptions of Condition
$\beta(\delta)$ and (\ref{cond:ConsistencyTau}) with  $\beta_{X_M}^*(r)$ as defined in (\ref{bet.sta.one}). It seems to be difficult to identify models where these differences are relevant. 

\section{Proofs}
\label{sec:proofs}

\subsection{Proof of Theorem~\ref{the.asy.nor}}

By the Cram\'{e}r-Wold technique, the multivariate CLT stated in (\ref{asy.equ.nor}) 
is equivalent to 
\begin{equation}
\label{asy.cra.wol}
 s^\top {\bm Y}_k = s_1\,Y_k(C_1) + \ldots + s_\ell\,Y_k(C_{\ell}) \Lonk \;{\mathcal N}_1(0\,,\,\sigma^2)
\quad\mbox{with}\quad \sigma^2 = s^\top\,{\bm \Sigma}\,s
\end{equation}
for any $s = (s_1,\ldots,s_\ell)^\top\in\R^\ell \ne {\o}_{\ell}\,$.

To prove (\ref{asy.cra.wol}) we extend  Bernstein's classical blocking method for weakly 
dependent random fields over a cubic index set of $\mathbb{Z}^d$, see e.g. \cite{Bul07},
\cite{Hei94} or \cite{Naha91}, to $\beta$-mixing fields indexed by elements of $H_k = \{z \in \mathbb{Z}^d: 
E_z \subset W_k \}\,$, where $E_z = [-1/2,1/2)^d +z$ for $z \in \mathbb{Z}^d$ and $\{W_k\}$ 
is an arbitrary CAS. The proof of (\ref{asy.cra.wol}) is divided  into four  steps. 

\medskip\noindent
{\bf Step 1. Bounds and asymptotics for the variance of the sum} \\[1mm]
In view of (\ref{yps.kah.til}) we may write  
\be 
s^\top {\bm Y}_k = \frac{1}{\sqrt{|W_k|}}\,\bigl(\, V_k + V_k^\prime\,\bigr)\,,\quad\mbox{where}
\quad V_k = \sum_{z\in H_k} U_z\quad\mbox{and}\quad V_k^\prime = 
\sum_{z\in \partial H_k}U_z^{(k)}
\nonumber
\ee
with 
\be
U_z^{(k)} = \sum_{n \ge 1} \ind_{E_z\cap W_k}(X_n)\,g(M_n)\quad,\quad
U_z = \sum_{n \ge 1} \ind_{E_z}(X_n)\,g(M_n)\quad\mbox{for}\quad z \in \mathbb{Z}^d\,,
\nonumber
\ee
$\partial H_k = \{z\in {\mathbb Z}^d\setminus H_k :|\,E_z\cap W_k\,| > 0\,\}$ and 
$g(M_n)=\sum_{i=1}^\ell s_i\bigl(\ind_{C_i}(M_n)-\pom(C_i)\bigr)\,$.
Clearly, $\E U_z^{(k)} = \E U_z = 0\,$ and $\max\{|U_z^{(k)}|,|U_z|\} \le c(s)\,
X_M(E_z \times \mathbb{M})$ 
for $z \in \mathbb{Z}^d\,$,  since $|g(M_n)| \le c(s)=|s_1|+\cdots+|s_{\ell}|\,$. 
Hence, by stationarity of $X_M$ and Condition $\beta(\delta)$,

\be
\max\{\|\,U_z^{(k)}\,\|_{2+\delta}\,,\,\|\,U_z\,\|_{2+\delta}\} \le c(s) \;
\|X_M([0,1)^d \times \mathbb{M})\|_{2+\delta}\quad\mbox{for}\quad z\in \mathbb{Z}^d\,.
\nonumber 
\ee

 In the following we use the maximum norm $|z|=\max_{1\leq i \leq d}|z_i|$ to express the 
distance  of $z=(z_1,...,z_d)\in\mathbb{Z}^d$ to the origin $\o$. By applying the covariance 
inequality (\ref{cov.ine.bet2}) together with Condition $\beta(\delta)\,$, we obtain
{\allowdisplaybreaks \begin{align*} \begin{split}
{\rm Var}(V_k^\prime) &= \sum_{y,z \in \partial H_k}\,\E\,U_y^{(k)}\,U_z^{(k)}   \le 
\sum_{y,z \in \partial H_k}\,\,\E\,| U_y^{(k)}\,U_z^{(k)} \,|  \le 
\#\partial H_k \sum_{z \in \mathbb{Z}^d}\,\,\E\,| U_{\bo}^{(k)}\,U_z^{(k)}
\,| 
\\
&\le 2\,c(s)^2\,\|X_M([0,1)^d \times \mathbb{M})\|^2_{2+\delta}\, \#\partial H_k
\,\sum_{z \in \mathbb{Z}^d}
\Bigl(\,\beta
\bigl({\calA_{X_M}(E_{\bo})},{\calA_{X_M}}(E_z)\bigr)\,\Bigr)^{\frac{\delta}{2+\delta}} 
\\
\end{split}\end{align*}}
{\allowdisplaybreaks \begin{align} \begin{split}
&\le 2\,c(s)^2\,\|X_M ( [0,1)^d \times \mathbb{M} )\|^2_{2+\delta}\, \#\partial H_k\,\Bigl(3^d + \sum_{z:\,|z|\geq 2} \bigl(\beta_{X_M}^*(|z|-1)\bigr)^{\frac{\delta}{2+\delta}}\,\Bigr)
\\
&\le  2\,c(s)^2\,\|X_M([0,1)^d \times \mathbb{M})\|^2_{2+\delta}\,
\#\partial H_k \Bigl(\,3^d + 2\,d\,\sum_{n \ge 1}(2n+3)^{d-1}\,
\bigl(\beta_{X_M}^*(n)\bigr)^{\frac{\delta}{2+\delta}}\,\Bigr)\,
\\
&\le c_1\,\#\partial H_k,
\label{estvkprime}
\end{split}\end{align}}%
$\mbox{for some constant }
c_1=c_1(s,d,\delta) > 0 $ where the relation $\#\{z\in\mathbb{Z}^d: |z| = n\}
=(2n+1)^d-(2n-1)^d \le 2\,d\,(2n+1)^{d-1}$ has been used. 
A simple geometric argument shows that each unit cube $E_z$ hitting the boundary  $
\partial W_k$ is contained in the annulus  $\partial W_k \oplus
B(\o,\sqrt{d})$ implying that

\be
\#\partial H_k  \le |\partial W_k \oplus B(\o,\sqrt{d})| \le 
2\,\bigl(\,|W_k \oplus B(\o,\sqrt{d})| - |W_k|\,\bigr)\,.   
\nonumber
\ee

Steiner's formula (cf. \cite{ScW08}, p. 600) applied to the convex
body $W_k$ reveals that the volume $|W_k \oplus B(\o,\sqrt{d})| - |W_k|$
 does not decrease when  $W_k$ is replaced by a larger convex body, e.g. by  $d^{3/2}\,R_k$ from  relation (\ref{John}) below, where  the hyper-rectangle
$R_k$ has edge lengths $a_1^{(k)},\ldots,a_d^{(k)}\,$. Replacing additionally
$B(\o,\sqrt{d})$ by the cube $[-\sqrt{d},\sqrt{d}]^d$ we get 
\bea
|W_k \oplus B(\o,\sqrt{d})| - |W_k| &\le & |d^{3/2}\,R_k \oplus [-\sqrt{d},\sqrt{d}]^d|
- |d^{3/2}\,R_k|\nonumber\\
&= & 2\, d^{(3d-2)/2}\,\sum_{i=1}^d a_1^{(k)}\cdots a_{i-1}^{(k)}
\Bigl(a_{i+1}^{(k)}+\frac{2}{d}\Bigr)\cdots \Bigl(a_d^{(k)}+\frac{2}{d}\Bigr)
\nonumber\\
&\le& 2^{d-1}\,d^{(3d-2)/2}\,\calH_{d-1}(\partial R_k),\quad\mbox{if}\quad
\min_{1 \le i \le d}a_i^{(k)} \ge \frac{2}{d}\,.
\nonumber\eea
Hence, since (\ref{John}) implies $\calH_{d-1}(\partial R_k) \le
\calH_{d-1}(\partial W_k)$ and $d^{3/2}\,\min_{1 \le i \le d}a_i^{(k)}\ge
2\,\varrho(W_k)\,$,  
it follows 
that 
$\#\partial H_k \le \,2^{d}\,d^{(3d-2)/2}\,\calH_{d-1}(\partial W_k)$ 
if $\varrho(W_k) \ge \sqrt{d}\,$, 
which in turn by combining (\ref{ineq:HeinrichPawlas}),
(\ref{estvkprime}) and the inclusion  $\# H_k \le |W_k| \le \# H_k + 
\#\partial H_k$ implies that

\be 
\label{var.vau.kah}  
\frac{{\rm Var}(V_k^\prime)}{|W_k|} \le c_2\,\frac{\calH_{d-1}(\partial W_k)}{|W_k|}
\le \frac{c_2\,d}{\varrho(W_k)} \lonk 0\qquad\mbox{and}\qquad \frac{\# H_k}{|W_k|} \lonk 1
\ee

for any CAS $\{W_k\}\,$.  Thus, by a standard Slutsky argument, (\ref{asy.cra.wol}) 
is equivalent to

\be 
\label{nte.ref.pro}
\frac{V_k}{\sqrt{\# H_k}}\;\Lonk \;{\mathcal{N}}_1(0\,,\,\sigma^2) \,.
\ee
 
The technique used above to estimate ${\rm Var}(V_k^\prime)$ will in the following be 
applied to show that 

\be
\label{sigma}
\sigma^2 \,=\, \lim_{k \to \infty}{\rm Var}(s^\top {\bm Y}_k)\; = \;
\lim_{k \to \infty}\frac{{\rm
Var}(V_k)}{\#H_k}\; =\; \sum_{z \in\mathbb{Z}^d} \E\,(U_{\o}\,U_z)\,.
\ee
The series on the rhs of $(\ref{sigma})$ converges absolutely
as immediate consequence of the estimate  
\be 
{\rm Var}(V_k) \;\le \; \#H_k\,\sum_{z \in\mathbb{Z}^d} |\,\E\,(U_{\o}\,U_z)\,|\;
 \le \; c_1\,\#H_k\,,
\nonumber
\ee
where the positive constant $c_1$ is the same as in (\ref{estvkprime}).
The Cauchy-Schwarz inequality and the previous estimates of ${\rm Var}(V_k)$ and 
${\rm Var}(V_k^\prime)$ show that

\be
\left|\,{\rm Var}(s^\top {\bm Y}_k) - \frac{{\rm Var}(V_k)}{|W_k|}\,\right|
\,\le\, 2\, \frac{|{\rm Cov}(V_k,V_k^\prime)|}{|W_k|} + \frac{{\rm Var}(V_k^\prime)}{|W_k|} 
\,\le\, \frac{2c_1\sqrt{\#H_k\,\#\partial H_k}}{|W_k|} + 
\frac{c_1\#\partial H_k}{|W_k|}\,
\nonumber 
\ee

proving the second equality in (\ref{sigma}). To prove the third equality in 
(\ref{sigma}) we use the identity
\bea
\frac{{\rm Var}(V_k)}{\#H_k} &=& \frac{1}{\#H_k}\sum_{y,z\in H_k} \E\,(U_{\o}\,U_{z-y}) = 
\sum_{z\in {\mathbb{Z}}^d } \frac{\#(H_k \cap (H_k-z))}{\#H_k}\E\,(U_{\o}\,U_z)
\nonumber
\eea
and the geometric inequality (following from the very definition of $H_k$ and 
$\partial H_k$)  
\be
\# (\,H_k \cap (H_k-z)\,) \le |\,W_k \cap (W_k-z)\,| \le \# (\,H_k \cap (H_k-z)\,) + 
\# \partial H_k + \# \partial(H_k-z)
\nonumber
\ee
for $z \in {\mathbb{Z}}^d\,$. This fact combined with (\ref{ineq:HeinrichPawlas}) and 
(\ref{var.vau.kah}) shows that
\be
\frac{\#(H_k \cap (H_k-z))}{\#H_k} \lonk 1\quad\mbox{for any fixed}\quad 
z \in {\mathbb{Z}^d}
\nonumber
\ee
proving the third equality in (\ref{sigma}) by applying the dominated convergence 
theorem.

\bigskip\noindent
{\bf Step 2. Passage to bounded random variables by truncation}\\[1mm]
\medskip\noindent
For any fixed $a > 0$ we define the random field $\{U_z(a)\,,\,z\in H_k\}$ of the 
truncated (and centered) random variables and the sum $V_k(a)$ by
\be
\label{tru.ran.var}
U_z(a) = U_z\ind_{\{|U_z|\le a\}} - \E\bigl(U_z\ind_{\{|U_z|\le a\}}\bigr)\quad
\mbox{and}\quad V_k(a) = \sum_{z\in H_k} U_z(a)
\ee
so that, for any $z\in H_k\,$,
\bee
|\,U_z(a)\,| \le 2\,a\quad\mbox{and}\quad \bigl(\,\E|\,U_z - U_z(a)\,|^{2+\delta}\,
\bigr)^{\frac{1}{2+\delta}}
= \|\,U_{\o} - U_{\o}(a)\,\|_{2+\delta} \lona 0\,. 
\eee 
By quite the same arguments as used in Step 1 based on the covariance 
inequality (\ref{cov.ine.bet2}) and  Condition $\beta(\delta)\,$, we find that 
\bee
{\rm Var}\bigl(\,V_k - V_k(a)\,\bigr) 
\le  2\,\# H_k \,\|\,U_{\o} - U_{\o}(a)\,\|_{2+\delta}^2\,
\Bigl(\,3^d + 2\,d\,\sum_{n \ge 1}(2n+3)^{d-1}\,
\bigl(\beta_{X_M}^*(n)\bigr)^{\frac{\delta}{2+\delta}}\,\Bigr)
\eee
for $k \ge 1\,$. Hence, by Slutsky's lemma, the weak limits of $V_k/\sqrt{\# H_k}$ and
$V_k(a)/\sqrt{\# H_k}$ as $k \to \infty$ are arbitrarily close  whenever  $a > 0$ 
is large enough. It therefore remains to prove the CLT in (\ref{nte.ref.pro}) 
 for the bounded random variables in (\ref{tru.ran.var}), i.e., for any fixed $a > 0\,$, 

\be
\label{clttrunc}
\frac{V_k(a)}{\sqrt{\# H_k}} \; \Lonk \;{\mathcal{N}}_1\bigl(\,0\,,\,\sigma^2(a)\,\bigr)
\quad\mbox{with}\quad \sigma^2(a) = \sum_{z\in\mathbb{Z}^d} \E\,U_{\o}(a)\,U_z(a)
\,. \ee

\bigskip\noindent
\newpage
{\bf Step 3. Adaptation of Bernstein's blocking method to non-cubic index sets}\\[1mm]

We start with some preliminary considerations. A well-known result from convex 
geometry first proved by F.~John, see e.g. \cite{Ball92}, asserts that there 
exists a unique ellipsoid $\mathcal{E}_k$ (called {\it John ellipsoid}) of maximal volume
contained in $W_k$ with midpoint $c(\mathcal{E}_k)$ and semi-axes of lengths $e_1^{(k)},\ldots,e_d^{(k)}$ 
such that $\mathcal{E}_k \subseteq W_k \subseteq c(\mathcal{E}_k) + d\,(\,\mathcal{E}_k - c(\mathcal{E}_k)\,)\,$.

Further, it is easy to determine a unique hyper-rectangle $R_k$ centered at the origin
$\o$ circumscribed by $\mathcal{E}_k-c(\mathcal{E}_k)$ with edge-lengths $a_i^{(k)} =  
2\,e_i^{(k)}/\sqrt{d}$ for $i=1,\ldots,d$ such that $\mathcal{E}_k -c(\mathcal{E}_k) \subseteq \sqrt{d}\,R_k$
and finally
\be 
\label{John} 
R_k \subseteq \mathcal{E}_k -c(\mathcal{E}_k) \subseteq W_k -
c(\mathcal{E}_k) \subseteq d\,(\,{\mathcal E}_k - c(\mathcal{E}_k)\,)
\subseteq d^{3/2}\,R_k\,.
\ee
Since  the  MPP $X_M$ observed in the CAS $\{W_k\}$ is stationary, we may
assume that $c(\mathcal{E}_k) = \o$ and without loss of generality let the edge lengths of $R_k$
be arranged in ascending order $a_1^{(k)} \le \cdots \le a_d^{(k)}$
(possibly after renumbering of the edges). Note that $R_k$ is not
necessarily in a position parallel to the coordinate axes. 
But there is an orthogonal matrix $O_k$  such that
\be
\label{rotation}
O_k\,R_k = \timesop_{i=1}^d \Bigl[-\frac{a_i^{(k)}}{2}\,,\,\frac{a_i^{(k)}}{2}\Bigr]\,.
\ee  

Let $\{p_k\}$ and $\{q_k\}$ be two sequences of positive integers (which will be specified 
later) satisfying $p_k \ge q_k \lonk \infty$ and $q_k/p_k
\lonk 0$. We define two types of pairwise disjoint cubes 
\bee
P_y^{(k)} = P_{\o}^{(k)} + (2\,p_k+q_k+1)\,y\quad\mbox{and}\quad
Q_y^{(k)} = Q_{\o}^{(k)} + (2\,p_k+q_k+1)\,y\quad\mbox{for}\quad y \in \mathbb{Z}^d\,,
\eee
where $P_{\o}^{(k)}= \{-p_k,\ldots,0,\ldots,p_k\}^d$ and $Q_{\o}^{(k)} = 
\{-p_k,\ldots,0,\ldots,p_k+q_k\}^d$ for $k \ge 1\,$. 

Now, we describe how to modify {\em Bernstein's blocking method} in order to 
prove the CLT stated in (\ref{clttrunc}). For the  family of  {\em block sums}

\bee
V_y^{(k)}(a) = \sum_{z \in P_y^{(k)}\cap H_k} U_z(a)\quad\mbox{for}\quad y \in G_k = 
\{z \in\mathbb{Z}^d: P_z^{(k)} \cap H_k \ne \emptyset\}
\eee
we shall show in Step 4 that 
\be
\label{cltVyk}
\frac{1}{\sqrt{\# H_k}}\,\sum\limits_{y \in G_k} V_y^{(k)}(a) \; \Lonk 
\;{\mathcal{N}}_1\bigl(\,0\,,\,\sigma^2(a)\,\bigr)
\ee
by assuming the mutual independence of the random variables $V_y^{(k)}(a)\,,\,y \in G_k\,$, 
which can be justified by Condition $\beta(\delta)$. Moreover, it will be sufficient to prove (\ref{cltVyk}), since
we will show below that
\be
\label{PconvVyk}
\frac{1}{\sqrt{\# H_k}}\Bigl(\,V_k(a) - \sum\limits_{y \in G_k} V_y^{(k)}(a)\,\Bigl) \; 
\PLonk \;0\,.
\ee

Next, we specify the choice of $p_k$ and $q_k$ in dependence on   
the edge lengths of $R_k$ and the supposed decaying rate of
$\beta^{**}_{X_M}(r)\,$. In view of 
$\varrho(W_k) \longrightarrow \infty$ and (\ref{John}) it follows that
$\min\{a_1^{(k)},\ldots,a_d^{(k)}\} \longrightarrow \infty$ as $k \to \infty\,$.
Note that the choice $p_k = \lfloor \e_k\,|W_k|^{1/2d}\rfloor$ as in case
of a cubic observation window with a certain null sequence $\{\e_k\}\,$, 
see \cite{Hei94}, does not always imply (\ref{cltVyk}) and
(\ref{PconvVyk}) if at least one of the first $d-1$ ordered edge lengths of $R_k$
increases very slowly to infinity. So one has to choose $p_k$  large enough
 but much smaller than $a_d^{(k)}$. For this purpose put
$r_k(s) = (\,a_{s+1}^{(k)}\cdot\ldots\cdot a_d^{(k)}\,)^{1/(2d-s)}$
for each $s \in \{0,1,\ldots,d-1\}\,$. Because of
$r^{2d-1}\,\beta^{**}_{X_M}(r) \longrightarrow 0$ as $r \to \infty$, there
exist non-increasing sequences  $\e_k(s)$ of positive numbers  such that
\be
\label{betaqk}
\e_k(s) \lonk 0,\quad\quad\e_k(s)\,r_k(s)\lonk \infty,\quad\mbox{and}\quad 
\frac{(r_k(s))^{2d-1}}{\e_k(s)}\,\beta^{**}_{X_M}(\e_k(s)\,r_k(s))\lonk 0 \,.
\ee 

Let $\e_k = \max\{\e_k(0),\ldots,\e_k(d-1)\}$ and $p_k(s) = 
\e_k^{1/(2d-s)}\,r_k(s)$ and select $s_k$ to be the smallest number 
 $s\in\{0,1,...,d-1\}$ such that $a_{s+1}^{(k)} \ge 2\,p_k(s)+1$ for 
$k \ge k_0\,$, where $k_0$ is a sufficiently large positive integer.
Thus, we define the integer sequences $p_k$ and $q_k$ by
\be
\label{pkqk}
p_k = \lfloor p_k(s_k)\rfloor = \lfloor\e_k^{1/(2d-s_k)}\,r_k(s_k)\rfloor
\quad \mbox{and}\quad q_k = \lfloor \e_k\,r_k(s_k)\rfloor\;\;\mbox{for}
\;\;k \ge k_0\,. 
\ee

Further, we need lower and upper bounds for the number $N_k$ of cubes 
${\widetilde Q_y}^{(k)} = [-p_k-\frac{1}{2},p_k+q_k+\frac{1}{2})^d+(2\,p_k+q_k+1)\,y$ 
hitting  $H_k\,$, i.e., $N_k = \#\{y\in \mathbb{Z}^d : {\widetilde
  Q_y}^{(k)} \cap H_k \ne \emptyset\}\,$.  For this put 
$N_j^{(k)}(c,w) = \#\{y\in \mathbb{Z}^d : {\widetilde Q_y}^{(k)} 
\cap (L_j^{(k)}(c) + w) \ne \emptyset\}$ for $w \in \R^d$ and some real $c
> 0\,$, where
$L_j^{(k)}(c) = O_k^T\,\{(x_1,\ldots,x_d)\in \R^d: -c\,a_j^{(k)}/2 \le x_j \le 
c\,a_j^{(k)}/2\,,\,x_i = 0\quad\mbox{for}\quad i \ne j\}\,$. 
The following rough estimates of $N_j^{(k)}(c,w)$ from below and above can be obtained  
by elementary geometric arguments: 
\be 
\Big\lfloor \frac{c\,a_j^{(k)}/\sqrt{d}}{2\,p_k+q_k+1}\Big\rfloor +
1\;\le\; N_j^{(k)}(c,w)\;\le\;
d\,\Bigl(\,\Big\lfloor \frac{c\,a_j^{(k)}}{2\,p_k+q_k+1}\Big\rfloor + 2\,\Bigr)
\quad\mbox{for any}\quad w \in \R^d\,.
\nonumber
\ee
Hence, by (\ref{John}) and (\ref{rotation}) the minimal number
$N_{\min}^{(k)}$ and 
the maximal number $N_{\max}^{(k)}$ of cubes ${\widetilde Q_y}^{(k)}$
hitting  $H_k$ satisfy the inequality
\be 
\prod_{j=1}^d\Bigl(\,\Big\lfloor
\frac{a_j^{(k)}/(2\,\sqrt{d})}{2\,p_k+q_k+1}\Big\rfloor +
1\,\Bigr) \;\le\; N_{\min}^{(k)}\; \le \; N_k \;\le
\;N_{\max}^{(k)}\;\le\; d^d\,\prod_{j=1}^d\Bigl(\,\Big\lfloor
\frac{d^{3/2}\,a_j^{(k)}}{2\,p_k+q_k+1}\Big\rfloor + 2\,\Bigr)\,.
\nonumber
\ee
  
In view of the  above choice of $s=s_k$ and (\ref{pkqk}), the number 
$N_k$ allows the estimate
\be 
\label{defNk} 
c_3\,\frac{a_{s_k+1}^{(k)} \cdots a_d^{(k)}}{p_k^{d-s_k}}\; \le \; N_k \; \le \; c_4\,
\frac{a_{s_k+1}^{(k)} \cdots a_d^{(k)}}{p_k^{d-s_k}}    
\quad\mbox{for all}\quad k \ge k_0
\ee
with positive constants $c_3, c_4$ only depending on the dimension $d\,$.
Combining the obvious fact that
$\#\,G_k \le N_k$ with (\ref{betaqk}), (\ref{pkqk}) and (\ref{defNk}) 
(with $p_k \ge 1$ and $\e_k \le 1$)  we arrive at
\bee
\#\,G_k\; p_k^{d-1}\,\beta_{X_M}^{**}(q_k) 
\le c_4\,\frac{(r_k(s_k))^{2d-1}}{\e_k^{1/2d}}\,\beta_{X_M}^{**}(q_k) \lonk 0\,.
\eee
Likewise, by (\ref{John}) and $a_i^{(k)} \le 2\,p_k+3$ for $i=1,...,s_k$,

\be
\label{sqrtHk}
 \frac{p_k^{d-s_k}}{\sqrt{\#\,H_k}}\,\prod_{j=1}^{s_k} a_j^{(k)} \le c_5\,
\left(\frac{p_k}{r_k(s_k)}\right)^{(2d-s_k)/2} \le c_5\,\sqrt{\e_k} \lonk 0\,.
\ee
Finally, we show that
\be
\label{GkPk} 
\frac{1}{\#\,H_k} \sum\limits_{y\in G_k} \#\,(P_y^{(k)}\cap H_k) \lonk  1\,,
\ee 
which, by the results of Step 1, is equivalent to
\be 
\label{GkPktilde} 
\frac{1}{|\,W_k\,|} \sum\limits_{y \in \mathbb{Z}^d} 
\bigl|\,\bigl(\,{\widetilde Q_y}^{(k)} \setminus {\widetilde P_y}^{(k)}\,\bigr) 
\cap W_k\,\bigr| \lonk  0\;,
\ee  
where ${\widetilde P_y}^{(k)} = [-p_k - \frac{1}{2},p_k + \frac{1}{2})^d + 
(2 p_k+q_k+1)\,y$ for $y \in \mathbb{Z}^d\,$. To estimate the volume of the space in 
$W_k$ outside the union of cubes ${\widetilde P_y}^{(k)}$ we introduce equidistant  
{\em slices} $S_{ij}^{(k)}$ in $\R^d$ of thickness $q_k$ and distance $2 p_k+1$
defined by
\be
S_{ij}^{(k)} = \Bigl\{(y_1,\ldots,y_d) \in \R^d : (2j+1)\,\bigl(p_k + \frac{1}{2}\bigr)
+j\,q_k \le y_i < (2j+1)\,\bigl(p_k + \frac{1}{2}\bigr) + (j+1)\,q_k\Bigr\}
\nonumber
\ee
for $i=1,...,d$ and $j \in \mathbb{Z}^1\,$. By (\ref{John}), (\ref{rotation}) 
and the choice of $p_k$ and $q_k$ it might happen that, for at most $s_k$
coordinates  $i\in \{1,...,d\}$,
$S_{ij}^{(k)} \cap W_k = \emptyset$ for all intergers $j\,$. 
For the remaining coordinates  $i \in \{1,...,d\}$  
there exist sequences of integers $n_k(i)$ (at least one of them tends to infinity
as $k \to \infty$) 
such that  
$S_{ij}^{(k)} \cap R_k \ne \emptyset$ for $|j|\le n_k(i)$ (and $S_{ij}^{(k)} 
\cap R_k = \emptyset$ for $|j| > n_k(i)$) and 
\be 
\frac{1}{|\,R_k\,|}\,\sum_{|j| \le n_k(i)} \,\bigl|\,S_{ij}^{(k)} \cap d^{3/2}\,R_k\,\bigr| 
\le c_6\, \frac{q_k}{p_k}\quad\mbox{for}\quad k\,\ge \,k_0\,,
\nonumber
\ee
where $c_6$ depends only on $d$. This estimate and the evident 
inequalities $|\,R_k\,| \le |\,W_k\,|$ and
\be  
 \sum\limits_{y \in \mathbb{Z}^d} 
\bigl|\,\bigl(\,{\widetilde Q_y}^{(k)} \setminus {\widetilde P_y}^{(k)}\,\bigr) 
\cap W_k\,\bigr| \le \sum_{i=1}^d \,\sum_{j\in \mathbb{Z}^1}\,\bigl|\,S_{ij}^{(k)} 
\cap d^{3/2}\,R_k\,\bigr|
\nonumber
\ee
show that the lhs of (\ref{GkPktilde}) is bounded by a constant multiple
of $q_k/p_k$ so that (\ref{GkPk}) is finally proved by (\ref{pkqk}).  

\bigskip\noindent
{\bf Step 4. Approximation by sums of independent random variables}\\[1mm]
\medskip
For brevity put $P_k = \bigcup_{y\in G_k}(P_y^{(k)} \cap H_k)\,$. Again  
by applying the 
covariance  inequality (\ref{cov.ine.bet2}) and Condition $\beta(\delta)$ 
to the stationary random field  $\{U_z(a), z\in H_k\}$ (with $|U_z(a)| \le 2\,a$ and thus 
$\delta = \infty$), we find in analogy 
to (\ref{estvkprime}) that
\bea
\frac{1}{\# H_k}\, \E\Bigl(\,V_k(a) - \sum_{y\in G_k}V_y^{(k)}(a)\Bigr)^2 &=& \frac{1}{\# H_k}\,\sum_{y,z\in H_k \setminus P_k} 
\E\bigl( U_y(a)\,U_z(a)\,\bigr)\nonumber\\
&& \nonumber \\
&\le& 8\,a^2\,\Bigl(\,3^d+2\,d\,\sum_{n \ge 1} (2n+3)^{d-1}\,\beta_{X_M}^*(n)\,\Bigr)\,
\frac{\#(H_k\setminus P_k)}{\#H_k}\,.
\nonumber
\eea

From (\ref{GkPk}) 
it is immediately clear that the ratio in the 
latter line disappears as $k \to \infty$, which confirms (\ref{PconvVyk}). Thus, in view of 
Slutsky's lemma, it remains to prove (\ref{cltVyk}).
We will do this under the assumption of  mutual independence of the block
sums $V_y^{(k)}(a), y\in G_k\,$. For this reason we show that the characteristic function
$\E\exp\{{\rm i}t\,\sum_{y\in G_k}V_y^{(k)}(a)\}$  differs from the product 
$\prod_{y\in G_k}\E\exp\{{\rm i}t\,V_y^{(k)}(a)\}$ uniformly in $t \in \R^1$ by
certain sequences tending to zero as $k \to \infty\,$.

Setting  $n_k = \# G_k$ we may write
\be
\xi_j = \exp\{{\rm i}t\,V_{y_j}^{(k)}(a)\}\quad\mbox{for}\quad y_j \in G_k 
\quad\mbox{with}\quad j=1,\ldots,n_k\,. 
\nonumber\ee
Using the algebraic identity 
\be
\E\prod_{j=1}^{n_k} \xi_j - \prod_{j=1}^{n_k} \E\xi_j = \sum_{j=1}^{n_k-1}
\E\xi_1\cdots \E\xi_{j-1}\,\Bigl(\,\E\xi_j\xi_{j+1}\cdots \xi_{n_k} -
\E\xi_j\,\E\bigl(\xi_{j+1}\cdots \xi_{n_k}\bigr)\,\Bigr)
\nonumber
\ee
and $|\xi_j| \le 1$ for $j=1,...,n_k$ we get 
\be 
\Big|\,\E\exp\{{\rm i}t\,\sum_{y\in G_k}V_y^{(k)}(a)\} - \prod_{y\in G_k}\E\exp\{{\rm i}t\,V_y^{(k)}(a)\}\,\Big|
\le \sum_{j=1}^{n_k-1}\Big|\,{\rm Cov}(\xi_j,\xi_{j+1}\cdots \xi_{n_k})\,\Big|\,.
\nonumber
\ee

By the stationarity of $X_M$ we may assume that the real as well as the imaginary part 
of $\xi_j$ is measurable w.r.t. the $\sigma-$algebra $\calA_{X_M}(\cube_{p_k+1/2})$
and the  product $\xi_{j+1}\cdots \xi_{n_k}$ is measurable w.r.t. 
$\calA_{X_M}(\R^d \setminus \cube_{p_k+q_k+1/2})$. By applying the covariance inequality
(\ref{cov.ine.bet2}) with $\delta = \infty$ (to the real and imaginary part of
$\xi_j$ resp. $\xi_{j+1}\cdots \xi_{n_k}$)  and using (\ref{bet.sta.one})
we find that 
\bea
\Big|\,{\rm Cov}(\xi_j,\xi_{j+1}\cdots \xi_{n_k})\,\Big| &\le&  8\,
\beta\bigl(\calA_{X_M}(\cube_{p_k+1/2}),\calA_{X_M}(\R^d\setminus \cube_{p_k+q_k+1/2})\bigr)
\nonumber\\ 
&\le&  8\,(p_k+1/2)^{d-1}\,\beta_{X_M}^{**}(q_k)\,.
\nonumber
\eea
Since $n_k = \# G_k \le N_k$ it follows with (\ref{defNk}) 
that 
\be 
\sup_{t \in \R^1}\,\Big|\,\E\exp\{{\rm i}t\,\sum_{y\in G_k}V_y^{(k)}(a)\} - \prod_{y\in G_k}\E\exp\{{\rm i}t\,V_y^{(k)}(a)\}\,\Big|
\le 8\,n_k\,(\,p_k + 1/2\,)^{d-1}\,\beta_{X_M}^{**}(q_k) \lonk 0\,.
\nonumber
\ee

The latter relation and the Berry-Esseen bound in the CLT for independent
random variables (which can be expressed by the third-order Lyapunov ratio, see e.g. \cite{Bul07}, p. 204, and references therein) reveal that (\ref{cltVyk}) holds if 
\be 
\label{Lyap} 
L_3^{(k)}(a) = \frac{1}{(\sigma_k^2(a))^{3/2}}\sum_{y \in G_k} \E|\,V_y^{(k)}(a)\,|^3 
\lonk 0 \qquad\mbox{and}\qquad \frac{\sigma_k^2(a)}{\# H_k} \lonk \sigma^2(a)\,,
\ee
where $\sigma^2(a)$ is defined by (\ref{clttrunc}) and $\sigma_k^2(a) = \sum_{y \in G_k}
\E(\,V_y^{(k)}(a)\,)^2$ coincides with the variance 
of $V_k(a)$ in case  of independent block sums $V_y^{(k)}(a),\,y \in G_k\,$.

It is easily seen that $|\,V_y^{(k)}(a)\,| \le 2\,a\,\#\,(P_y^{(k)}\cap H_k) \le 
 2\,a\,(2p_k+1)^{d-s_k}\,\prod_{i=1}^{s_k}(d^{3/2}\,a_i^{(k)}+1)$ and therefore
\be  
L_3^{(k)}(a) \le 2\,a\,\prod_{i=1}^{s_k}\bigl(d^{3/2} a_i^{(k)}+1\bigr)\,
\frac{(2\,p_k+1)^{d-s_k}}{(\sigma_k^2(a))^{3/2}}\sum_{y \in G_k} 
\E(\,V_y^{(k)}(a)\,)^2 \le c_7\,
\frac{2\,a\,p_k^{d-s_k}}{(\sigma_k^2(a))^{1/2}}\,\prod_{i=1}^{s_k}a_i^{(k)}
\nonumber
\ee 
with some positive constant $c_7$ only depending on $d\,$. In combination with (\ref{sqrtHk}) the second relation in (\ref{Lyap}) for $\sigma^2(a) > 0$  yields 
 $L_3^{(k)}(a) \lonk 0$. Hence, the first part of (\ref{Lyap}) is proved.

To accomplish the proof of (\ref{Lyap}) we remember that $\sigma^2(a)$ is the
asymptotic variance (\ref{sigma}) with  $V_k(a)$ from (\ref{tru.ran.var}) instead 
of $V_k$. Taking into account (\ref{GkPk}) or (\ref{GkPktilde}) we may replace $V_k(a)$
by the reduced sum $\sum_{y \in G_k}V_y^{(k)}(a)$ so that the second part of (\ref{Lyap})
 is a consequence of  
\be
\frac{1}{\# H_k}\,\Big|\,\E\Bigl(\,\sum_{y \in G_k}V_y^{(k)}(a) \,\Bigr)^2 - \sigma_k^2(a)\,\Big|
\le c_8\,a^2\,\frac{\# P_k}{\# H_k}\,\sum_{n \ge q_k} (2n+3)^{d-1}\,\beta_{X_M}^*(n) \lonk 0\,.
\nonumber
\ee

Here we have again used the notation $P_k$ and the standard covariance estimates from the very 
beginning of Step 4. Summarizing all Steps 1 - 4 completes the proof of 
Theorem~\ref{the.asy.nor}. \\ \halmos


%
%


\subsection{Proof of Lemma~\ref{lem.int.gam}}

By definition of the  signed measures $\gamma^{(2)}_X$  and 
$\gamma^{(2)}_{X,red}$  in Section~$\ref{sec:mom_cum_meas}$ and using algebraic 
induction, for any bounded Borel-measurable function $g : (\R^d)^2 \to \R^1\,$ we obtain the relation 
\be 
\label{def.gam.unm} 
\lambda\,\int\limits_{\R^d}\int\limits_{\R^d}
g(x,y)\,\gamma^{(2)}_{X,red}({\rm d}y)\,{\rm d}x = \int\limits_{(\R^d)^2} g(x,y-x)\,\gamma^{(2)}_X({\rm d}(x,y)).  
\ee 
Let $H^+, H^-$ be a Hahn decomposition of $\R^d$ for $\gamma^{(2)}_{X,red}$, i.e.,
\be \gamma^{(2)+}_{X,red}(\cdot)= \gamma^{(2)}_{X,red}(H^+ \cap (\cdot))\qquad\mbox{and}
\qquad \gamma^{(2)-}_{X,red}(\cdot)= -\gamma^{(2)}_{X,red}(H^-\cap (\cdot))\,.
\nonumber
\ee
We now apply (\ref{def.gam.unm}) for $g(x,y) = \ind_{E_{\o}}(x)\, \ind_{H^+\cap E_z}(y)\,$, 
where 
$E_z =
[-\frac{1}{2},\frac{1}{2})^d + z$ for $z\in \mathbb{Z}^d\,$. Combining this with
the definition of the (reduced) second factorial moment measures
$\alpha^{(2)}_X\,$ and $\alpha^{(2)}_{X,red}$ of the unmarked PP
$X = \sum_{i\ge 1} \delta_{X_i}\,$, see (\ref{def.red.alf}) for $m=2$, and 
\bee 
\gamma_X^{(2)}(A \times B) = \alpha_X^{(2)}(A \times B) - \lambda^2 \,|A|\, |B|\quad \mbox{for all bounded} 
\quad A,B\in \B\,,
\eee
leads to
\bea
\lambda\,\gamma^{(2)}_{X,red}(H^+\cap E_z) 
&=& \int\limits_{(\R^{d})^2}  \ind_{E_{\o}}(x)\ind_{H^+ \cap
  E_z}(y-x)\alpha^{(2)}_X({\rm d}(x,y)) - \lambda^2\,|E_{\o}|\,|H^+ \cap E_z|\nonumber\\
&&\nonumber\\
&=& \E \Su_{i,j\ge 1} \ind_{E_{\o}}(X_i)\ind_{H^+\cap
  E_z}(X_j-X_i) - \E X(E_{\o})\,\E X(H^ + \cap E_z).\nonumber\\
\nonumber
\eea
Since $\o \notin H^+ \cap E_z$ for $z\in \mathbb{Z}^d$ with $|z|\ge 2$ we
may continue with
\bea
\lambda\,\gamma^{(2)}_{X,red}(H^+ \cap E_z) &=&
\E \sum_{i\ge 1}\delta_{X_i}(E_{\o})\, X\bigl((H^+ \cap E_z) + X_i\bigr) - 
\E X(E_{\o})\,\E X(H^+\cap E_z) \nonumber \\
&=&  \E f(Y,Y^\prime_z) - \E f(\widetilde Y,\widetilde Y^\prime_z)
\quad\mbox{for}\quad |z| \ge 2\,,
\label{gam.equ1}
\eea
where 
\be
\label{gam.equ2}
f(Y,Y^\prime_z) = \sum_{i\ge 1}\delta_{X_i}(E_{\o})\, X\bigl((H^+ \cap E_z) + X_i\bigr)
\le X\bigl(E_{\o}\bigr)\, X\bigl(E_z\oplus E_{\o}\bigr)
\ee
with $Y(\cdot) = \sum_{i \ge 1}\delta_{X_i}\bigl((\cdot)\cap E_{\o}\bigr)$ resp. 
$Y^\prime_z(\cdot)=\sum_{j\ge 1}\delta_{X_j}\bigl(((\cdot) \cap E_z)\oplus E_{\o}\bigr)$
being restrictions of the stationary PP $X = \sum_{i \ge 1}\delta_{X_i}$ 
to $E_{\o}$ resp. $E_z\oplus E_{\o} = [-1,1)^d+z\,$. Further, let $\widetilde Y$ 
and $\widetilde Y^\prime_z$ denote copies of the PPs $Y$ and 
$Y^\prime_z\,$, respectively, which are assumed to be independent implying that 
$\E f(\widetilde Y, \widetilde Y'_z) = \E X(E_{\o})\,\E X(H^ + \cap E_z)\,$.
Since $Y$ is measurable w.r.t. $\mathcal{A}_X(E_{\o})$, whereas $Y'_z$ is
$\mathcal{A}_X(\R^d\setminus [-(|z|-1),|z|-1]^d )$-measurable,
we are in a position to apply Lemma~$\ref{lem.cov.bet}$ with
$\beta\bigl(\mathcal{A}_X(E_{\o}),\mathcal{A}_X(\R^d \setminus [-(|z|-1),|z|-1]^d \bigr)
\le \beta_{X_M}^*(|z|-\frac{3}{2})$ for $|z|\ge 2$. Hence, the estimate
(\ref{cov.ine.bet}) together with (\ref{gam.equ1}) and  (\ref{gam.equ2}) yields
\be
\bigl| \lambda\,\gamma^{(2)}_{X,red}(H^+\cap E_z) \bigr| \le
2\,\Bigl(\,\beta^*_{X_M}(|z|-\frac{3}{2})\,\Bigr)^{\frac{\eta}{1+\eta}}\,
\Bigl(\,\max\Bigl\{\,\E f^{1+\eta}(Y,Y^\prime_z)\,,\,\E
f^{1+\eta}(\widetilde Y,\widetilde Y^\prime_z)\,\Bigr\}\Bigr)^{\frac{1}{1+\eta}}\,,
\nonumber
\ee
where the maximum term on the rhs has the finite upper bound $2^{d(1+\eta)}\,\E(X(E_{\o}))^{2+2\eta}$ for $\delta = 2\,\eta > 0$ in accordance with our assumptions. This is seen from (\ref{gam.equ2}) using the Cauchy-Schwarz inequality and the stationarity of $X$ 
giving
\be
\E f^{1+\eta}(Y, Y'_z) \le \Bigl(\,\E(X(E_{\o})^{2+2\eta})\,
\E(X([-1,1]^d)^{2+2\eta})\,\Bigr)^{1/2} \le 2^{d(1+\eta)}\,\E(X(E_{\o})^{2+2\eta})
\nonumber
\ee
and the same upper bound for $\E f^{1+\eta}(\widetilde Y,\widetilde Y^\prime_z)\,$.
By  combining all the above estimates with $\lambda\,\gamma^{(2)}_{X,red}(H^+ \cap 
[-\frac{3}{2},\frac{3}{2})^d) \le 3^d \,\E X(E_{\o})^2$  we arrive at 
\be 
\lambda\,\gamma^{(2)}_{X,red}(H^+) \le 3^d \, \E X(E_{\o})^2 + 2^{d+1}\,
\Bigl(\E(X(E_{\o}))^{2+\delta}\Bigr)^{\frac{2}{2+\delta}}\,\sum_{z\in \mathbb{Z}^d : |z| \ge 2}
\Bigl(\,\beta^*_{X_M}(|z|-\frac{3}{2})\,\Bigr)^{\frac{\delta}{2+\delta}}.
\nonumber
\ee

By the assumptions of Lemma~\ref{lem.int.gam} the moments and the series on the rhs 
are finite and the same bound can be derived  for 
$-\lambda\,\gamma^{(2)}_{X,red}(H^-)$ which shows the validity of (\ref{int.cum.mea}). 

\bigskip\noindent
The proof of (\ref{int.con.alf}) resembles that of (\ref{int.cum.mea}). First we extend 
the identity (\ref{def.gam.unm}) to the (reduced) second factorial moment measure of the 
MPP $X_M$ defined by  (\ref{def.alf.isk}) and (\ref{cond:RadonNikodym}) for $m=2$ 
which reads as follows:

\bea
\lambda\,\int\limits_{\R^d} \int\limits_{\R^d} g(x,y)\,P_M^{\o,x}(C_1 \times C_2)\,
 \alpha_{X, red}^{(2)}({\rm d}y){\rm d}x &=& \int\limits_{(\R^d)^2} g(x,y-x)\,P_M^{x,y}(C_1\times C_2)\,\alpha_X^{(2)}\bigl({\rm d}(x,y)\bigr)\nonumber\\
&&\nonumber\\
&=& \E \Su_{i,j \ge 1} g(X_i,X_j-X_i) \ind_{C_1}(M_i)\ind_{C_2}(M_j)\,.
\nonumber
\eea 
For the disjoint Borel sets $G^+$ and $G^-$ defined by
\be
G^{+(-)} = \bigl\{\,x \in \R^d : \,P_M^{\o,x}(C_1\times C_2)\, \ge \,(<)\, 
P_M^{\o}(C_1)\,P_M^{\o}(C_2) \,\bigr\}
\nonumber\ee 
we replace  $g(x,y)$ in the above relation by  $g^{\pm}(x,y) =
\ind_{E_{\o}}(x)\, \ind_{E_z^{\pm}}(y)\,$, where $E_z^{\pm} = G^{\pm} 
\cap E_z$ for $|z| \ge 2\,$, and consider the 
restricted MPPs 
$Y_{\o}(\cdot) = X_M\bigl((\cdot)\cap (E_{\o}\times C_1)\bigr)$,   
$Y'_{z,\pm}(\cdot) = X_M\bigl((\cdot)\cap ((E_z^{\pm} \oplus E_{\o})\times C_2) \bigr)$ and their
copies  $\widetilde Y_{\o}$ and $\widetilde Y'_{z,\pm}\,$, which are 
assumed to be stochastically 
independent. Further, in analogy to (\ref{gam.equ2}), define
\be
f(Y_{\o},Y'_{z,\pm}) = \sum_{i \ge 1}\delta_{(X_i,M_i)}(E_{\o}\times C_1)\, 
X_M\bigl(\,(E_z^{\pm} + X_i) \times C_2\bigr)
\le X\bigl(E_{\o}\bigr)\, X\bigl( E_z \oplus E_{\o}\bigr)\,.
\nonumber
\ee

It is rapidly seen that, for $|z| \ge 2\,$, 
\begin{eqnarray*} 
\E f(Y_{\o},Y'_{z,\pm}) &=& \lambda\,\int\limits_{E_z^{\pm}} 
P_M^{\o,x}(C_1 \times C_2)\,\alpha^{(2)}_{X,red}({\rm d}x)\quad\mbox{and}\quad\\
\E f(\widetilde Y_{\o},\widetilde Y'_{z,\pm}) &=& 
\E X_M(E_{\o}\times C_1)\,\E X_M\bigl(E_z^{\pm} \times C_2\bigr) = 
\lambda^2\,P_M^{\o}(C_1)\,P_M^{\o}(C_2)\,|E_z^{\pm}|
\end{eqnarray*}
and in the same way as in the foregoing proof we find that, for $|z|\ge 2\,$,
\be 
|\,\E f(Y_{\o},Y'_{z,\pm}) - \E f(\widetilde Y_{\o},\widetilde Y'_{z,\pm})\,| 
\le 2^{d+1}\,\bigl(\E X(E_{\o})^{2+\delta}\bigr)^{\frac{2}{2+\delta}}\,
\bigl(\,\beta^*_{X_M}(|z|-\frac{3}{2})\,\bigr)^{\frac{\delta}{2+\delta}}\,.
\nonumber
\ee
 
Finally, the decomposition $\alpha_{X, red}^{(2)}(\cdot) = \gamma_{X,
  red}^{(2)}(\cdot) + \lambda\,|\cdot|$ together with the previous estimate
leads to 
\begin{eqnarray*}
&& \lambda\int\limits_{E_z} \Bigl| P_M^{\o,x}(C_1 \times C_2) - P_M^{\o}(C_1)\,P_M^{\o}(C_2)\Bigr|\,\alpha^{(2)}_{X,red}({\rm d}x)
= \E f(Y_{\bf o},Y'_{z,+}) - \E f(\widetilde Y_{\bf o},\widetilde Y'_{z,+})\\
&& - \bigl(\,\E f(Y_{\bf o},Y'_{z,-}) - 
\E f(\widetilde Y_{\bf o},\widetilde Y'_{z,-})\,\bigr) -
\lambda P_M^{\o}(C_1)\,P_M^{\o}(C_2)\,\Bigl(\gamma^{(2)}_{X,red}(\,E_z^+\,)
 - \gamma^{(2)}_{X,red}(\,E_z^-\,)\Bigr)\\
&&\\
&&  \le 2^{d+2}\,\bigl(\E(X(E_{\o})^{2+\delta})\bigr)^{\frac{2}{2+\delta}}\,
   \bigl(\beta_{X_M}^*(|z|-\frac{3}{2})\bigr)^{\frac{\delta}{2+\delta}}
+ \lambda \,|\gamma^{(2)}_{X, red}|(E_z)\quad\mbox{for}\quad |z| \ge 2\,.
\end{eqnarray*}

Thus, the sum over all $z \in \mathbb{Z}^d$ 
is finite in view of  our assumptions and the above-proved relation  
(\ref{int.cum.mea}) which completes the proof of Lemma~\ref{lem.int.gam}.
\halmos

\subsection{Proof of Theorem \ref{the.rep.tau}}

It suffices to show (\ref{dar.tau.til}), since independent marks
imply that  $P_M^{{\bf o},x}(C_1\times C_2)
= P^{\bf o}_M(C_1)\,P^{\bf o}_M(C_2)$ for $x \ne {\bf o}$ and any $C_1,C_2 \in \calB(\mathbb M)$ so that the integrand on the rhs of (\ref{dar.tau.til}) disappears which yields 
(\ref{tau.ind.mar}) for stationary independently MPPs. 
By the very definition of $Y_k(C)$ we obtain that
\bea
\quad {\rm Cov}\,\bigl( Y_k(C_i),  Y_k(C_j)\bigr) =
 \frac{1}{|W_k|} \;\E\sum_{p\ge 1} \ind_{W_k}(X_p)\bigl(\ind_{C_i}(M_p)-\pom(C_i)\bigr)
\bigl(\ind_{C_j}(M_p)-\pom(C_j)\bigr)&&\nonumber\\
+\; \frac{1}{|W_k|} \;\E{\sum\limits_{p,q \geq 1}}^{\neq}
\ind_{W_k}(X_p)\ind_{W_k}(X_q)\bigl(\ind_{C_i}(M_p)-\pom(C_i)\bigr)
\bigl(\ind_{C_j}(M_q)-\pom(C_j)\bigr)\,.\qquad &&
\label{streu2}
\eea
Expanding the difference terms in the parentheses leads to eight
expressions which, up to constant factors, take either the form
\beastar 
\E \sum_{p \ge 1} \ind_{W_k}(X_p) \ind_C(M_p) = \lambda |W_k|\,P^{\bo}_M(C) \quad\mbox{or}\quad\E{\sum\limits_{p,q \geq 1}}^{\neq}
\ind_{W_k}(X_p)\ind_{W_k}(X_q)\ind_{C_i}(M_p)\ind_{C_j}(M_q)&& \\
= \int\limits_{(\R^d)^2}\ind_{W_k}(x)\ind_{W_k}(y)\,P_M^{\bo,y-x}(C_i \times C_j)\,\alpha_X^{(2)}({\rm d}(x,y)) = \lambda\int\limits_{{\R}^d}P_M^{\bo,y}(C_i \times
C_j)\,\gamma_k(y) \, \alpha_{X,red}^{(2)}({\rm d}y)\,,&& 
\eeastar
where $y \mapsto \gamma_k(y)=|W_k \cap (W_k-y)|$ denotes the set covariance 
function of $W_k\,$. Summarizing all these terms gives
\begin{eqnarray*}
{\rm Cov}\,\bigl(Y_k(C_i),Y_k(C_j)\bigr)
 = \lambda \Bigl(\pom(C_i \cap C_j)-\pom(C_i)\pom(C_j)\Bigr)
+ \lambda\int\limits_{{\mathbb R}^d} \,\frac{\gamma_k(x)}{|W_k|}\;\Bigl( P_M^{\bo,x}(C_i \times
C_j)&&\\
\qquad - \,\, \pom(C_i)\,P_M^{\bo,x}(C_j \times \mathbb{M})
 -\pom(C_j)\,P_M^{\bo,x}(C_i \times \mathbb{M}) +
\pom(C_i)\,\pom(C_j)\Bigr)\;\alpha_{X,red}^{(2)}({\rm d}x)\,.&&
\end{eqnarray*}
The integrand in the latter formula is dominated by the sum
\be
\bigl| P_M^{\bo,x}(C_i \times
C_j)-\pom(C_i)\pom(C_j)\bigr|+ \bigl|P_M^{\bo,x}(C_j \times \mathbb{M})-\pom(C_j)\bigr| +\bigl|P_M^{\bo,x}(C_i \times
\mathbb{M})-\pom(C_i)\bigr|\,,
\nonumber\ee

which, by (\ref{int.con.alf}), is integrable w.r.t. $\alpha_{X,red}^{(2)}\,$. 
Hence, (\ref{dar.tau.til}) follows by (\ref{ineq:HeinrichPawlas}) and Lebesgue's 
dominated convergence theorem. \halmos

\subsection{Proof of Theorem \ref{the.con.var} }
We again expand the parentheses in the second term of the
estimator $\cec$ defined by
(\ref{est.tau.one}) and express the expectations in terms of $P_M^{\bo,y}$ and $\alpha_{X,red}^{(2)}$. Using the obvious relation $\gamma_k(y)=\int_{\R^d} \ind_{W_k}(x)
\ind_{W_k}(y+x)\,{\rm d}x$ we find that, for any  $C_i, C_j \in \calB(\mathbb{M})\,$, 

\beastar
\E\Su\limits_{p,q \geq 1}\frac{
\ind_{W_k}(X_p)\ind_{W_k}(X_q)\ind_{C_i}(M_p)\ind_{C_j}(M_q)}{|(W_k-X_p)\cap
(W_k-X_q)|}
= \int\limits_{(\R^d)^2} \frac{\ind_{W_k}(x) \ind_{W_k}(y)
P_M^{x,y}(C_i \times C_j)}{\gamma_k(y-x)}\; \alpha_X^{(2)}({\rm d}(x,y)) &&\\
= \lambda\int\limits_{\R^d} \frac{P_M^{\bo,y}(C_i \times
C_j)}{\gamma_k(y)} \int\limits_{\R^d} \ind_{W_k}(x)
\ind_{W_k}(y+x)\, {\rm d}x\, \alpha_{X,red}^{(2)}({\rm d}y)
=\lambda\int\limits_{\R^d}P_M^{\bo,y}(C_i \times C_j) \,
\alpha_{X,red}^{(2)}({\rm d}y)\,.&&
\eeastar

As in the proof of Theorem~\ref{the.rep.tau}  after summarizing all terms we obtain that

\smallskip\noindent
\begin{eqnarray*} \E \cec &=& \lambda
\Bigl(\pom(C_i \cap C_j)-\pom(C_i)\pom(C_j)\Bigr)
+\lambda \int\limits_{\R^d}\Bigl(P_M^{\bo,x}(C_i\times
C_j)\\
&-& \,P_M^{\bo,x}(C_i\times \mathbb{M})\pom(C_j)
- P_M^{\bo,x}(C_j\times \mathbb{M})\,\pom(C_i) + 
\pom(C_i)\,\pom(C_j)\Bigr)\,\alpha^{(2)}_{X,red}({\rm d}x)\,,
\end{eqnarray*}
which, by comparing with (\ref{dar.tau.til}), yields that $\E\cec = \sigma_{ij}\,$.
The asymptotic unbiasedness of $\cau$ is rapidly seen by the equality 
$\E\cau = {\rm Cov}\bigl(\,Y_k(C_i),Y_k(C_j)\,\bigr) = \E Y_k(C_i)Y_k(C_j)\,$, which follows
directly from (\ref{streu2}), and (\ref{cov.mat.tau}). 
\halmos

\subsection{Proof of Theorem~\ref{thm:estTauThree}} \label{sec:proofCovEst}

Since $\E\bigl( \sigma_{ij} - \cc \bigr)^2 = \Var \cc + 
\bigl( \sigma_{ij} - \E\cc \bigr)^2$
we have to show that 
\be  
\label{form:VarE}
\E \cc\lonk \sigma_{ij} \qquad\mbox{and}\qquad \Var \cc  \lonk 0\,. 
\ee
For notational ease, we put $\,m(u,v) = \bigl(\ind_{C_i}(u) -
\pom(C_i)\bigr) \bigl(\ind_{C_j}(v) - \pom(C_j)\bigr)\,,\,a_k = b_k |W_k|^{1/d}\,$,\\ 
\be
r_k(x,y) = \frac{\ind_{W_k}(x) \ind_{W_k}(y)}{\gamma_k(y-x)}\,w\left(\frac{\|y-x\|}{a_k}\right)
\quad\mbox{and}\quad \tau_k = {\sum_{p,q\ge 1}}^{\neq}\,r_k(X_p, X_q)\,m(M_p, M_q)\,.
\nonumber\ee\\
Hence, together with (\ref{erg.emp.mark}) and (\ref{yps.kah.til}) we may rewrite $\cc$ as follows:
\be
\label{streu3}
\cc =\frac{1}{\sqrt{|W_k|}}\;Y_k(C_i\cap C_j)+ \widehat{\lambda}_k \Bigl(\pom(C_i \cap C_j)-\pom(C_i)\pom(C_j)\Bigr)+\tau_k\,.
\ee

Using the definitions and relations (\ref{def.alf.isk}) -- (\ref{cond:RadonNikodym})
and $\int_{\R^d}r_k(x,y+x){\rm d}x = w\bigl(\|y\|/a_k\bigr)$ we find that the
expectation $\E\,\tau_k$ can be expressed by
\be
\int\limits_{(\R^d\times{\mathbb M})^2}r_k(x,y) m(u,v)
\alpha_{X_M}^{(2)}\bigl({\rm d}(x,u,y,v)\bigr) 
= \lambda\int\limits_{\R^d}
\int\limits_{{{\mathbb M}^{\phantom{}}}^2} m(u,v)P_M^{\bo,y}\bigl({\rm d}(u,v)\bigr)
w\Bigl(\frac{\|y\|}{a_k}\Bigr) \alpha_{X,red}^{(2)}\bigl({\rm d}y\bigr)\,. 
\nonumber\ee

The inner integral $\int_{{\mathbb M}^2} m(u,v)\,P_M^{\bo,y}\,\bigl({\rm d}(u,v)\bigr)$ coincides with the integrand occurring in (\ref{dar.tau.til}) and this term is integrable w.r.t. $\alpha_{X,red}^{(2)}$ due to 
(\ref{int.con.alf}) which in turn is a consequence of (\ref{cond:ConsistencyTau}) as shown in 
Lemma~\ref{lem.int.gam}. Hence, by Condition $(wb)$ and the dominated convergence theorem,
we arrive at 
\be
\E\,\tau_k \lonk \lambda\int\limits_{\R^d}\int\limits_{{{\mathbb M}^{\phantom{}}}^2}m(u,v)\,
P_M^{\bo,y}\bigl({\rm d}(u,v)\bigr)\,\alpha_{X,red}^{(2)}\bigl({\rm d}y\bigr)
= \sigma_{ij} - \lambda\bigl( \pom(C_i \cap C_j)- \pom(C_i)\pom( C_j) \bigr)\,.
\nonumber\ee
The definitions of $\widehat{\lambda}_k$ and $Y_k(\cdot)$ by (\ref{erg.emp.mark}) and 
(\ref{yps.kah.til}), respectively, reveal that 
$\E\,\widehat{\lambda}_k = \lambda$ and $\E\,Y_k(C_i\cap C_j) = 0\,$. This
combined with the last limit and (\ref{streu3}) proves the first relation
of (\ref{form:VarE}).
To verify the second part of (\ref{form:VarE}) we apply the Minkowski
inequality to the rhs of (\ref{streu3}) which yields the estimate 

\be 
\bigl(\Var\,\cc\bigr)^{1/2} \le |W_k|^{-1/2}\,\bigl(\Var\,Y_k(C_i\cap C_j)\bigr)^{1/2} 
+ \bigl(\Var\,\widehat{\lambda}_k\bigr)^{1/2} + \bigl(\Var\,\tau_k\bigr)^{1/2}\,.
\nonumber
\ee

The first summand on the rhs tends to 0 as $k \to \infty$ since $\E\,Y_k(C)^2$ 
has a finite limit for any $C\in {\mathcal B}({\mathbb M})$ as shown in Theorem~\ref{the.rep.tau}
under condition (\ref{int.con.alf}). The second summand is easily seen to
disappear as $k \to \infty$ if (\ref{int.cum.mea}) is fulfilled, see
e.g. \cite{Hei94}, \cite{Heinrich08} or \cite{Hei10}. 
Condition (\ref{cond:ConsistencyTau}) implies both (\ref{int.cum.mea}) and
(\ref{int.con.alf}), see Lemma~\ref{lem.int.gam}. Therefore, it remains to show that
$\Var\,\tau_k \longrightarrow 0$  as $k \to \infty\,$.

For this purpose we employ the variance formula (\ref{streuung}) stated in Lemma~\ref{Lemma:FormelMonster} in the special case $f(x,y,u,v) = r_k(x,y)\,m(u,v)\,$.  In this way we get the decomposition
$\Var\,\tau_k = I_k^{(1)} + I_k^{(2)} + I_k^{(3)}\,$, where $I_k^{(1)}$, $I_k^{(2)}$ and 
$I_k^{(3)}$ denote the three multiple integrals on the rhs of (\ref{streuung}) with $f(x,y,u,v)$
replaced by the product  $r_k(x,y)\,m(u,v)\,$. We will see that the integrals 
$I_k^{(1)}$ and $I_k^{(2)}$ are easy to estimate only by using (\ref{int.cum.mea}) and (\ref{int.con.alf}) while in order to show that  $I_k^{(3)}$ goes to  $0$ as $k \to \infty$,  the full strength of the mixing condition (\ref{cond:ConsistencyTau}) must be exhausted. Among others we use repeatedly the estimate 

\smallskip
\be 
\frac{1}{\gamma_k(a_k y)}\le \frac{2}{|W_k|}\quad\mbox{for}\quad y \in  B(\bo, r_w)\,, \label{form:gamma_k} 
 \ee

\smallskip
which follows directly  from (\ref{ineq:HeinrichPawlas}) and the choice of 
$\{b_k\}$ in (\ref{cond:bk1}). The definition of $I_k^{(1)}$  together with 
(\ref{form:gamma_k}) and $\alpha_{X, red}^{(2)}({\rm d}x)= \gamma_{X, red}^{(2)}({\rm d}x)+\lambda \, {\rm d}x$ yields
\bea
|I_k^{(1)}| &\le& 2\,\int\limits_{(\R^d)^2} \bigl(r_k(x_1,x_2)\bigr)^2 
\alpha_X^{(2)}\bigl({\rm d}(x_1,x_2)\bigr) = 2 \lambda \int\limits_{\R^d}  
\frac{1}{\gamma_k(y)} w^2\Bigl(\frac{\|y\|}{a_k}\Bigr)\alpha_{X,red}^{(2)}({\rm d}y) \nonumber \\ 
&\le& \frac{4\,\lambda}{|W_k|}\,\Bigl(\, m_w^2\;|\gamma_{X,red}^{(2)}|(\R^d) + 
\lambda\,a_k^d\;\int\limits_{\R^d} w^2(\|y\|){\rm d}y\,\Bigr) \lonk 0\,,
\nonumber\eea
where the convergence results from Condition $(wb)$ and (\ref{cond:ConsistencyTau}), which implies $|\gamma_{X, red}^{(2)}|(\R^d)<\infty$ by virtue of  Lemma~\ref{lem.int.gam}. Analogously, using besides (\ref{form:gamma_k}) and Condition $(wb)$ the relations
\be
w\Bigl(\frac{\|x\|}{a_k}\Bigr) \le m_w\,\ind_{[-\lceil a_k r_w\rceil, 
\lceil a_k r_w\rceil]^d}(x)\quad\mbox{and}\quad 
W_k \subseteq \bigcup_{z\in {\overline H_k}}E_z\quad\mbox{with}\quad 
{\overline H_k}=H_k\cup \partial H_k,
\nonumber\ee
with the notation introduced at the beginning of the proof of Theorem~\ref{the.asy.nor}, 
we obtain that
\bea
|I_k^{(2)}| &\le& 4\int\limits_{(\R^d)^3} r_k(x_1,x_2)\;r_k(x_1, x_3)\; 
\alpha_X^{(3)}\bigl({\rm d}(x_1,x_2,x_3)\bigr)\nonumber\\
&\le& \frac{16\,m_w^2}{|W_k|^2}\sum_{z \in {\overline H_k}} \alpha_X^{(3)}\bigl((E_z\oplus[-\lceil a_k r_w\rceil, \lceil a_k r_w\rceil]^d) \times (E_z\oplus[-\lceil a_k r_w\rceil, \lceil a_k r_w\rceil]^d) \times E_z  \bigr)\,.
\nonumber\eea

Since the cube $E_z\oplus[-\lceil a_k r_w\rceil, \lceil a_k r_w\rceil]^d$ 
decomposes into $(2\lceil a_k r_w\rceil+1)^d$ disjoint unit cubes and  
 $\alpha_X^{(3)}(E_{z_1} \times E_{z_2} \times E_{z_3} )\le \E(X(E_{\bo}))^3$ by Hölder's inequality, we may proceed with 

\be
|I_k^{(2)}| \le \frac{16\,m_w^2}{|W_k|^2}\; \#{\overline H}_k\; (2\lceil a_k r_w\rceil+1)^{2d}\;\E(X(E_{\bo}))^3\, \le c_9 \; b_k^{2d}\,|W_k| \lonk 0\,.
\nonumber\ee

Here we have used the moment condition in (\ref{cond:ConsistencyTau}), (\ref{var.vau.kah}), and the assumptions (\ref{cond:bk1}) imposed on the sequence $\{b_k\}\,$.

In order to prove that $I_k^{(3)}$ vanishes as $k \to \infty$ we first 
evaluate the inner integrals over the product $m(u_1,u_2)\,m(u_3,u_4)$
with $m(u,v)=\bigl(\ind_{C_i}(u)-P^{\bo}_M(C_i)\bigr)\bigl(\ind_{C_j}(v)
-P^{\bo}_M(C_j)\bigr)$ so that $I_k^{(3)}$ can be written as linear combination
of 16 integrals taking the form \\
\bea
J_k &=& \int\limits_{(\R^d)^2}\int\limits_{(\R^d)^2} r_k(x_1,x_2)\,r_k(x_3,x_4)\Bigl[ P^{x_1,x_2, x_3, x_4}_M(\timesop_{r=1}^4D_r)\,
\alpha_X^{(4)}\bigl({\rm d}(x_1,x_2, x_3, x_4)\bigr)\nonumber  \\ 
&&\qquad \quad - \; P^{x_1,x_2}_M(D_1\times D_2)\,
P^{x_3,x_4}_M(D_3\times D_4)\,\alpha_X^{(2)}\bigl({\rm d}(x_1,x_2)\bigr)\,
\alpha_X^{(2)}\bigl({\rm d}(x_3,x_4)\bigr) \Bigr]\nonumber\\
&&\nonumber\\ 
&=& \int_{\timesop_{r=1}^4(\R^d\times D_r)} r_k(x_1,x_2)\,r_k(x_3,x_4)\,
\bigl( \alpha_{X_M}^{(4)} - \alpha_{X_M}^{(2)}\times \alpha_{X_M}^{(2)}\bigr)
\bigl({\rm d}(x_1,u_1,...,x_4,u_4)\bigr)\,,
\nonumber\eea\\
where the mark sets  $D_1, D_3 \in \{C_i, {\mathbb M}\}$ and $D_2, D_4 
\in \{C_j, {\mathbb M}\}$ are fixed in what follows and the signed measure 
$\alpha_{X_M}^{(4)} - \alpha_{X_M}^{(2)}\times \alpha_{X_M}^{(2)}$ on 
${\mathcal B}((\R^d\times {\mathbb M})^4)\,$ $\big($and its total variation
 measure $\big|\alpha_{X_M}^{(4)} - \alpha_{X_M}^{(2)}\times 
\alpha_{X_M}^{(2)}\big|\,\big)$ come into play by virtue of the definition 
(\ref{cond:RadonNikodym}) for the $m$-point Palm mark distribution in case $m=2$ 
and $m=4\,$.

\smallskip
As $|z_1-z_2| > \lceil a_k r_w \rceil$ (where, as above, $|z|$ denotes the maximum norm of $z\in {\mathbb Z}^d$) implies $\|x_2-x_1\| > a_k r_w$ and thus $r_k(x_1,x_2) = 0$  for all 
$x_1 \in E_{z_1}, x_2 \in E_{z_2}$, we deduce from (\ref{form:gamma_k}) together with Condition $(wb)$ that
\be
|J_k|\le \frac{4\, m_w^2}{|W_k|^2}\;\Bigg(\sum_{n=0}^{2\,\lceil a_kr_w\rceil}+ 
\sum_{n > 2\,\lceil a_kr_w\rceil}\,\Bigg)\sum_{(z_1,z_2)\in S_k}\,
\sum_{(z_3,z_4)\in S_{k,n}(z_1)} V_{z_1, z_2, z_3, z_4}\,,
\label{form:Ik3Sum1}
\ee
where $S_k=\{(u,v) \in {\overline H_k}\times{\overline H_k}: |u-v|\le 
\lceil a_k r_w \rceil\}\,,\,S_{k,n}(z)=\{(z_1,z_2)\in S_k: 
\min\limits_{i=1,2}|z_i-z|=n\}$  and  $V_{z_1, z_2, z_3, z_4} =
\big|\alpha_{X_M}^{(4)} - \alpha_{X_M}^{(2)}\times \alpha_{X_M}^{(2)}\big|
\bigl(\times_{r=1}^4 (E_{z_r}\times D_r)\bigr)$ for any $z_1,...,z_4 \in {\mathbb Z}^d\,$.
 
Obviously, for any fixed  $z \in {\overline H_k}$, at most 
$2\,(2\,\lceil a_k r_w \rceil+1)^d\;(4\,\lceil a_k r_w \rceil+1)^d$ pairs 
$(z_3,z_4)$ belong to $\bigcup_{n=0}^{2\,\lceil a_k r_w \rceil} S_{k,n}(z)$
and the number of pairs $(z_1,z_2)$ in $S_k$ does not exceed the product $\#{\overline H_k}\;(2\,\lceil a_k r_w \rceil+1)^d$. Finally, remembering that $a_k = b_k\,|W_k|^{1/d}$ and using 
the evident\\ estimate $V_{z_1, z_2, z_3, z_4} \le 2\;\E(X(E_{\bo}))^4$ 
together with (\ref{var.vau.kah}) and Condition $(wb)$, we arrive at \\    
\be 
\frac{4\,m_w^2}{|W_k|^2}\,\sum_{(z_1,z_2)\in S_k}  
\sum_{n=0}^{2\,\lceil a_k r_w \rceil}\sum_{(z_3,z_4)\in S_{k,n}(z_1)} 
V_{z_1, z_2, z_3, z_4} \le c_{10}\,\frac{\#{\overline H_k}}{|W_k|^2}\,\bigl(b_k^d\,|W_k|\bigr)^3 \lonk 0\,.
\nonumber\ee
It remains to estimate the sums on the rhs of (\ref{form:Ik3Sum1}) running over $n > 2\,\lceil a_k r_w \rceil$.
For the signed measure $\alpha^{(4)}_{X_M} - \alpha^{(2)}_{X_M}\times \alpha^{(2)}_{X_M}$  we consider the Hahn decomposition $H^+, H^- \in {\mathcal B}((\R^d \times {\mathbb M})^4)$ yielding positive (negative) values on subsets of $H^+$($H^-$). Recall that $K_a=[-a,a]^d$.
For fixed $z_1 \in {\overline H_k}$, $z_2\in {\overline H_k} \cap (\cube_{\lceil a_k r_w \rceil}+z_1)$  
and $(z_3,z_4)\in S_{k,n}(z_1)$, we now consider  the decompsition $V_{z_1, z_2, z_3, z_4} = 
V^{+}_{z_1, z_2, z_3, z_4} + V^{-}_{z_1, z_2, z_3, z_4}$ with 
\be
V^{\pm}_{z_1, z_2, z_3, z_4} = \pm\,\big(\alpha_{X_M}^{(4)} - \alpha_{X_M}^{(2)}\times \alpha_{X_M}^{(2)}\big)\bigl(H^{\pm} \cap \timesop_{r=1}^4 (E_{z_r}\times D_r)\bigr).\,
\nonumber\ee

Since $(z_3,z_4)\in S_{k,n}(z_1)$ means that  $z_3 \in {\overline H_k} \cap \bigl(\cube^c_n+z_1\bigr)$, where $K^c_a=\R^d \setminus K_a\,$, and $z_4\in {\overline H_k}\cap \bigl(\cube_{\lceil a_k r_w \rceil}+z_3\bigr)\cap\bigl(\cube^c_n+z_1\bigr)$, we define MPPs $Y_k$ and $Y_n'$ as the restrictions of $X_M$ to 
$(\cube_{\lceil a_k r_w \rceil+1/2}+z_1)\times {\mathbb M}$ and $(\cube^c_{n-1/2}+z_1)\times {\mathbb M}\,$, respectively. Let furthermore  $\widetilde Y_k$ and $\widetilde Y_n'$ be copies of $Y_k$ and $Y_n'$ which are independent.

Next we define  functions $f^{+}(Y_k,Y_n^{\prime})$ and 
$f^{-}(Y_k,Y_n^{\prime})$ by\\
\be 
f^{\pm}(Y_k,Y_n^{\prime}) = \Su\limits_{p,q\ge 1}\;\Su\limits_{s,t\ge 1} 
\ind_{\pm}(X_p,M_p,X_q,M_q,X_s',M_s',X_t',M_t')\,,
\nonumber\ee
where $\ind_{\pm}(\cdots)$ denote the indicator functions of the sets 
$H^{\pm} \cap \timesop_{r=1}^4 (E_{z_r}\times D_r)$ so that we get\\
\be 
V^{\pm}_{z_1, z_2, z_3, z_4} = \E f^{\pm}(Y_k,Y_n^{\prime})- \E f^{\pm}(\widetilde Y_k, \widetilde Y_n')\quad\mbox{for}\quad (z_1,z_2)\in S_k\,,\,
(z_3,z_4)\in S_{k,n}(z_1)\,.
\nonumber\ee\\
Hence, having in mind the stationarity of $X_M$, we are in a position to 
apply the covariance inequality (\ref{cov.ine.bet}), 
which provides for $\eta > 0$ and $n > 2\,\lceil a_k r_w \rceil\,$ that
\bea
V^{\pm}_{z_1, z_2, z_3, z_4} &\le&
 2\Bigl(\,\beta \bigl( \calA(\cube_{\lceil a_k r_w \rceil+1/2}+z_1), 
\calA(\cube^c_{n-1/2}+z_1)\bigr)\,\Bigr)^{\frac{\eta}{1+\eta}}\qquad\nonumber\\
&\times& \Bigl(\,\E\big(\prod_{r=1}^2X_M(E_{z_r}\times D_r)\big)^{2+2\eta}\;
\E\big(\prod_{r=3}^4X_M(E_{z_r}\times D_r)
\bigr)^{2+2\eta}\,\Bigr)^{\frac{1}{2+2\eta}}\nonumber\\
&\le& 2 \bigl(\beta^*_{X_M}(n-\lceil a_kr_w\rceil -1)\bigr)^{\frac{\eta}{1+\eta}}\;\bigl(\E (X(E_{\bo}))^{4+4\eta}\bigr)^\frac{1}{1+\eta}\,.
\label{Vestimate}\eea\\
In the last step we have used the Cauchy-Schwarz inequality and the definition 
(\ref{bet.sta.one}) of the  $\beta$-mixing rate $\beta^*_{X_M}$. Finally, 
setting $\eta = \delta/4$ with $\delta > 0$ from (\ref{cond:ConsistencyTau})   
the estimate (\ref{Vestimate}) enables us to derive the following bound of 
that part on the rhs of (\ref{form:Ik3Sum1}) connected with the series over
 $n > 2\,\lceil a_k r_w \rceil$:
\bea
 c_{11}\,\frac{\#{\overline H_k}}{|W_k|^2}\,(2\lceil a_kr_w \rceil +1)^{2d}
\sum_{n > 2\lceil a_kr_w \rceil}\bigl( (2n+1)^d - (2n-1)^d\bigr) \bigr(\beta^*_{X_M}(n-\lceil a_kr_w\rceil -1)\bigl)^{\frac{\delta}{4+\delta}}\,.
\nonumber \eea 
Combining $a_k = b_k|W_k|^{1/d}$ and $(\ref{var.vau.kah})$ with condition (\ref{cond:ConsistencyTau}) and the choice of $\{b_k\}$  in (\ref{cond:bk1}), it is
easily checked that the latter expression and thus $J_k$ tend to 0 as 
$k \to \infty\,$. This completes the proof of Theorem~\ref{thm:estTauThree}.
\halmos 
 
\section{Examples} \label{sec:examples}
\subsection{$m$-dependent marked point processes}
A stationary MPP $X_M$ is called $m$-{\em dependent} if, for any  
$B, B'\in \B$, the $\sigma-$algebras $\mathcal{A}_{X_M}(B)$ 
and $\mathcal{A}_{X_M}(B')$ are stochastically
independent if $\inf\{|x-y|: \,x \in B,\, y\in B'\} > m$ or, equivalently,
\bee 
\beta \bigl(\mathcal{A}_{X_M}(K_a), \mathcal{A}_{X_M}(K^c_{a+b})\bigr) = 0
\quad\mbox{for}\quad b > m\;\;\mbox{and}\;\; a > 0\,.
\eee
In terms of the corresponding mixing rates this means  that $\beta_{X_M}^*(r)=\beta_{X_M}^{**}(r)=0$ if $r > m\,$.
For $m$-dependent MPPs $X_M$ it is evident that Condition $\beta(\delta)$ in 
Theorem~\ref{the.asy.nor} is only meaningful for $\delta = 0\,$, that is, 
$\E(X([0,1]^d))^2 < \infty\,$. This condition also implies 
(\ref{int.cum.mea}) and (\ref{int.con.alf}). Likewise, the assumption 
(\ref{cond:ConsistencyTau}) of Theorem~\ref{thm:estTauThree} reduces to 
$\E(X([0,1]^d))^4 < \infty$ which suffices to prove the $L^2$-consistency  
of the empirical covariance matrix $\widehat{\mathbf \Sigma}_k^{(3)}\,$.

\subsection{Geostatistically marked point processes}
Let  $X = \sum_{n\ge 1}\delta_{X_n}$ be an unmarked simple PP on $\R^d$ and  
$M = \{M(x),\,x \in \R^d\}$ be a measurable random field on $\R^d$ taking 
values in the Polish mark space ${\mathbb M}$. Further assume that $X$ and $M$ 
are stochastically independent over a common probability space 
$(\Omega, {\mathcal A}, {\mathbb P})$. An MPP $X_M = \sum_{n\ge 1}\delta_{(X_n,M_n)}$ with atoms  $X_n$ of $X$ and marks $M_n = M(X_n)$  is called {\em geostatistically marked}. Equivalently, the random counting measure 
$X_M \in{\mathsf N}_{\mathbb M}$ can be represented by means of the Borel 
sets $M^{-1}(C) = \{x\in \R^d: M(x) \in C\}$ (if $C\in {\mathcal B}({\mathbb M})$) by
\be
\label{geostat}
X_M(B \times C) = X(B \cap M^{-1}(C))\quad\mbox{for}\quad B \times C \in 
{\mathcal B}(\R^d)\times {\mathcal B}({\mathbb M})\,.
\ee
Obviously, if both the PP $X$ and the mark field $M$ are stationary then so is
$X_M$ and vice versa. Furthermore, the $m$-dimensional distributions of 
$M$ coincide with the $m$-point Palm mark distributions of $X_M$. The following Lemma allows to estimate the $\beta$-mixing
coefficient (\ref{def.bet.emm}) by the sum of the 
corresponding coefficients of the PP $X$ and the mark field $M$.

\begin{Lemma} Let the MPP $X_M$ be defined by $(\ref{geostat})$ with an unmarked 
PP and a random mark field $M$ being stochastically independent of each other. 
Then, for any $B, B'\in \mathcal{B}(\R^d)\,$,  
\be 
\label{bet.mix.coef}
\beta\bigl(\calA_{X_M}(B),\calA_{X_M}(B^\prime)\bigr) \le \beta\bigl(\calA_X(B),\calA_X(B^\prime)\bigr) + \beta\bigl(\calA_M(B),\calA_M(B^\prime)\bigr)\,,  
\ee
where the $\sigma-$algebras $\calA_X(B),\calA_X(B')$ and $\calA_M(B),\calA_M(B')$  are
generated by the restriction of $X$ and $M$, respectively, to the sets $B\,,B'$.
\end{Lemma}

To sketch a proof for (\ref{bet.mix.coef}), we regard the differences $\Delta(A_i,A_j')={\mathbb P}(A_i\cap A_j')- {\mathbb P}(A_i)\,{\mathbb P}(A_j')$ for two finite partitions $\{A_i\}$ and $\{A'_j\}$ of $\Omega$ consisting of events of the form 
\be
A_i = \bigcap_{p=1}^{k}\{X_M(B_p\times C_p)\in \Gamma_{p,i}\}\;,\;
A_j' = \bigcap_{q=1}^{{\ell}}\{X_M(B_q'\times C_q')\in \Gamma'_{q,j}\}
\quad\mbox{with}\quad \Gamma_{p,i},\Gamma'_{q,j} \subseteq {\mathbb Z}^1_+\,,
\nonumber\ee
with pairwise disjoint bounded Borel sets $B_1,...,B_k \subseteq B$ and
$B_1',...,B_{\ell}' \subseteq B'$.  Making use of 
(\ref{geostat}) combined with the independence assumption yields the identity 
\bea
\Delta(A_i,A_j') &=&
\int\limits_{\Omega}\int\limits_{\Omega}\Bigl({\mathbb P}_{\calA_X(B)\otimes\calA_X(B^\prime)}-{\mathbb P}_{\calA_X(B)}\times {\mathbb P}_{\calA_X(B^\prime)}\Bigr)(A_i\cap A_j')\,{\rm d}{\mathbb P}_{\calA_M(B)\otimes\calA_M(B^\prime)}
\nonumber\\
&+& \int\limits_{\Omega}\int\limits_{\Omega}{\mathbb P}_{\calA_X(B)}(A_i)\,
{\mathbb P}_{\calA_X(B^\prime)}(A_j')\,{\rm d}\Bigl({\mathbb P}_{\calA_M(B)\otimes\calA_M(B^\prime)}- {\mathbb P}_{\calA_M(B)}\times {\mathbb P}_{\calA_M(B^\prime)}\Bigr)\,,
\nonumber\eea
which by (\ref{def.bet.emm}) and the integral form of the total variation confirms (\ref{bet.mix.coef}).

\subsection{Cox processes on the boundary of Boolean models}
Let $\Xi = \bigcup_{n \ge 1} (\Xi_n+Y_n)$ be a {\em Boolean model}, see e.g. \cite{molchanov97}, governed  by some stationary Poisson process $\sum_{n\ge 1} \delta_{Y_n}$ in $\R^d$ with intensity $\lambda>0$  and a sequence $\{\Xi_n\}_{n\ge 1}$ of independent copies of some random convex, compact set $\Xi_0$
 called {\em typical grain} (where we may assume that $\bo \in \Xi_0$). With the radius 
functional $\diam = \sup\{\|x\|:\,x \in \Xi_0\}$, the condition $\E \diam^d < \infty$
ensures that $\Xi$ is a well-defined random closed set.
We consider a marked Cox process $X_M$,
where  the unmarked Cox process $X=\sum_{n\ge 1}\delta_{X_n}$ is concentrated on the boundary $\partial\Xi$ of $\Xi$ with 
 random intensity   measure being proportional 
to the $(d-1)$-dimensional Hausdorff measure $\mathcal{H}_{d-1}$ on 
$\partial\Xi$.
As marks $M_n$ we take the outer unit normal vectors  at the points
 $X_n \in \partial\Xi$, which are (a.s.) well-defined for $n\ge 1$ due
to the assumed convexity of  $\Xi_0$. 
This example with marks given by the orientation of outer normals in random boundary points may occur rather specific. However, this way our asymptotic results may be used to construct asymptotic tests for the fit of a Boolean model to a given dataset w.r.t. its rose of directions. For instance, if the typical grain 
$\Xi_0$ is rotation-invariant (implying the isotropy of $\Xi$), then 
the Palm mark distribution $\pom$ of the stationary MPP $X_M = \sum_{n\ge 1}\delta_{(X_n,M_n)}$ is the uniform distribution on the unit sphere ${\mathbb S}^{d-1}$ in $\R^d$. We will now discuss assumptions ensuring that Condition $\beta(\delta)$ and (\ref{cond:ConsistencyTau}) hold, which are required for our CLT (\ref{asy.equ.nor}) and the consistent estimation of the covariances (\ref{cov.mat.tau}), respectively. 
Using slight modifications of the proofs for Lemmas 5.1 and 5.2 in \cite{Hei99} one can show that for $a,b>0$ 
\bee 
\beta \bigl(\calA_{X_M}(K_a),\calA_{X_M}(K^c_{a+b}) \bigr) \le \lambda\,2^{d+2}\,\left(3+\frac{4a}{b}\right)^{d-1} 
\E \bigl(\diam^d\ind\{\diam \ge b/4\}\bigr)\,.
\eee
According to (\ref{bet.sta.one}) we may thus define the $\beta$-mixing rates 
$\beta_{X_M}^*(r)$ and $\beta_{X_M}^{**}(r)$  for $r \ge \frac{1}{2}$ by
\bea 
\beta_{X_M}^*(r) &=& c_{12}\,  
\E \bigl(\diam^d\ind\{\diam \ge r/4\}\bigr) \ge \sup_{1/2\le a \le r}\beta \bigl(\calA_{X_M}(K_a),\calA_{X_M}(K^c_{a+r}) \bigr)\,,\nonumber\\ 
&&\nonumber\\
\beta_{X_M}^{**}(r) &=& \frac{c_{12}}{r^{d-1}}\,  
\E \bigl(\diam^d\ind\{\diam \ge r/4\}\bigr) \ge \sup_{a \ge r}a^{-(d-1)}\,\beta \bigl(\calA_{X_M}(K_a),\calA_{X_M}(K^c_{a+r}) \bigr)\,, 
\nonumber \eea
where $c_{12}=\lambda 2^{d+2}7^{d-1}\,$. 

\smallskip\noindent
It is easily seen that $\E \diam^{2d} < \infty$ implies $r^{2d-1}\beta_{X_M}^{**}(r) \lonr 0$. Moreover, $\E\diam^{2d(p+\delta)/\delta+\e} < \infty$ for some $\e>0$
ensures $\int_1^\infty r^{d-1} \bigl( \beta^*_{X_M}(r)
\bigr)^{\delta/(2p+\delta)} {\rm d}r < \infty$ for $p\ge 0$. 
Since the 
random intensity measure of $X$ on $E_{\bo}$ and thus also $X(E_{\bo})$ 
has moments of any order by virtue of $\E \diam^d < \infty\,$, the
parameter $\delta > 0$ in Condition $\beta(\delta)$ and 
(\ref{cond:ConsistencyTau}) can be chosen arbitrarily large. This results in the following lemma.
 
\begin{Lemma} For the above-defined stationary marked Cox process $X_M$ on the 
boundary of a Boolean model $\Xi$ with typical grain $\Xi_0$ the assumptions of 
the Theorems~$\ref{the.asy.nor}$ and~$\ref{thm:estTauThree}$  are satisfied 
whenever 
\be
\label{typ.grain} 
\E \|\Xi_0\|^{2d+\e} < \infty \quad\mbox{for some}\quad \e > 0\,.  
\ee
\end{Lemma}

{\bf  Remark:} The marked Cox process $X_M$ is $m$-dependent if $\|\Xi_0\|$ 
is bounded by some constant. By using  approximation techniques with 
truncated grains as suggested in \cite{Hei99}, pp. 299-302, it can 
be shown that (\ref{typ.grain}) is just needed for $\e = 0\,$. Moreover, the statistical analysis of roses of directions  via marked Cox 
processes applies also in case of non-Boolean $\beta$-mixing fibre processes, see
e.g. \cite{Hei94} for Voronoi tessellations.  

\section{Simulation study}\label{sec.sta.app}

\subsection{Moving average model in $\R^2$}\label{sub.sta.seg}

In this section we introduce an $m$-dependent MPP model, which was used for our simulations since it allows to control
the range of spatial dependence for a fixed Palm mark distribution. 
The locations of this MPP are given by a homogeneous Poisson process $\sum_{n\ge 1} \delta_{X_n}$ in $\R^2$. Each point is marked by a direction
in the upper half ${\mathbb S}_+^{1}$ of the unit circle. In order to construct the Palm mark distribution, we
consider the projected normal distribution PN$_2(a,
\bm\kappa)$ on ${\mathbb S}^1$. By definition, $Y\sim{\rm}{\rm PN}_2(a,
\bm\kappa)$ means that $Y=\frac{Z}{\|Z\|}$ for some Gaussian random vector $Z\sim\mathcal{N}_2(\bm a,\bm\kappa)$ in $\R^2$ with an invertible covariance matrix $\bm\kappa$.
Note that  PN$_2(\bo,
\sigma^2 I_2)$ is the uniform distribution on ${\mathbb S}^1$ for all $\sigma^2>0$, where $I_2$ is the identity matrix. Formulas for the density of a projected normal distribution can be found in \cite{Mardia00}. 
Let $\{M_n^{(1)}\}_{n \ge 1}$ be iid N$_2(\bo,\bm\kappa)$-distributed random vectors. The stability of the normal distribution w.r.t. convolution yields 
\bee M_n^{(2)}=\frac{\sum_{i=1}^\infty M_i^{(1)} \ind_{\{\|X_i-X_n\| \le \rho\}}}{\left\|\sum_{i=1}^\infty M_i^{(1)} \ind_{\{\|X_i-X_n\| \le \rho\}}\right\|}\sim {\rm PN}_2(\bo,\bm\kappa),\eee
for any $\rho\ge 0$ controlling the range of dependence. 
The marks of our model are finally defined as the axial
versions $M_n=M_n^{(2)} \ind_{\mathbb{S}^+_1}(M_n^{(2)})-M_n^{(2)}\ind_{\mathbb{S}^-_1}(M_n^{(2)})$ of the averages $M_n^{(2)}$, i.e., points on the lower half-circle ${\mathbb S}_-^1$ are rotated by $\pi$. 
Due to the moving average approach defining the preliminary marks $\{M^{(2)}_n\}$, we call the MPP $X_M=\sum_{n \ge 1} \delta_{(X_n, M_n)}$ the {\it moving average model } (MAM). The MAM is clearly $m$-dependent, where the range of dependence is controlled by the averaging parameter $\rho$.

\subsection{Tests}\label{sec:tests}

By simulations of the MAM we investigated the performance of the asymptotic $\chi^2$-goodness-of-fit test, which is based on the test statistic
\bee
T_k={\bf Y}_k^\top \widehat{\covM}_k^{-1}{\bf Y}_k \Lonk \chi^2_\ell.
\eee
If $(\widehat{\covM})_k$ is chosen as the $L^2$-consistent estimator $\cc$, and $(\pom)_{0}^{\phantom{o}}$ denotes a hypothetical Palm mark distribution, the hypothesis $H_0:\pom=(\pom)_{0}^{\phantom{o}}$ is rejected, if
$T_k>\chi^2_{\ell,1-\alpha}$, where $\alpha$ is the level of significance, and $\chi^2_{\ell,1-\alpha}$ denotes the $1-\alpha$-quantile of the $\chi^2_\ell$-distribution.
This test will be referred to as {\it \lq test for the typical mark distribution\rq}\, (TMD).
The construction of $\cc$ involves the sequence of bandwidths $\{b_k\}$. We set
\be b_k=c |W_k|^{-\frac{3}{4d}} \quad \mbox{ for a constant } c>0,\label{form:b_k_kappa}\ee
and chose $c$ such that condition ($wb$) for the $L^2$-consistency of $\cc$ was satisfied. 
The specific choice of the constant $c$ is crucial for test performance, as discussed below. 
This choice of $c$ can be avoided if $\covM$ is not estimated
from the data to be tested but incorporated into $H_0$. This means, we specify an MPP as null model, such that $\covM_0$ is either theoretically known or otherwise can be approximated by Monte-Carlo simulation.
By means of the combined null hypothesis $H_0: \pom=(\pom)_{0}^{\phantom{o}}$ and ${\covM}={\covM}_0$,
the test exploits
not only information on the distribution of the typical mark but additionally considers asymptotic effects of
spatial dependence. The test can thus be used to investigate if a given point pattern differs from the MPP null model w.r.t. the Palm mark distribution.  We will therefore refer to it as {\it \lq test for mark-oriented goodness of model fit\rq}
(MGM). By the  strong law of large numbers and the asymptotic unbiasedness of $\cau$,
a strongly consistent Monte-Carlo estimator for $\covM_0$ in an MPP model $X_M$ is given by
\bee \widehat{\covM}_{k,n}=\frac{1}{n}\sum_{\nu=1}^n \cau( X_M^{(\nu)}),\eee
where $X_M^{(1)}, \ldots, X_M^{(n)}$ are independent realizations of $X_M$. 
Thus, for large $k$ and $n$ the test statistic $T_{k,n}={\bf Y}_k^\top \widehat{\covM}_{k,n}^{-1}{\bf Y}_k $ has an approximate $\chi^2_\ell$ distribution. If $\alpha$ is the level of significance, the MGM test rejects $H_0$, if $T_{k,n}>\chi^2_{\ell,1-\alpha}$.
The estimator $\widehat{\covM}_{k,n}$ can also be used to construct a test for the typical mark distribution if independent replications of a point patterns are to be tested. In that case $X_M^{(1)}, \ldots, X_M^{(n)}$ are the replications. Note that for replicated point patterns, $H_0$ does not incorporate an assumption on $\covM$ and hence the corresponding test differs from the MGM test. The edge-corrected unbiased estimator $\cec$ was not used for the Monte-Carlo estimates in our simulation study, since $\cau$ can be computed more efficiently. \\

\subsection{Model parameters}\label{sec:parameters}
The MAM was simulated on the observation window $W_{1500}=[-1500, 1500]^2$. The expected number of points
was set to  $\E X(W_{1500})=3125$.
The asymptotic behavior of the test was studied by considering smaller observation windows 
corresponding to an expected number of $300$, $600,\ldots,3000$ points.
Spatial stochastic dependence of marks was varied by the parameter $\rho \in \{0,50,\ldots,300\}$.
In the MAM, marks of points with distance no larger than $2\rho$ in general exhibit stochastic dependence. If, on the contrary, two points are separated by more than $2\rho$, their marks are independent. Thus, $\rho =0$ corresponds to independent marking.
Deviations of the projected normal distribution from the uniform distribution on ${\mathbb S}_+^1$ were controlled by 
 varying $\kappa_{12}\in\{0, 0.1, 0.2, 0.4, 0.8\}$, where 
$\kappa_{12}=0$ represents the uniformly distributed case. The parameter $\kappa_{11}=\kappa_{22}=1$ was kept constant.
The bins
$C_1,\ldots, C_\ell\in\mathcal{B}({\mathbb S}_+^1)$ for the $\chi^2$-goodness-of-fit test were chosen as \bee C_i=\left\{(\cos\theta,
\sin \theta)^T :\; \theta \in\left[(i-1)\frac{\pi}{\ell+1}, i\frac{\pi}{\ell+1}\right)\right\},\; i=1,\ldots,\ell.\eee
We will discuss  the case $\ell=8$, where
the bins had a width of $20^\circ$. Simulations for $\ell=17$ did not reveal different general effects.

\subsection{Simulation results}

\begin{figure} %
  \centering 
  \subfigure[Type I error]
  {\includegraphics[width=5.9cm, angle=-90]{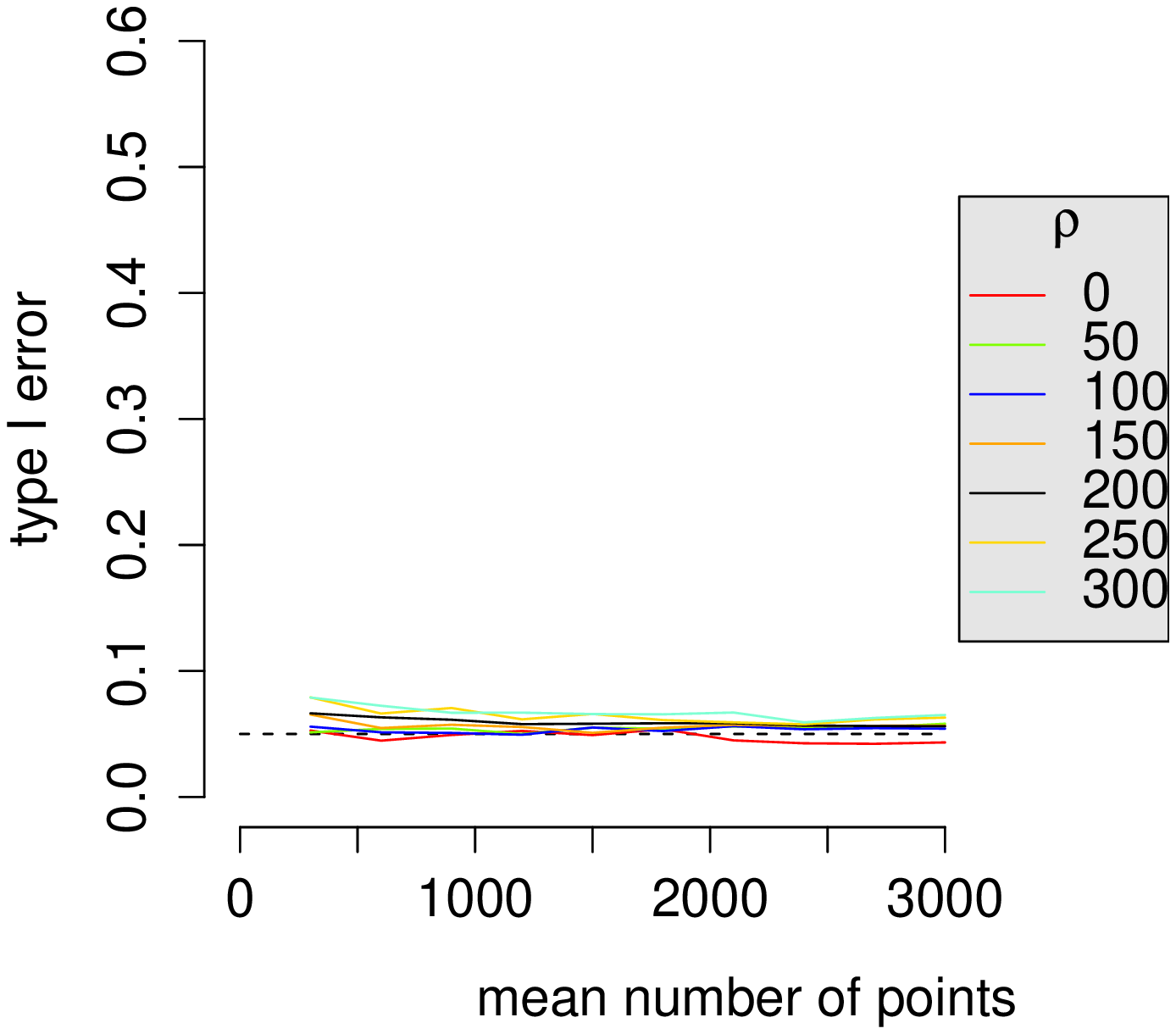}}
\subfigure[Type II error]
{\includegraphics[width=5.9cm, angle=-90 ]{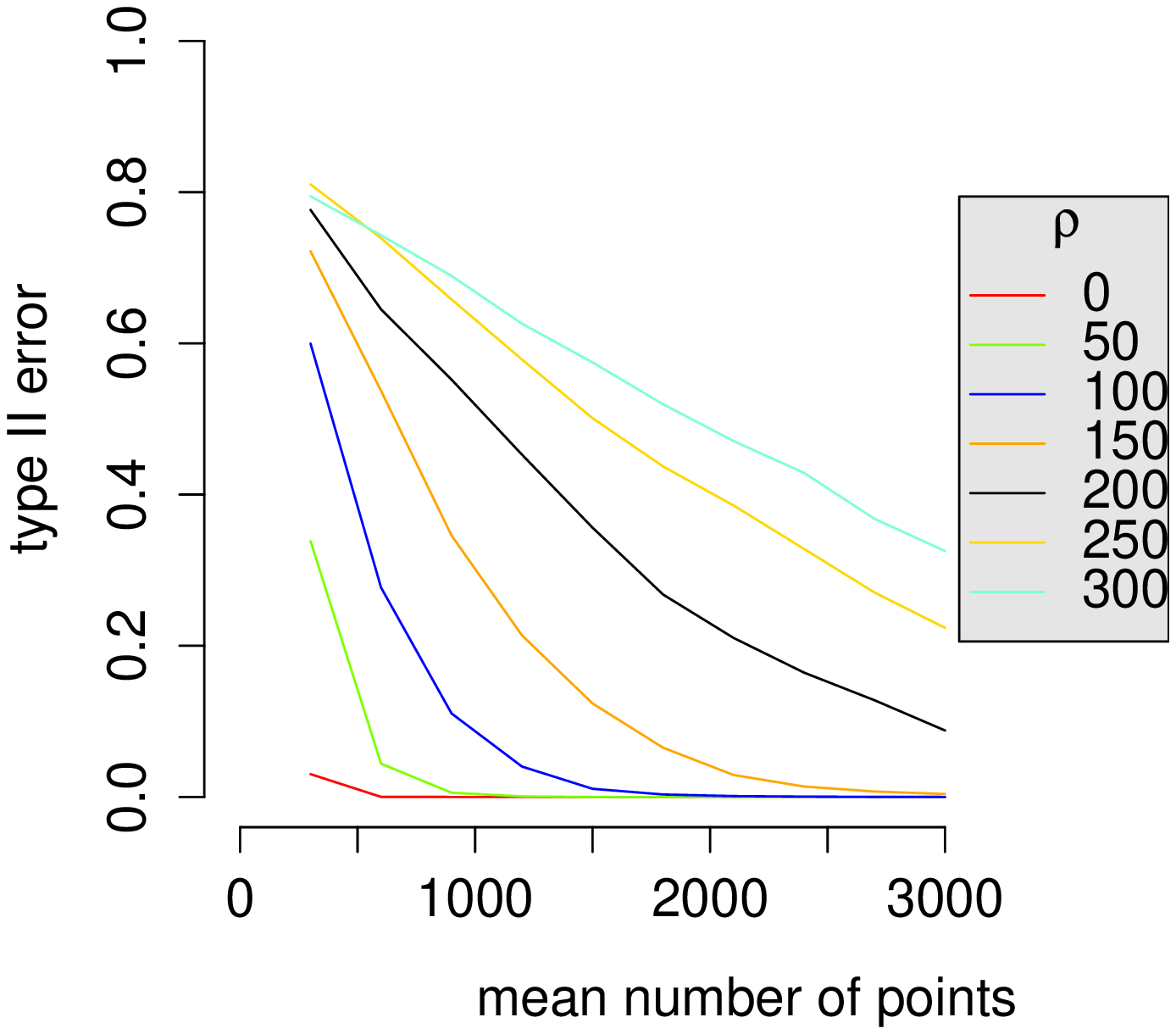}}\\

   \caption{Empirical errors of types I and II for the MGM test
  plotted against the mean number of points in the observation window ($\kappa_{12}=0.4$, $\alpha=0.5$). Different colors correspond to different values of the dependence parameter $\rho$.
}\label{2error_model_m_vs_pts}
\end{figure}

All simulation results are based on $10000$ model realizations per scenario. Type II errors were computed for
realizations where $\kappa_{12}\neq 0$, which means that the mark distribution was not uniform on ${\mathbb S}_+^1$, whereas $H_0:  \pom=U({\mathbb S}_+^1)$ hypothesized a uniform Palm mark distribution on ${\mathbb S}_+^1$ (corresponding to $\kappa_{12}=0$).\\
The performance of the MGM test is visualized in Tab.~\ref{table:2error_model_one-sided}. Empirical type I errors of the MGM test were close to the theoretical levels of significance
for $\alpha=0.025, 0.05,$ and $0.1$ with maximum deviations of around $0.015$. They were hardly affected by the 
dependence parameter $\rho$. Type II errors increased with $\rho$.
Under independent marking ($\rho=0$) as well as for $\rho=50$, error levels were close to $0$ for $\kappa_{12}\in\{ 0.2, 0.4, 0.8\}$. However, for an extreme range of dependence ($\rho=300$) even for a strong deviation of the data from a uniform Palm mark distribution ($\kappa_{12}=0.4$), rejection rates were only between $30$ and $40\%$.  For $\rho=300$ the range of dependence corresponds to $20\%$ of the sidelength of $W$.\\
Experiments with the TMD test revealed that the choice of the bandwidth parameter $c$ in (\ref{form:b_k_kappa}) is critical for test performance (Tab.~\ref{table:errors_kernel_m_vs_pts}). Whereas large values of $c$ result in small type I errors, they decrease the power of the test. On the other hand, small values for $c$ lead to superior power but increase type I errors (Tab. \ref{table:errors_kernel_m_vs_pts}). The empirical errors in Fig.~\ref{table:2error_kernel_one-sided} were computed for $c=50$ which yielded a reasonable compromise with respect to the two error types. In comparison to the MGM test the TMD test exhibits a higher sensitivity of empirical type I errors for varying values of $\kappa_{12}$, i.e., w.r.t. deviations from the uniform distribution on ${\mathbb S}_+^1$. Moreover, type II errors of the TMD test were up to $20\%$ higher than for the MGM test.\\
Tab.~\ref{table:errors_kernel_m_vs_pts} and Fig. \ref{2error_model_m_vs_pts} illustrate test performance
w.r.t. the mean number of points in $W$ and the dependence parameter $\rho$. The simulation experiments were conducted for $\alpha=0.05$. For power analysis, the tested data was simulated for $\kappa_{12}=0.4$, and thus the Palm mark distribution strongly differed from the uniform distribution on $\mathbb{S}^1$ of $H_0$.  
At a mean number of $3000$ observed points, $H_0$ was reliably rejected by the TMD test once $\rho \leq 150$ (for $c=50$). For $\rho \leq 100$ already $2000$ expected points were sufficient to reject $H_0$ for almost all realizations.
The MGM test required around $500$ points less than the TMD test in order to achieve comparable rejection rates (Fig.~\ref{2error_model_m_vs_pts}). \\
In summary, our simulation results indicate that the MGM test outperforms the TMD test especially with respect to power. This result is plausible since the additional information incorporated into $H_0$ by specification of a model covariance matrix can be expected to result in a more specific test.
It seems difficult to derive a general rule of thumb relating the required size of the observation window  to the dependence structure of the data and the intensity of the point pattern. However, Fig.~\ref{2error_model_m_vs_pts} and Tab.~\ref{table:errors_kernel_m_vs_pts} provide an idea on the practical requirements for asymptotic testing.

\bibliographystyle{ims}
\bibliography{literaturePaper}

\begin{table}

\begin{center}

\begin{tabular}[t!]{|p{1.5cm}|p{6.5cm}|p{6.5cm}|}

 \hline

&  Type I error & Type II error  \\

\hline
$\rho=0$ & \includegraphics[width=6.0cm, angle=0]{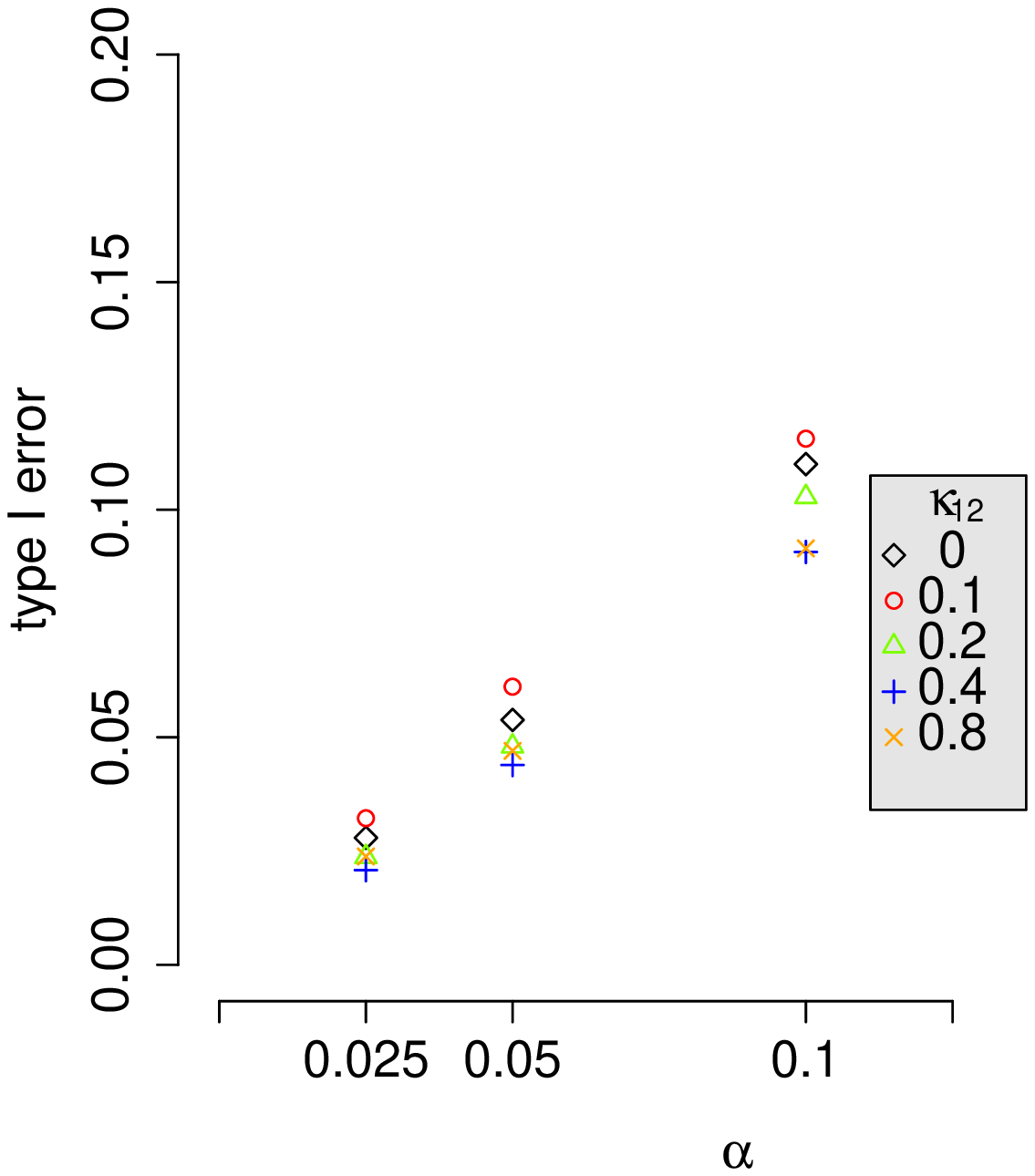}  &  \includegraphics[width=6.0cm, angle=0]{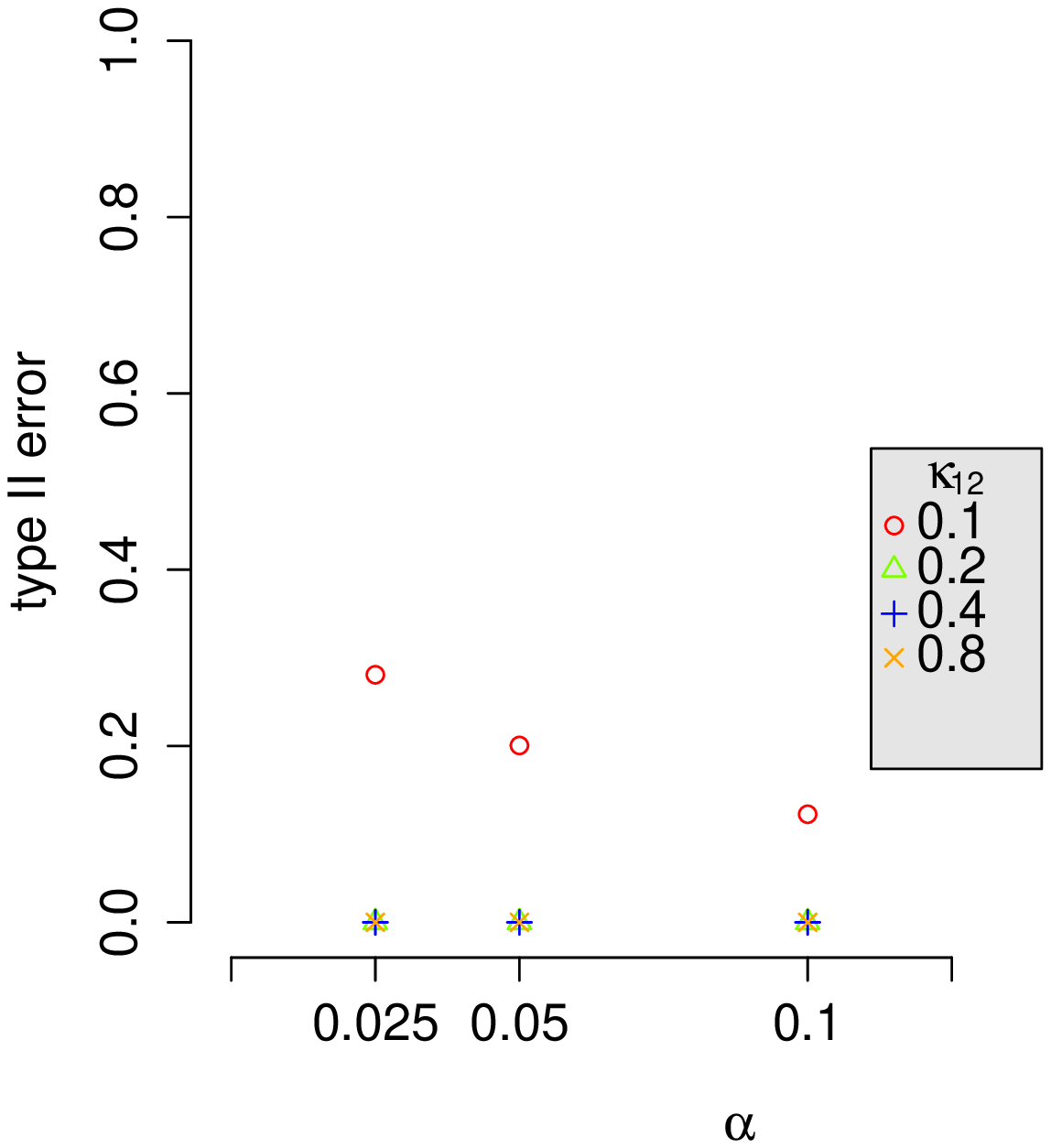}  \\
\hline

$\rho=50$ & \includegraphics[width=6.0cm, angle=0]{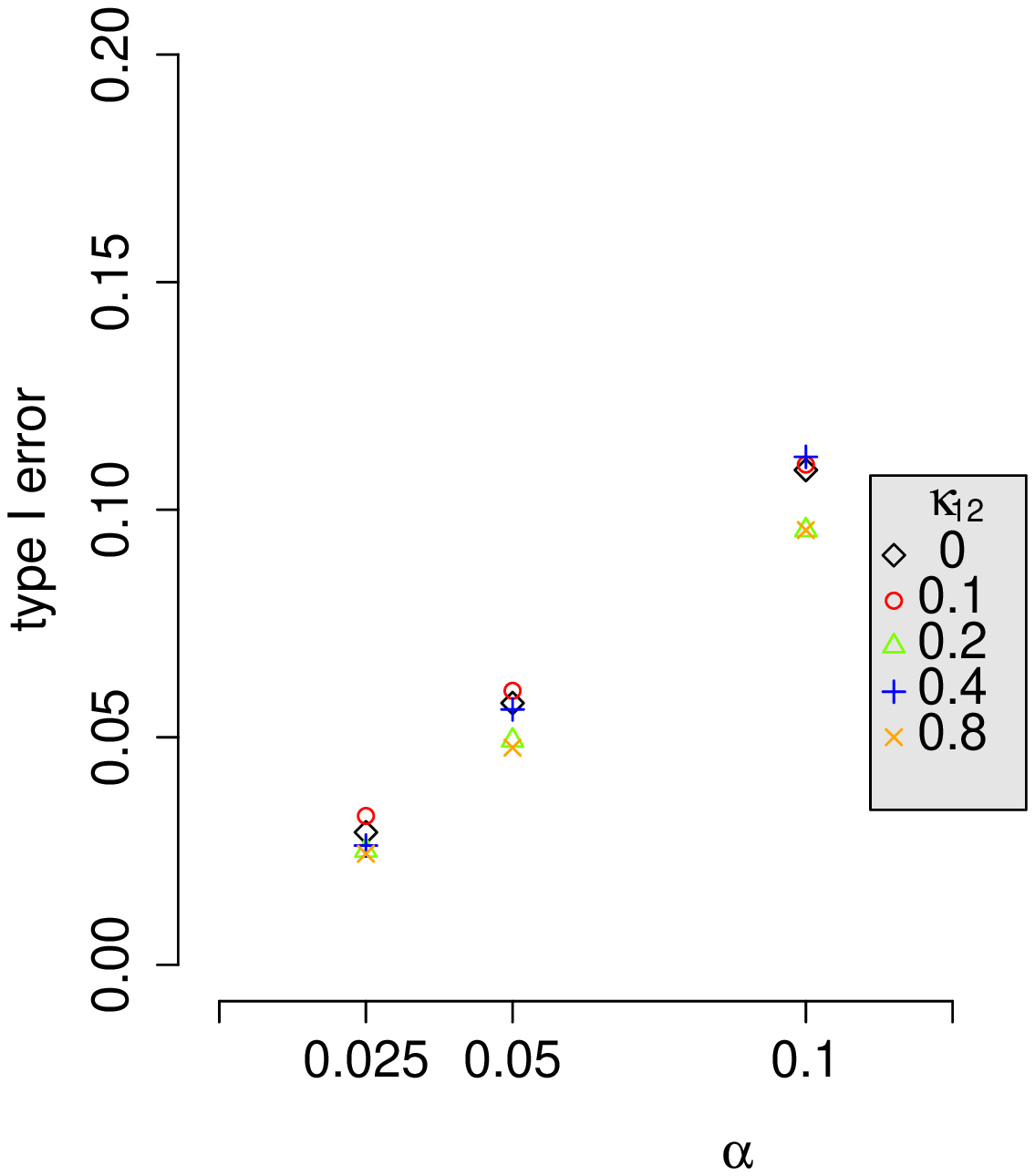}  &  \includegraphics[width=6.0cm, angle=0]{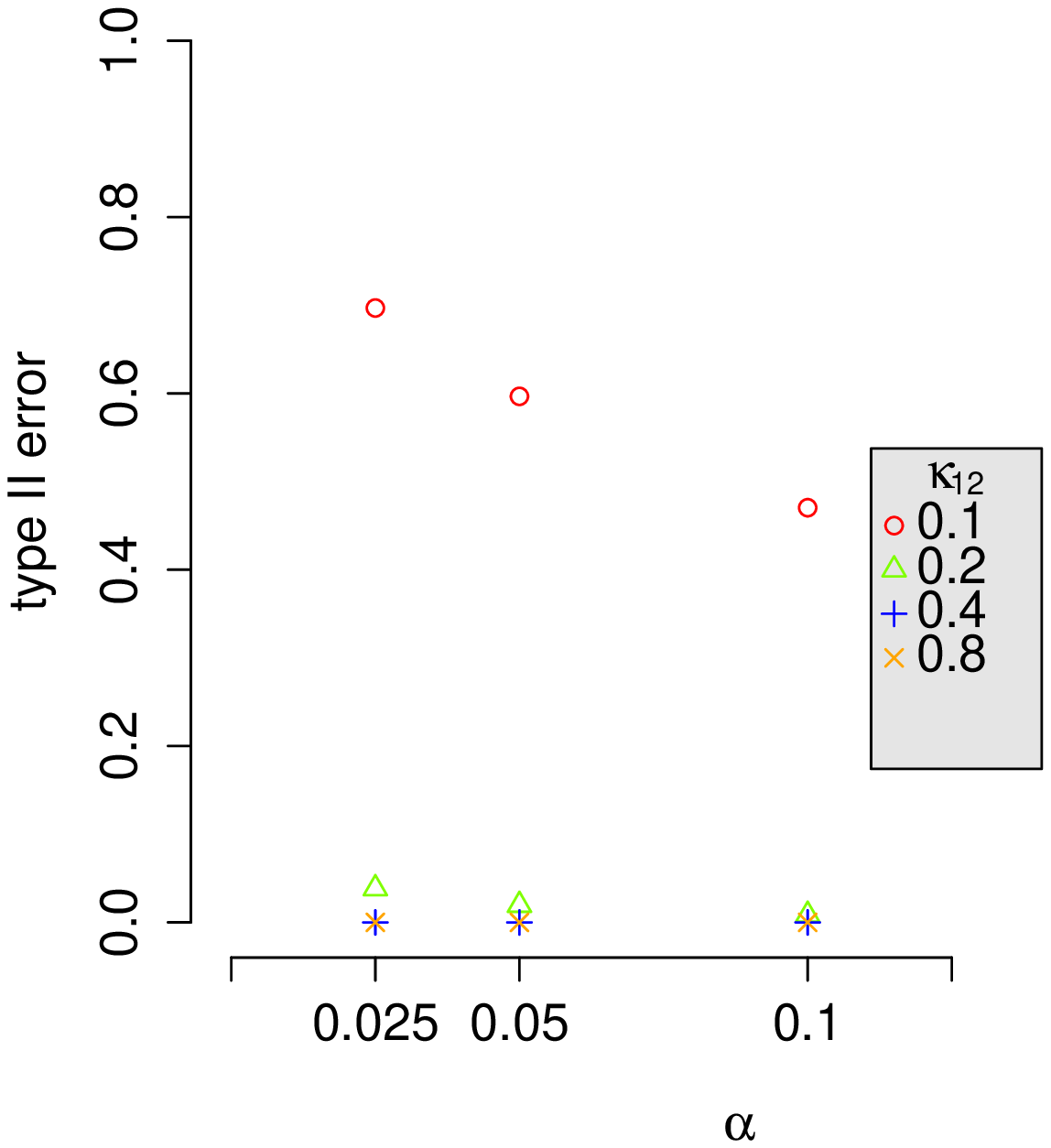}  \\
\hline

$\rho=300$ & \includegraphics[width=6.0cm, angle=0]{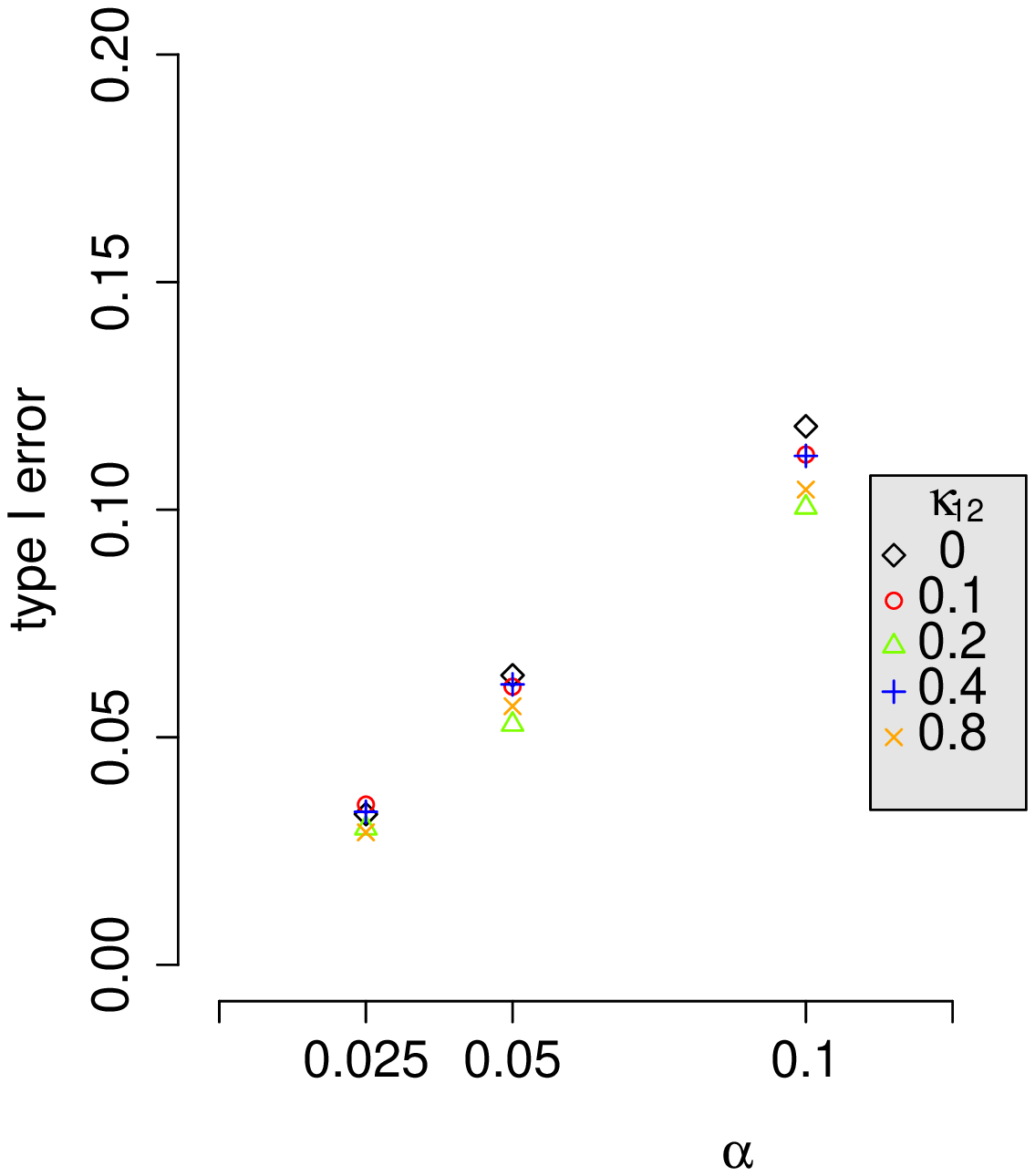}  &  \includegraphics[width=6.0cm, angle=0]{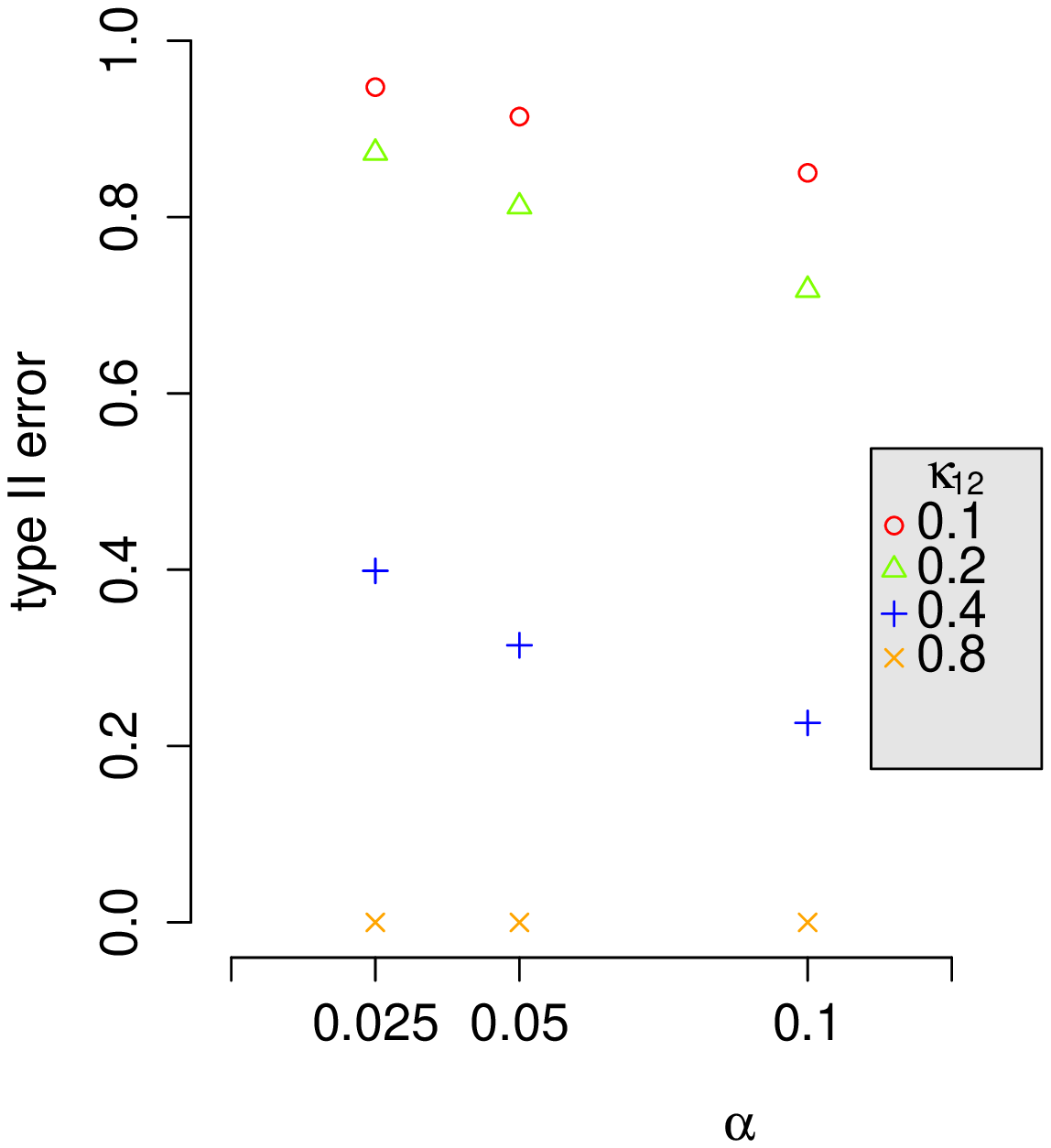}  \\
\hline


\end{tabular}

\caption{Empirical errors of types I and II for the MGM test.
}\label{table:2error_model_one-sided}

\end{center}
\end{table}

%

\begin{table}

\begin{center}

\begin{tabular}[t!]{|p{1.5cm}|p{6.5cm}|p{6.5cm}|}

 \hline

&  Type I error & Type II error \\

\hline
$\rho=0$ & \includegraphics[width=6.0cm, angle=0]{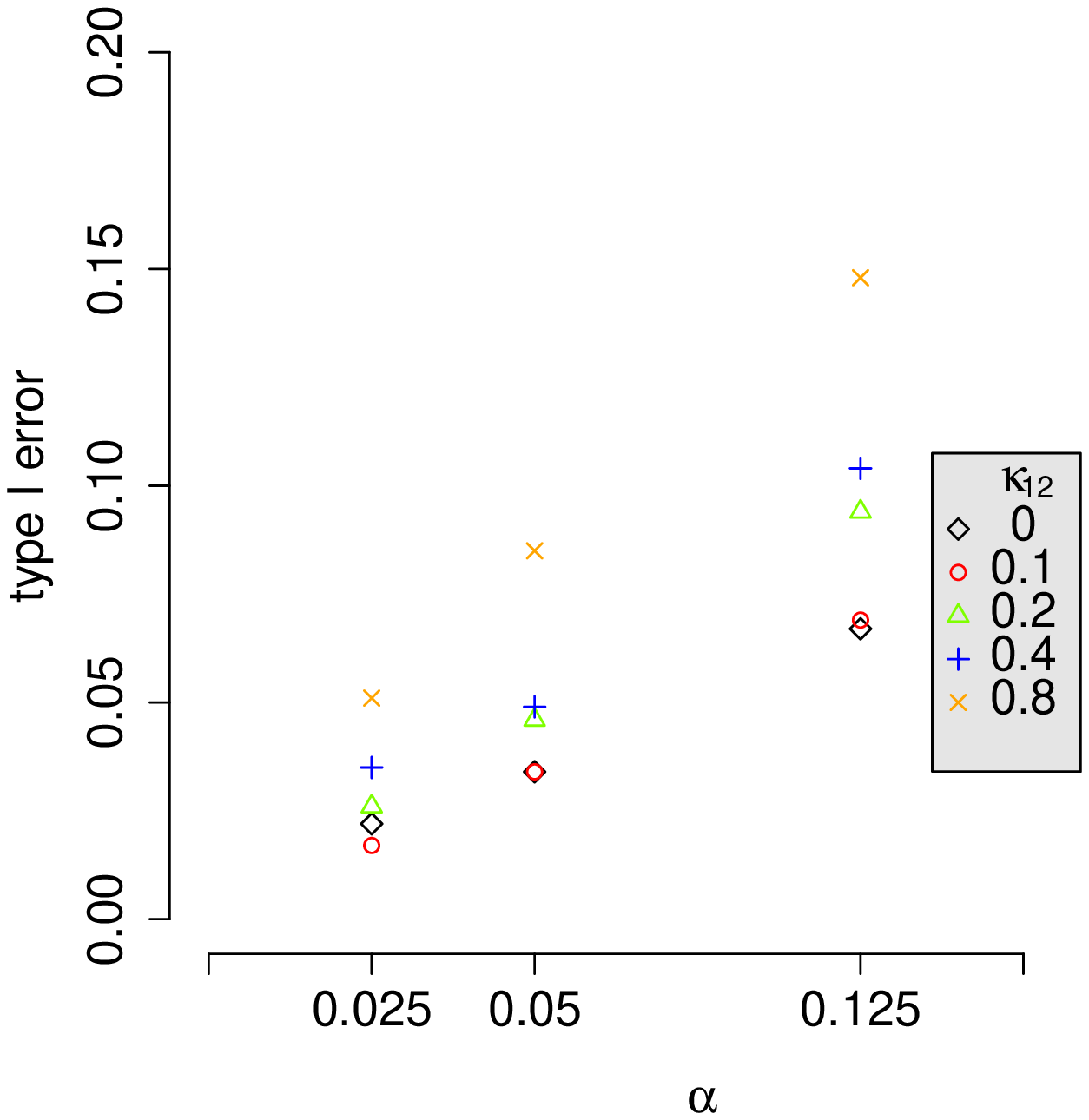}  &  \includegraphics[width=6.0cm, angle=0]{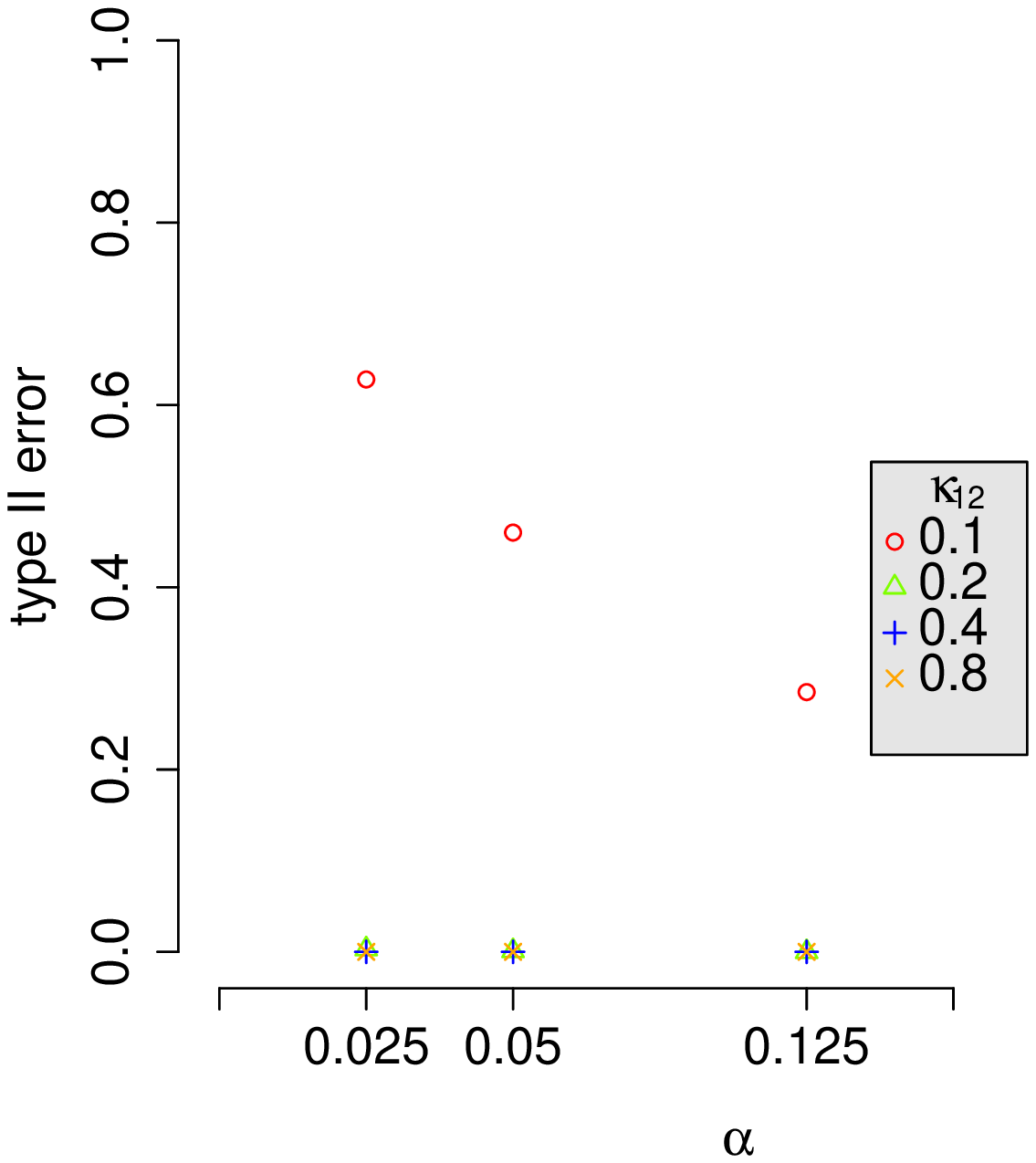}  \\
\hline

$\rho=50$ & \includegraphics[width=6.0cm, angle=0]{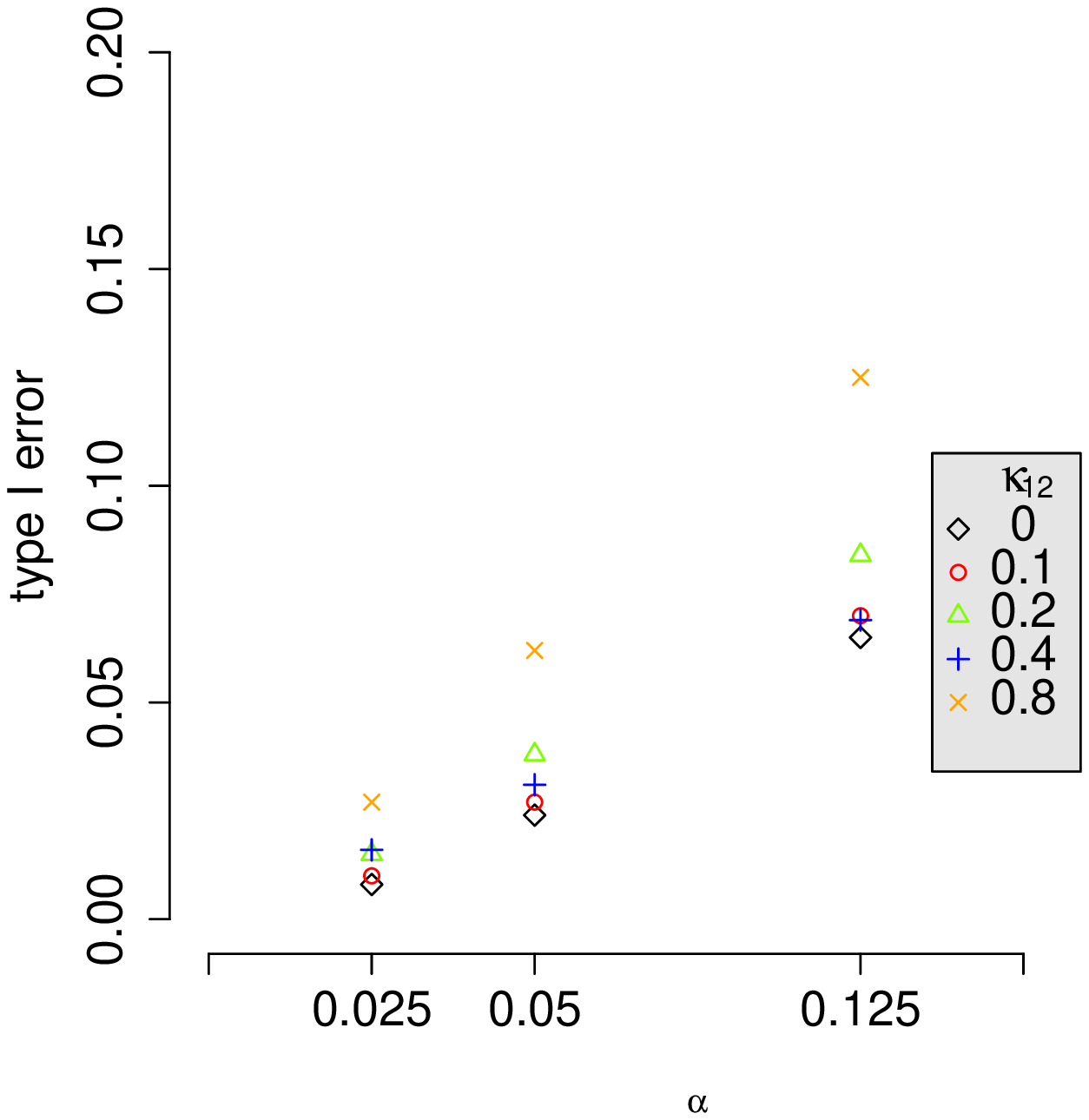}  &  \includegraphics[width=6.0cm, angle=0]{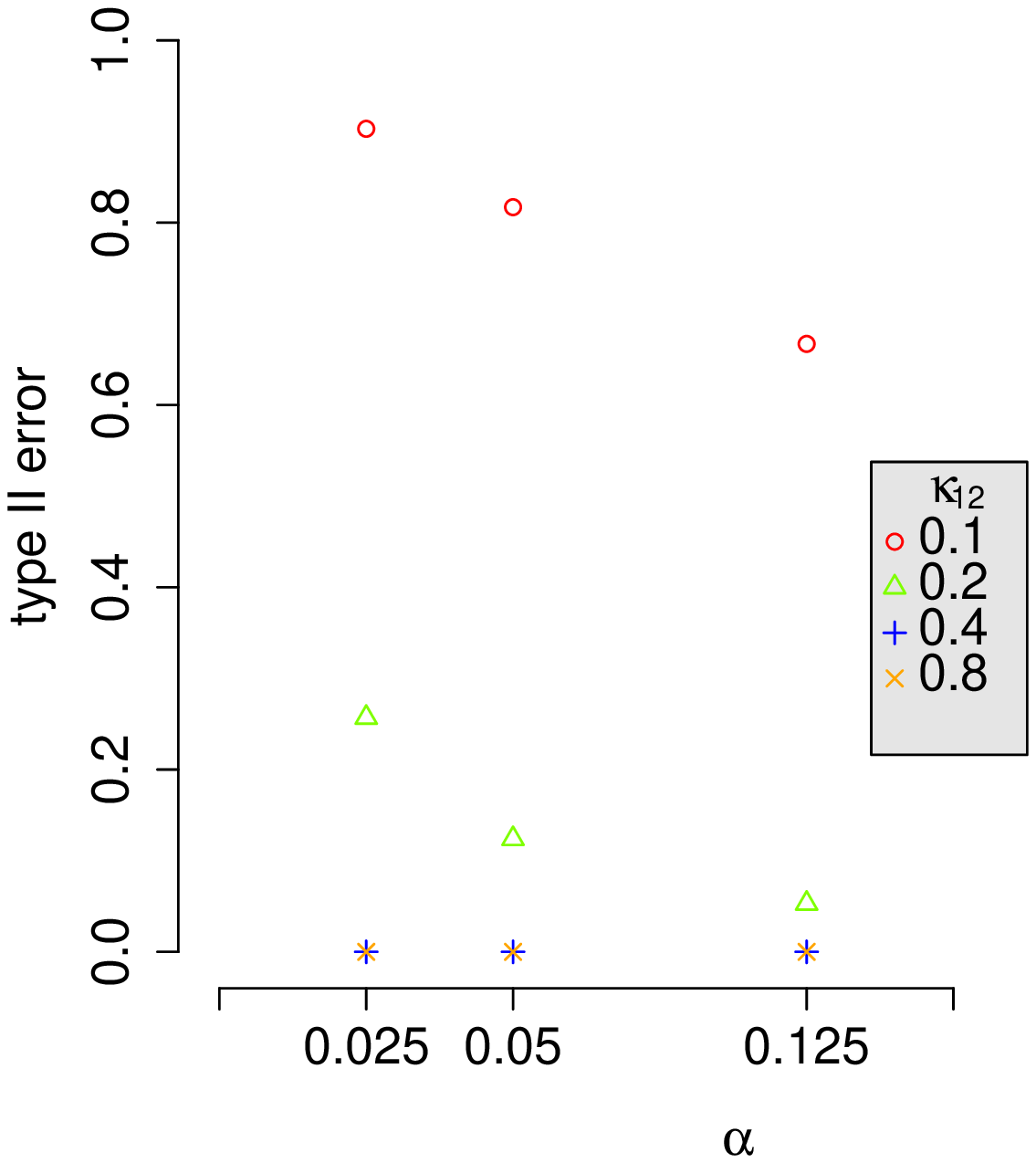}  \\
\hline

$\rho=300$ & \includegraphics[width=6.0cm, angle=0]{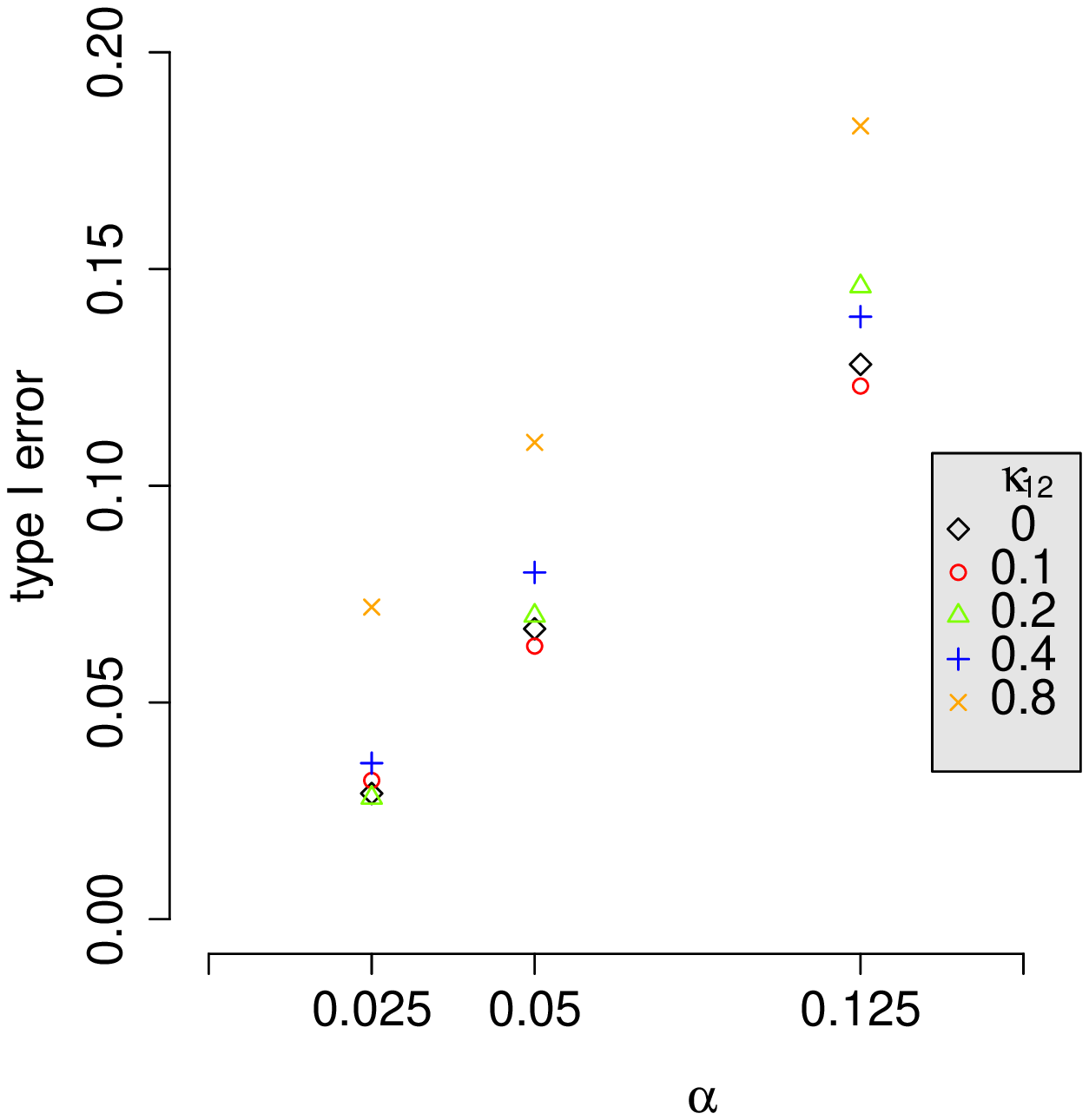}  &  \includegraphics[width=6.0cm, angle=0]{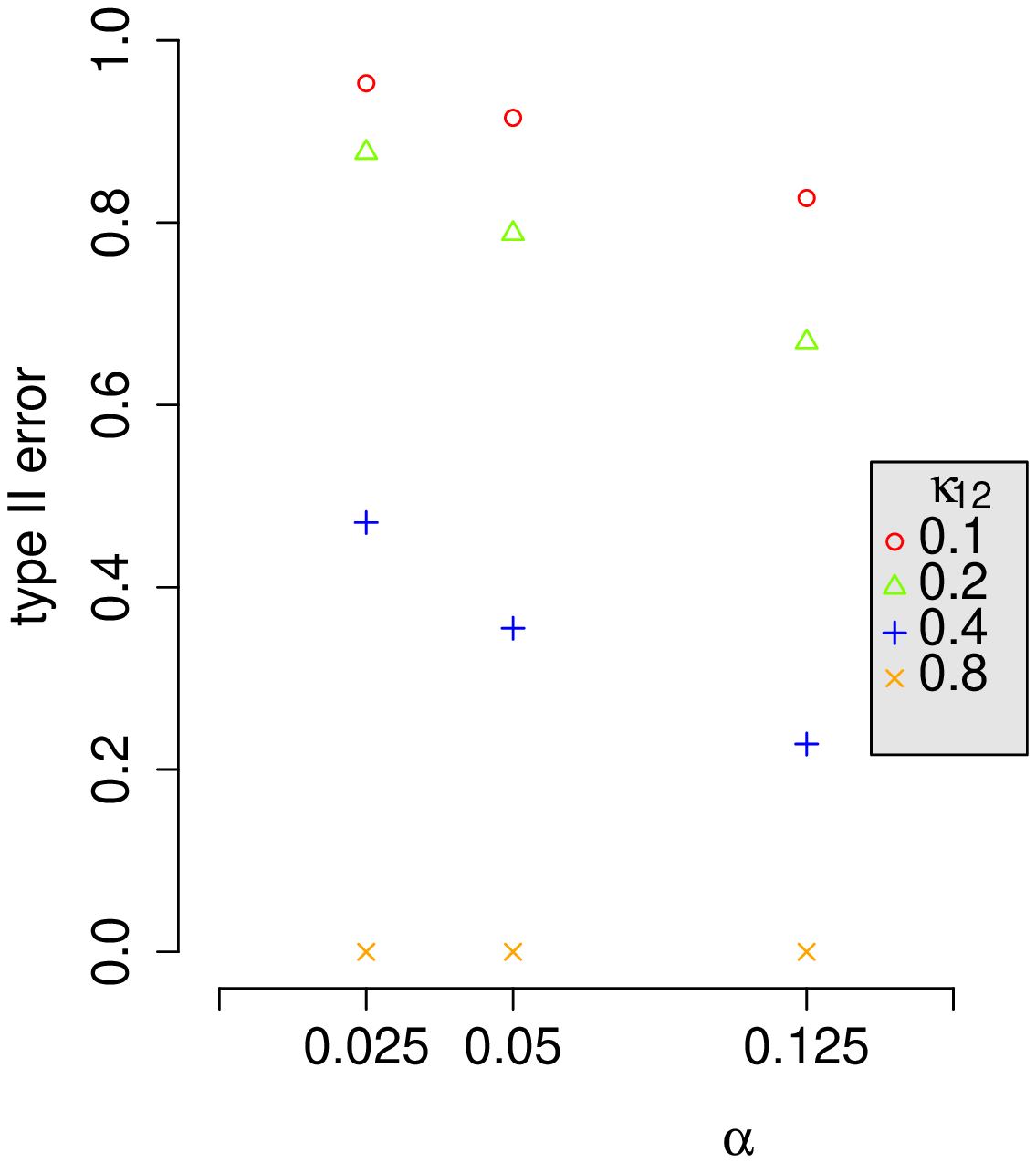}  \\
\hline


\end{tabular}

\caption{ Empirical errors of types I and II for the TMD test ({$c=50$}).
}\label{table:2error_kernel_one-sided}

\end{center}
\end{table}


\begin{table}

\begin{center}

\begin{tabular}[t!]{|p{1.5cm}|p{4.5cm}|p{4.5cm}|}

 \hline

&  Type I error & Type II error \\

\hline
$c=20$ & \includegraphics[width=4.0cm, angle=-90]{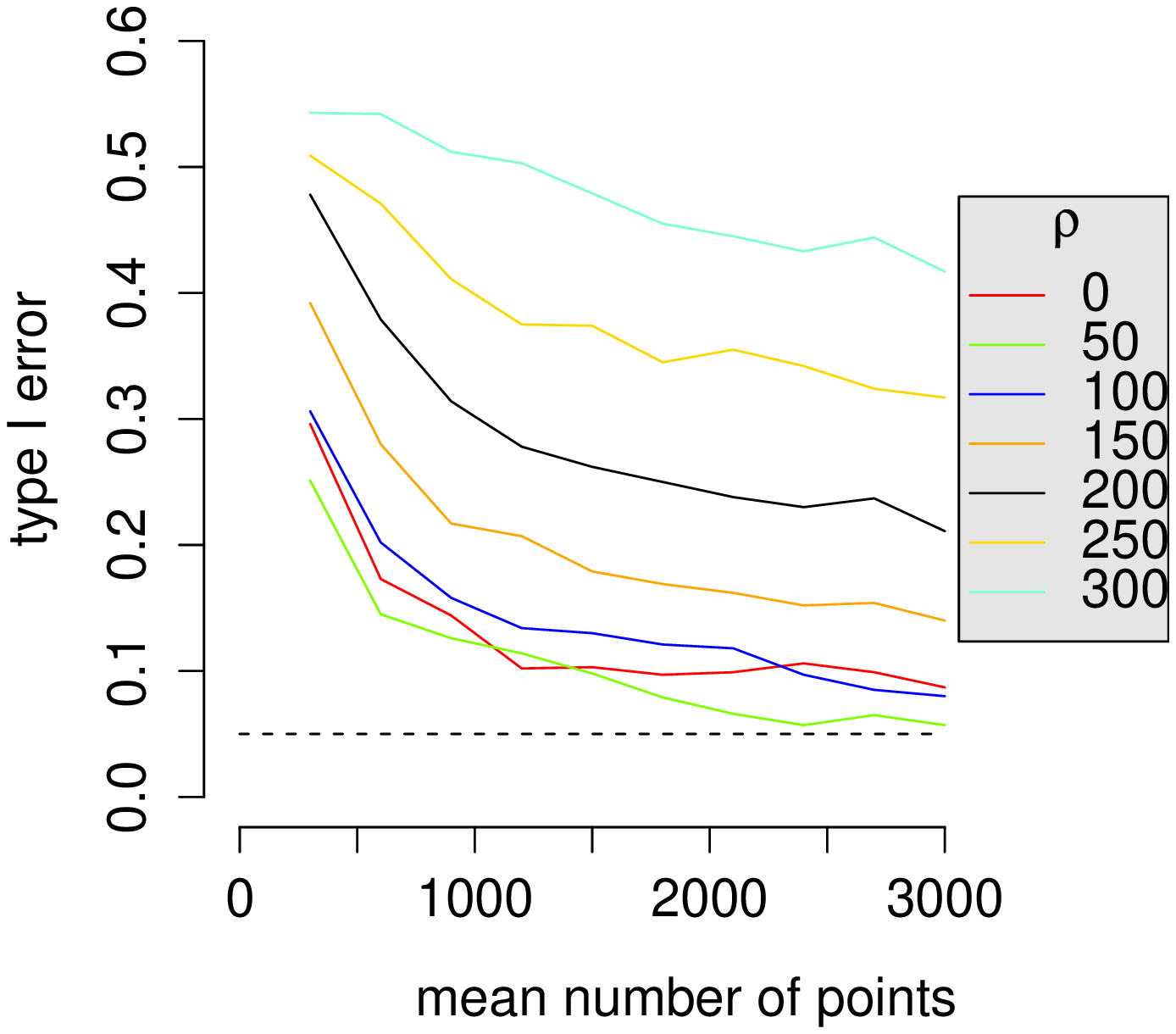}  &  \includegraphics[width=4.0cm, angle=-90]{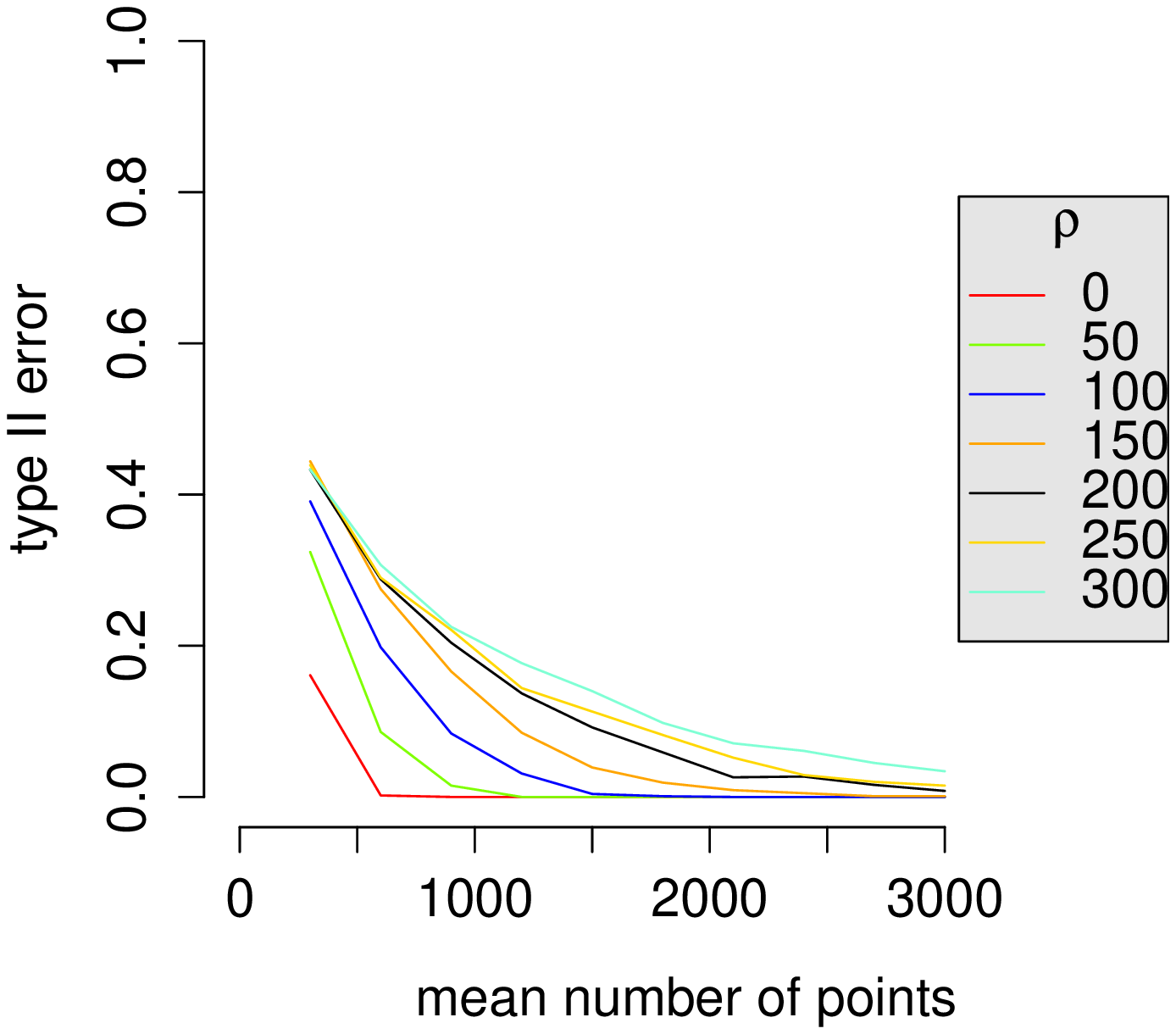}  \\
\hline

$c=30$ & \includegraphics[width=4.0cm, angle=-90]{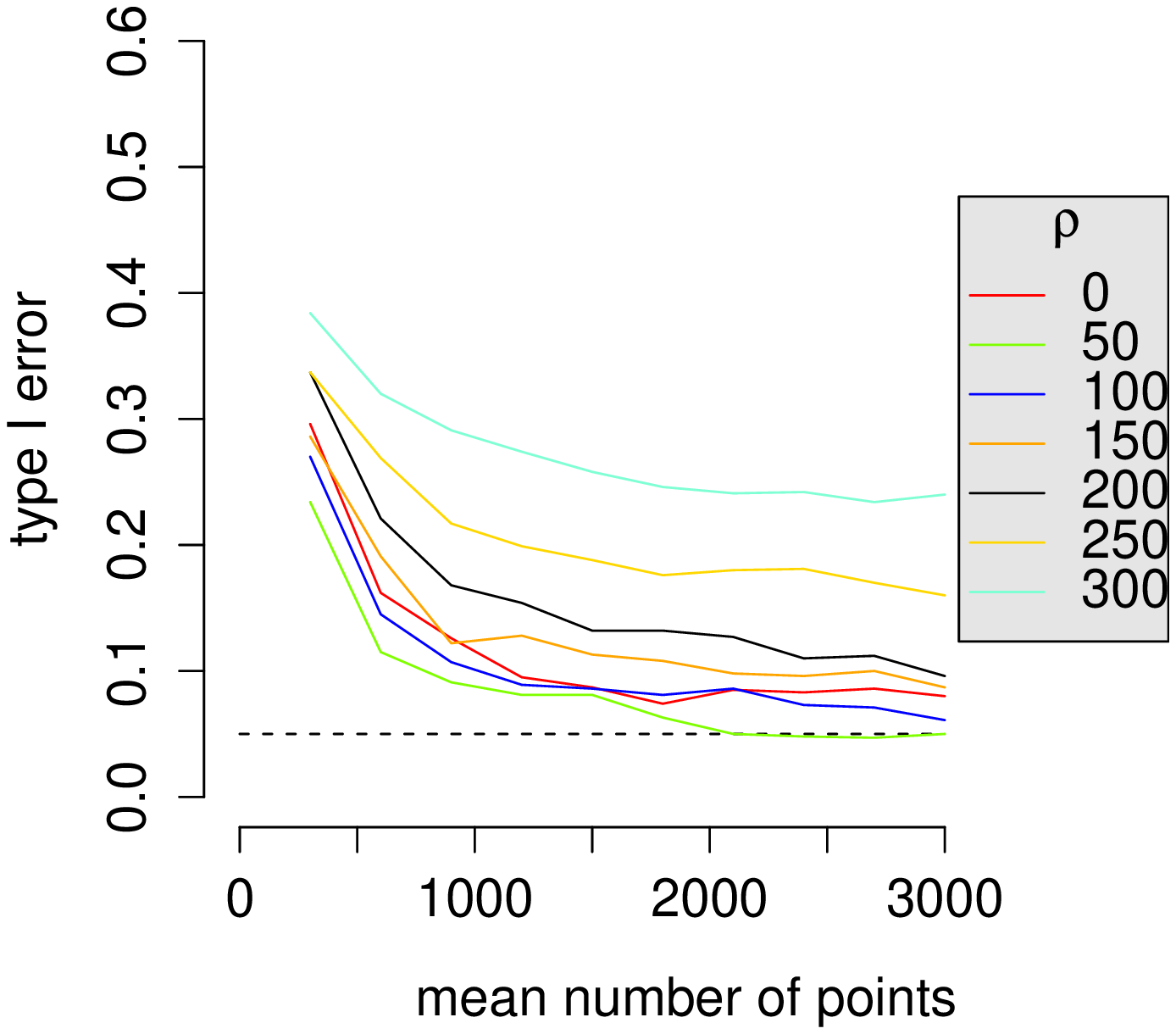}  
&  
\includegraphics[width=4.0cm, angle=-90]{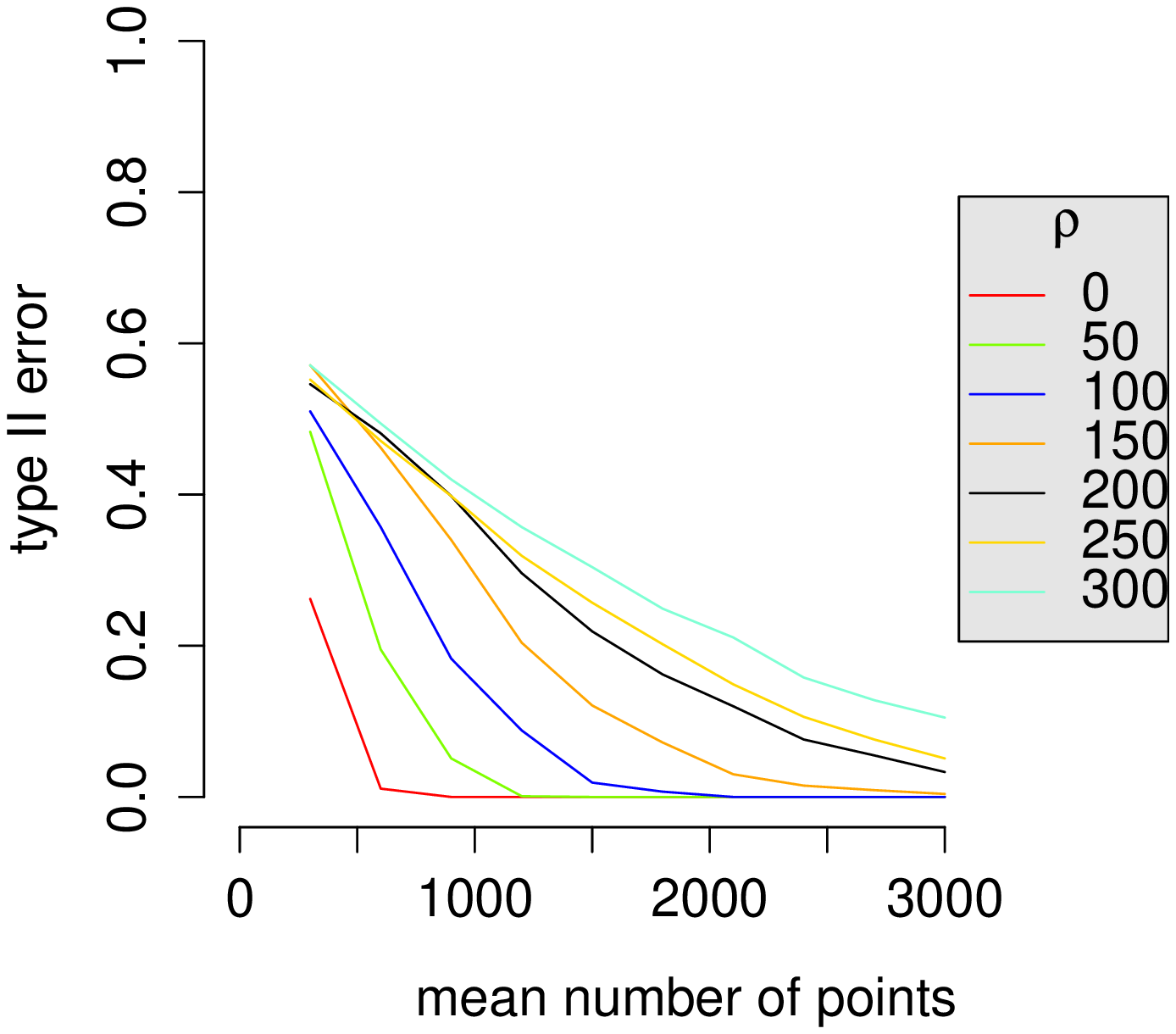}  \\
\hline

$c=40$ & \includegraphics[width=4.0cm, angle=-90]{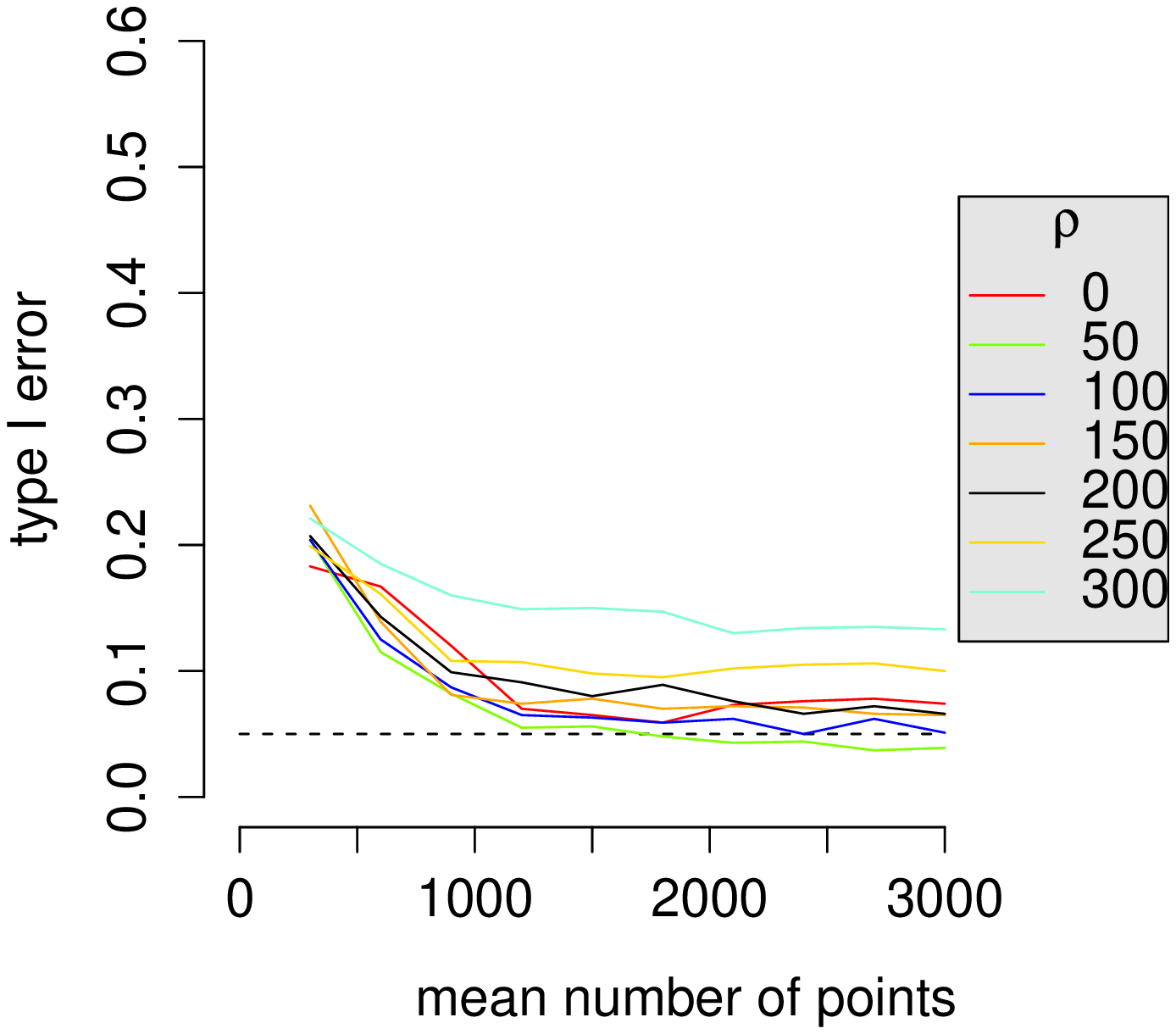}  
&  
\includegraphics[width=4.0cm, angle=-90]{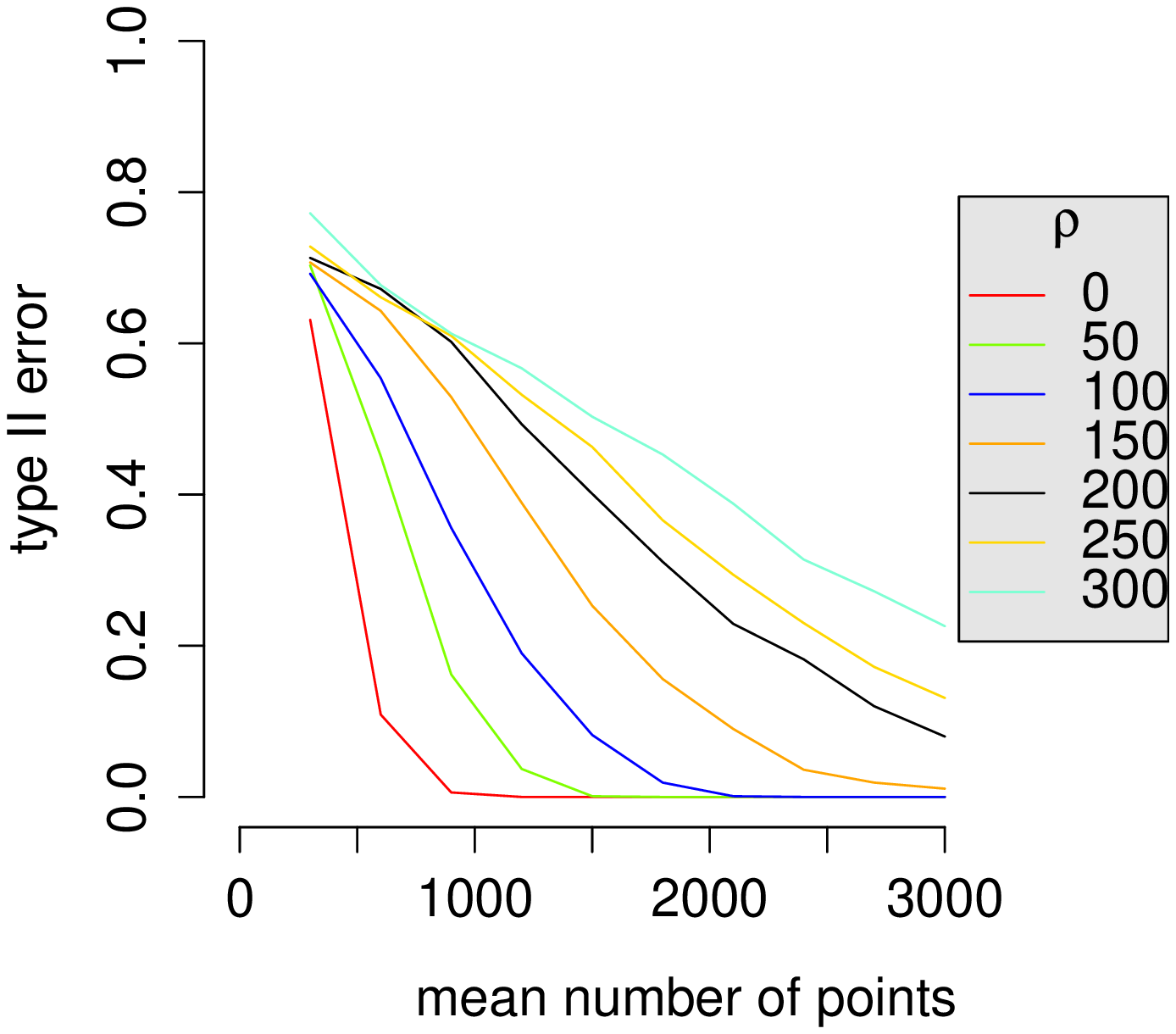}  \\
\hline

$c=50$ & \includegraphics[width=4.0cm, angle=-90]{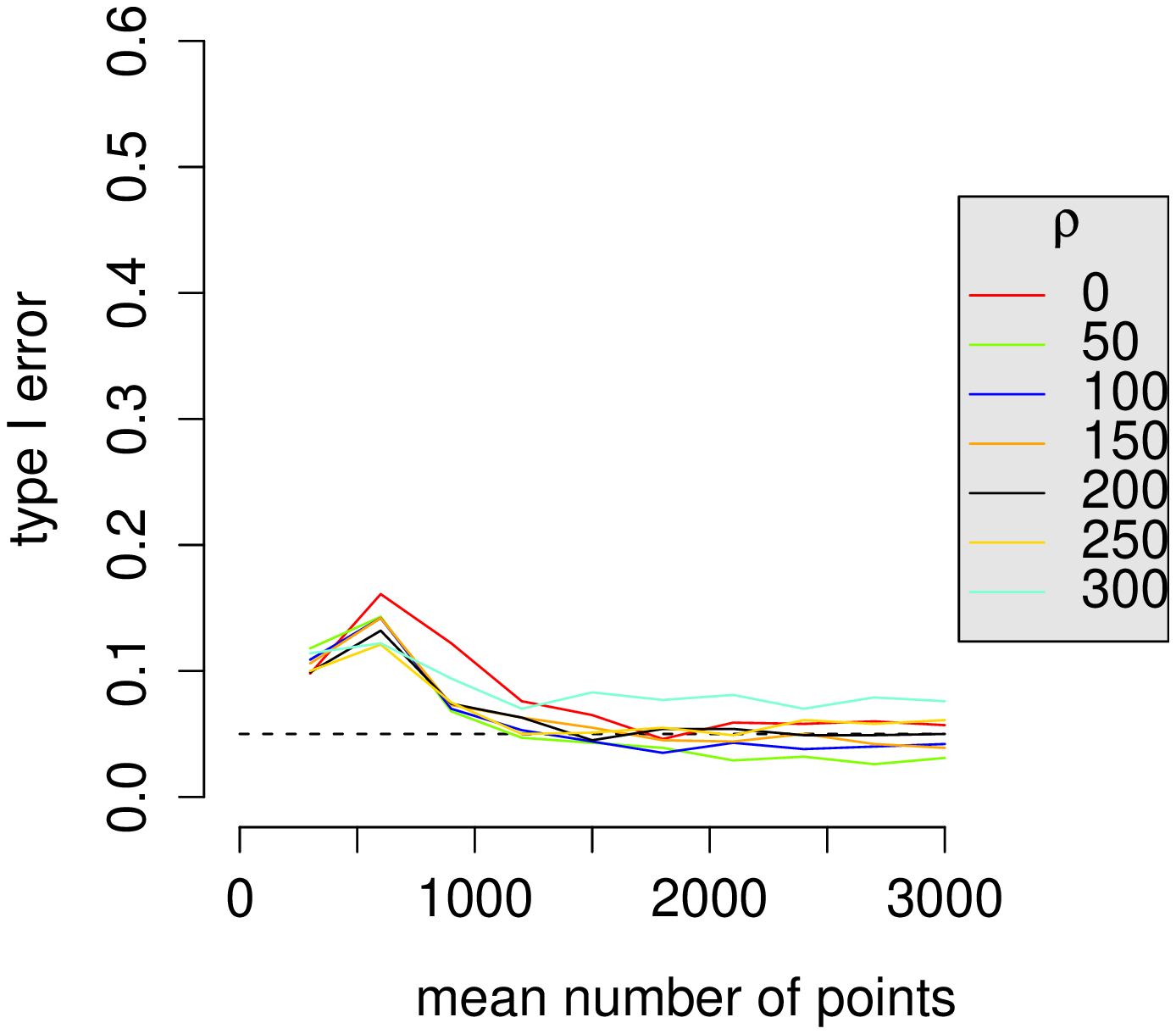}  &  \includegraphics[width=4.0cm, angle=-90]{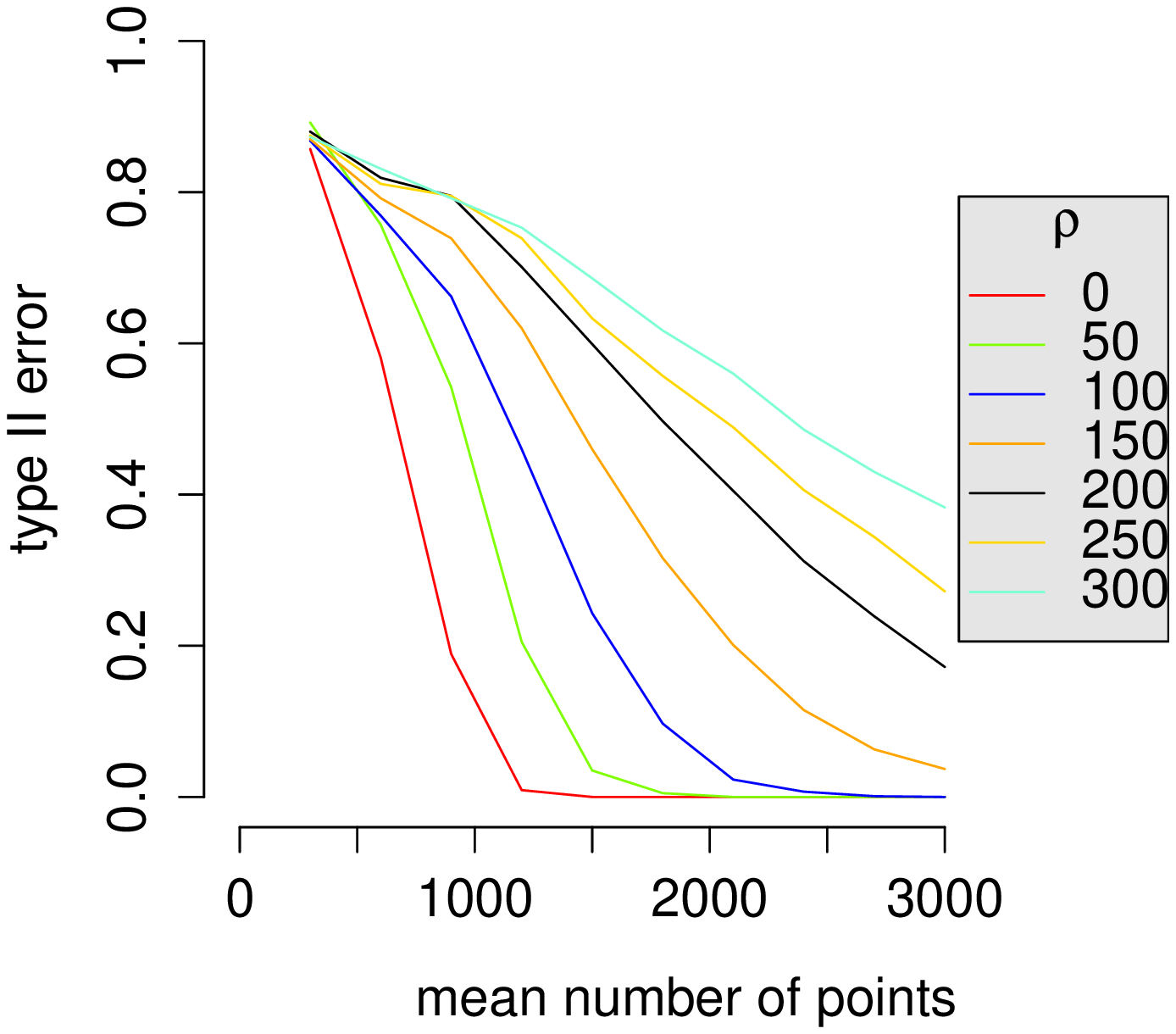}  \\
\hline

$c=60$ & \includegraphics[width=4.0cm, angle=-90]{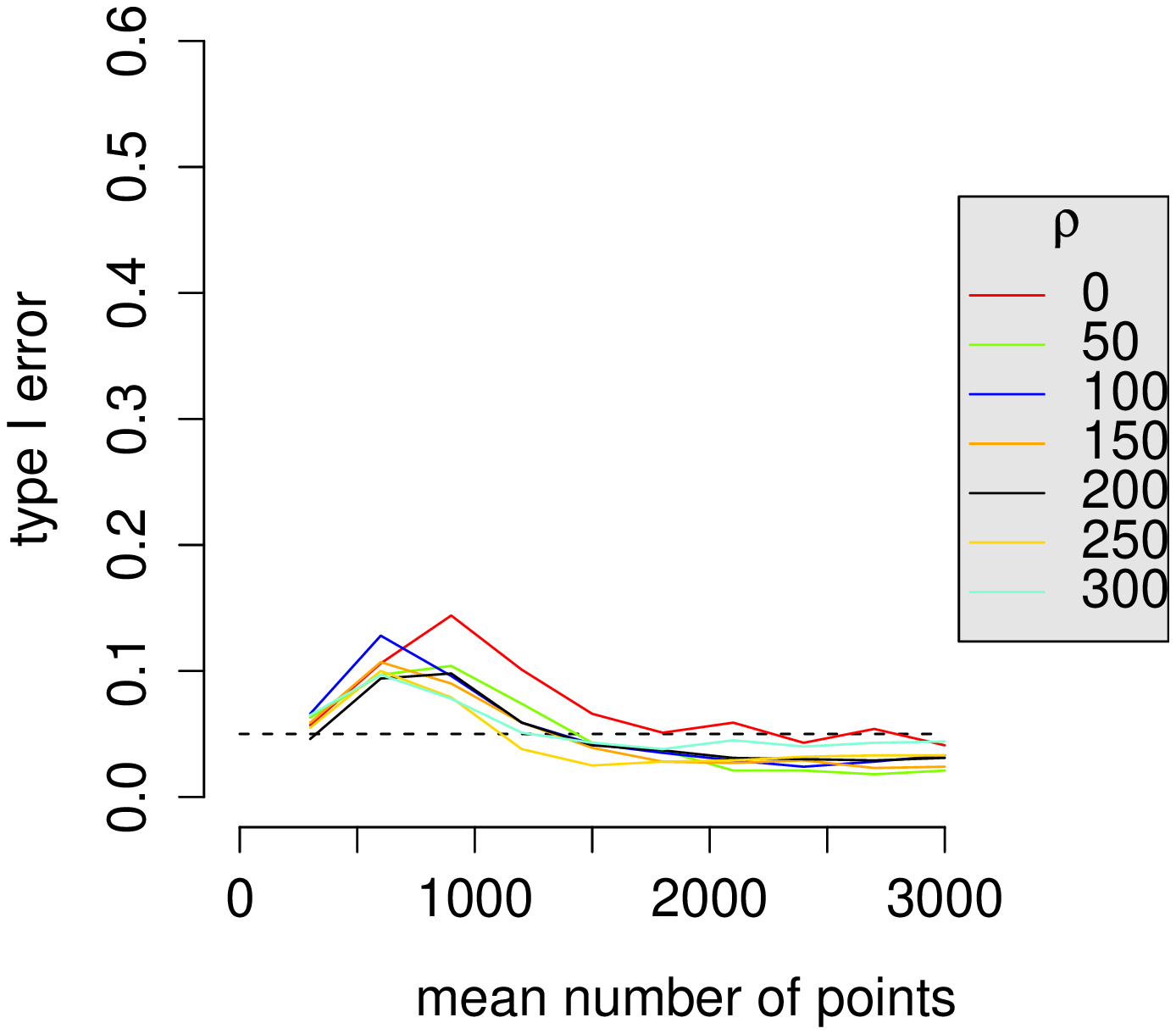}  &  \includegraphics[width=4.0cm, angle=-90]{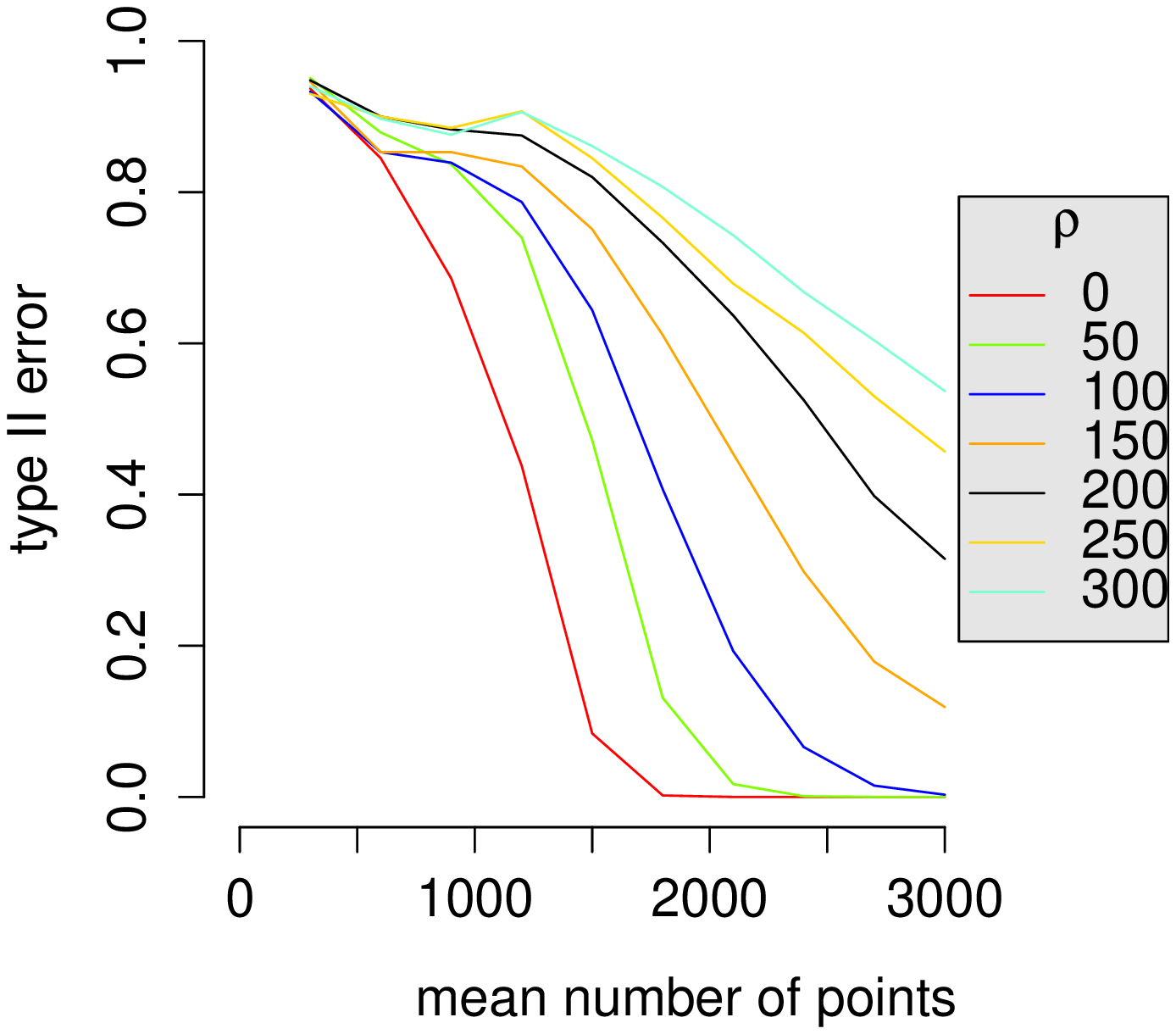}  \\
\hline


\end{tabular}

\caption{Empirical errors of types I and II for the TMD test plotted against the mean number of points in the observation window ($\kappa_{12}=0.4$, $\ell=8$, and $\alpha=0.05$). Different colors correspond to different values of the dependence parameter $\rho$.}\label{table:errors_kernel_m_vs_pts}
\end{center}
\end{table}

\end{document}